\documentclass[11pt,a4paper,twoside,leqno]{article}
\usepackage{graphicx}
\usepackage{geometry}
\usepackage{amssymb}
\usepackage{amsmath}
\usepackage[hyperfootnotes=false]{hyperref}
\usepackage{subfig}
\usepackage{caption}
\usepackage{algpseudocode}
\usepackage{algorithm}

\newcommand{\pa}[2]{\frac{\partial #1}{\partial #2}}
\title{GPU Accelerated Discontinuous Galerkin Methods for Shallow Water Equations}
\def\usFootnote{Department of Computational and Applied Mathematics, Rice University, 6100 Main Street - MS 134, Houston, TX 77005, USA}
\def\email{Corresponding author. \textit{Email addresses:
Rajesh.Gandham@rice.edu (R. Gandham), David.S.Medina@rice.edu (D. Medina),
timwar@rice.edu (T. Warburton)} }

\author{Rajesh Gandham \footnote{\email} \textsuperscript{,}\footnote{\usFootnote},\
David Medina\footnotemark[\value{footnote}],\
  Timothy Warburton\footnotemark[\value{footnote}]}
\date{}

\everymath{\displaystyle}
\usepackage[normalem]{ulem}
\numberwithin{equation}{section}
\numberwithin{figure}{section}

\numberwithin{algorithm}{section}
\numberwithin{table}{section}
\begin{document}

\maketitle

\begin{abstract}
We discuss the development, verification, and performance of a GPU accelerated discontinuous Galerkin method  for the solutions of two dimensional nonlinear shallow water equations. The shallow water equations are hyperbolic partial differential equations and are widely used in the simulation of tsunami wave propagations. 
Our algorithms are tailored to take advantage of the single instruction multiple data (SIMD)  architecture of graphic processing units. The time integration is accelerated by  local time stepping based on a multi-rate Adams-Bashforth  scheme. 
A total variational bounded  limiter is adopted for nonlinear stability of the numerical scheme. This limiter is coupled with a mass and momentum conserving positivity preserving limiter for the special treatment of a dry or partially wet element in the triangulation.
Accuracy, robustness and performance are demonstrated with the aid of test cases. We compare the performance of the kernels expressed in a portable threading language OCCA, when cross compiled with OpenCL, CUDA, and OpenMP at runtime.
\end{abstract}



\section{Introduction}
\label{sec:intro}

The shallow water equations (SWE) are of great interest in the modeling of tsunamis, storm surges and tidal waves. They are the simplest nonlinear models for water wave propagation. The shallow water assumptions simplify  three dimensional wave propagation to two dimensional hyperbolic partial differential equations, reducing the complexity of the model. The reduced complexity makes the shallow water equations attractive for tsunami modeling. These equations are valid for long waves but may represent the wave propagation of short waves or dispersive waves poorly. However,   these simplified equations provide satisfactory solutions of tsunami wave propagation \cite{leveque2011tsunami} over long distances.

The shallow water equations are two dimensional hyperbolic PDEs with velocity and fluid height as unknown quantities. 
These equations are complicated by the presence of largely varying length scales, varying bathymetry, and nonlinear effects near the shore. Stable, accurate and efficient algorithms are of great interest for these applications.

There is an extensive literature for finite difference \cite{casulli1990semi, merrill1997finite}, finite volume \cite{leveque2011tsunami, berger2011geoclaw} and finite element \cite{navon1988review} methods for shallow water equations. Recently, there has been growing interest in using discontinuous Galerkin methods (DG) for  solutions of the shallow water equations \cite{aizinger2002discontinuous,kubatko2006hp, eskilsson2004triangular, giraldo2008high}. DG methods are locally mass conservative like finite volume methods and can achieve high order accuracy on unstructured meshes like finite element methods. This allows flexibility in handling irregular boundaries with out compromising accuracy for problems with sufficiently smooth solutions. DG methods can achieve $\mathcal{O}(H ^{N+1/2})$ accuracy with a piecewise degree $N^{\text{}}$  polynomial approximation \cite{johnson1986analysis}. Where, $H$ is the largest length scale in the mesh. 

In DG formulations, elements are coupled using weak penalty terms, resulting in localized memory access. Furthermore, a high order polynomial representation of the solution in each element results in high arithmetic intensity per degree of freedom. Both of these features are well suited for the GPU hardware architecture \cite{klockner2009nodal, godel2010gpu, klockner2010high}. This motivated us to adopt GPUs along with a nodal DG discretization for large scale tsunami simulations. Furthermore, to alleviate the need to write kernels for thread models like OpenMP, CUDA, and OpenCL separately, we developed OCCA: A unified approach to multi-threading languages (unpublished). Kernels written in OCCA are cross compiled with any of these thread models at runtime.

This paper is organized as follows: In Sections \ref{sec:governing} and \ref{sec:discretization}, we outline the governing equations and nodal discontinuous Galerkin discretization. In Section \ref{sec:mrab}, we describe local time stepping using multi-rate Adams-Bashforth time integration. In Section \ref{sec:PP}, we explain the stabilization of numerical scheme using positivity preserving, wetting drying treatment and modified totally variational bounded (TVB) limiter. Several tests for verification of accuracy and robustness are presented in Section \ref{sec:tests}. We discuss GPU kernels and their performance in Section \ref{sec:gpu}.


\section{Governing Equations}
\label{sec:governing}

\begin{figure}[h!]
\begin{center}
\includegraphics[trim=2cm 18cm 1cm 0cm,clip=true,width=0.6\textwidth]{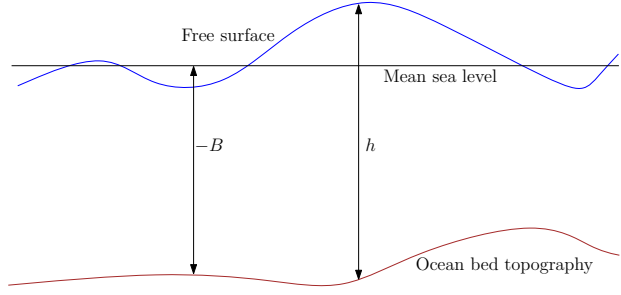} \hspace{1cm}
\caption{Diagram of notations.}
\label{fig:notation}
\end{center}
\end{figure}

The shallow water equations are depth averaged incompressible Navier-Stokes equations, and are given in conservative form by \cite{leveque2011tsunami},
\begin{eqnarray}
\label{eq:swe}
\pa{h}{t} + \pa{(hu)}{x} + \pa{(hv)}{y} &=& 0, \nonumber \\
\pa{}{t}(hu) + \pa{}{x}\left( hu^2 + \frac{1}{2}gh^2 \right) + \pa{}{y} (huv)&=& - gh\pa{B}{x} , \nonumber \\
\pa{}{t}(hv) + \pa{}{x}(huv) + \pa{}{y}\left( hv^2 + \frac{1}{2}gh^2\right) &=&  - gh\pa{B}{y} ,
\end{eqnarray}
where, $h$, $u$, and $v$ are water depth, depth averaged velocity components in longitudinal and latitudinal directions. $B$ is bathymetry and $g$ is the acceleration due to gravity (see Fig. (\ref{fig:notation})  for the notation). 

 In simplified form, these equations are represented as,
\begin{equation}
\label{eq:swe_conservative}
\pa{Q}{t} + \pa{F}{x} + \pa{G}{y} = S,
\end{equation}
\noindent The state vector $Q$, nonlinear flux vectors $F$, $G$ and the source vector $S$ are given by,
\begin{equation}
Q = \begin{pmatrix} h  \vspace{0.2cm}\\  hu \vspace{0.2cm}\\ hv\end{pmatrix}, \quad
F = \begin{pmatrix} hu \vspace{0.2cm}\\  hu^2 + \frac{gh^2}{2} \vspace{0.2cm}\\ huv \end{pmatrix}, \quad
G = \begin{pmatrix} hv \vspace{0.2cm}\\  huv \vspace{0.2cm}\\ hv^2 + \frac{gh^2}{2} \end{pmatrix}, \quad
S = \begin{pmatrix} 0  \vspace{0.2cm}\\  -gh \pa{B}{x} \vspace{0.2cm} \\ -gh \pa{B}{y} \end{pmatrix}.
\end{equation}
We use a high order discontinuous Galerkin method to obtain the solution of the Eq. (\ref{eq:swe_conservative}). We explain the method in the next section. 

\section{Discretization}
\label{sec:discretization}
We assume the domain   $\Omega \subset \mathbb{R}^2$, is partitioned in to a set of non-overlapping, conforming triangles $\{ \Omega = \displaystyle\cup_{k} D^{k} \}$. We approximate the solution $Q$ by $Q_H$, the components of which belong to the space of discontinuous piecewise polynomial of a given degree $N$ in each element ($P^{N}(D^{k})$). The PDEs in the Eq. (\ref{eq:swe_conservative}) are expressed in weak form. For each element we find $Q_H \in (P^{N}(D^{k}))^{3}$ such that,
\begin{equation}
\left( \pa{Q_{H}}{t} , \phi \right)_{D^{k}} + \left( \pa{F}{x} , \phi \right)_{D^{k}} + \left( \pa{G}{y} , \phi \right)_{D^{k}} - (S, \phi)_{D^{k}} = 0 \qquad \forall \phi \in P^{N}(D^{k}),
\end{equation}

  Integrating by parts and replacing the multi-valued fluxes  on the boundary of each element with stable numerical fluxes $F^*$ and $G^*$gives,
\begin{equation}
\label{eq:weak_form}
\left( \pa{Q_{H}}{t} , \phi \right)_{D^{k}} = \left( \pa{\phi}{x} , F \right)_{D^{k}} + \left( \pa{\phi}{y} , G \right)_{D^{k}} - \left( F^{*}n_{x}+G^{*}n_{y} , \phi \right)_{\partial D^{k}} + (S, \phi)_{D^{k}}.
\end{equation}

\begin{figure}[h]
 \begin{center}
  \subfloat[Interpolation nodes]{%
    \begin{minipage}[c]{0.3\linewidth}
      \centering%
      \includegraphics[trim=4cm 6cm 4cm 7cm,clip=true,width=\textwidth]{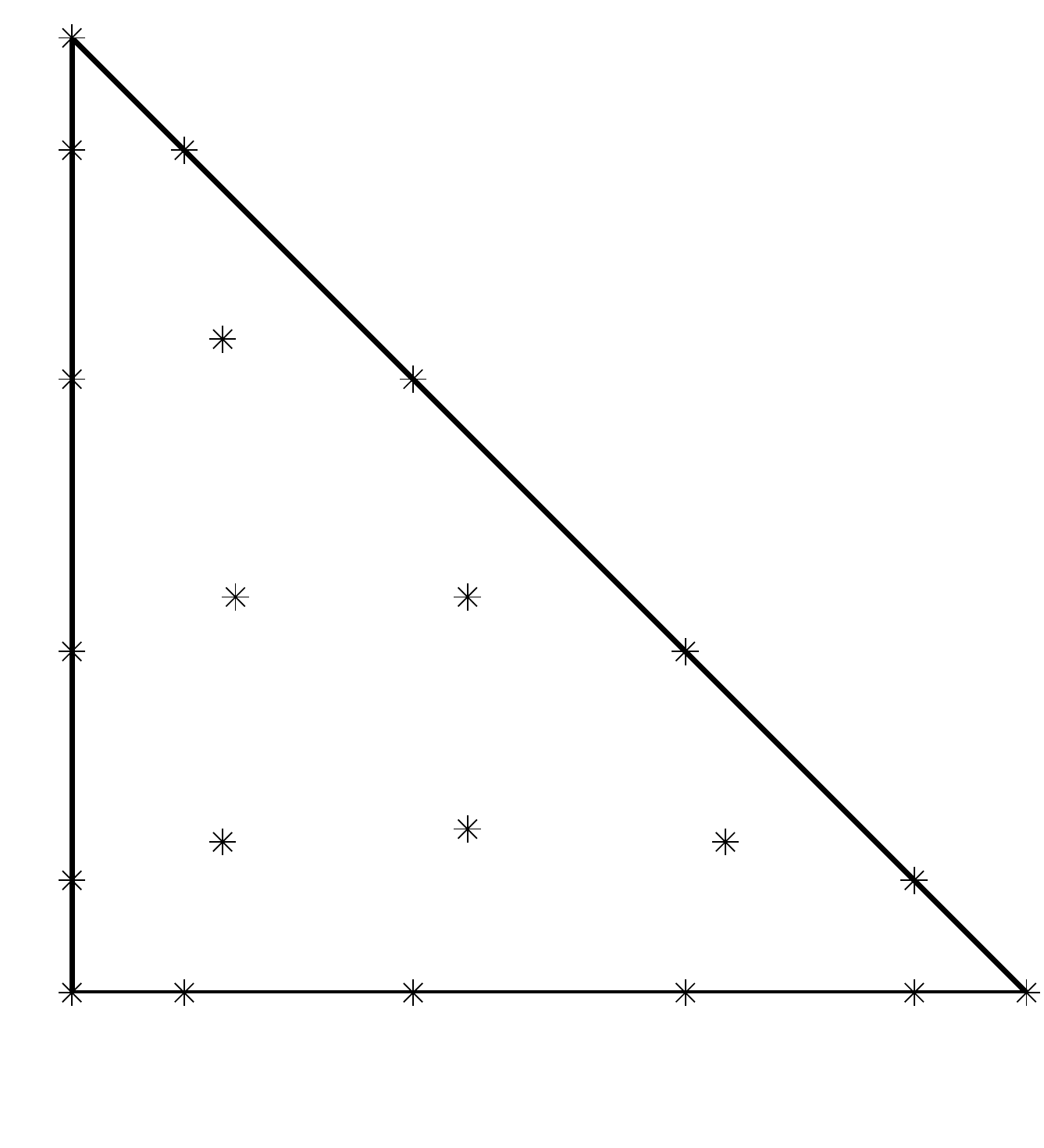}
    \end{minipage}}
  \hspace{0.3cm}
  \subfloat[Cubature nodes]{%
    \begin{minipage}[c]{0.3\linewidth}
      \centering%
      \includegraphics[trim=4cm 6cm 4cm 7cm,clip=true,width=\textwidth]{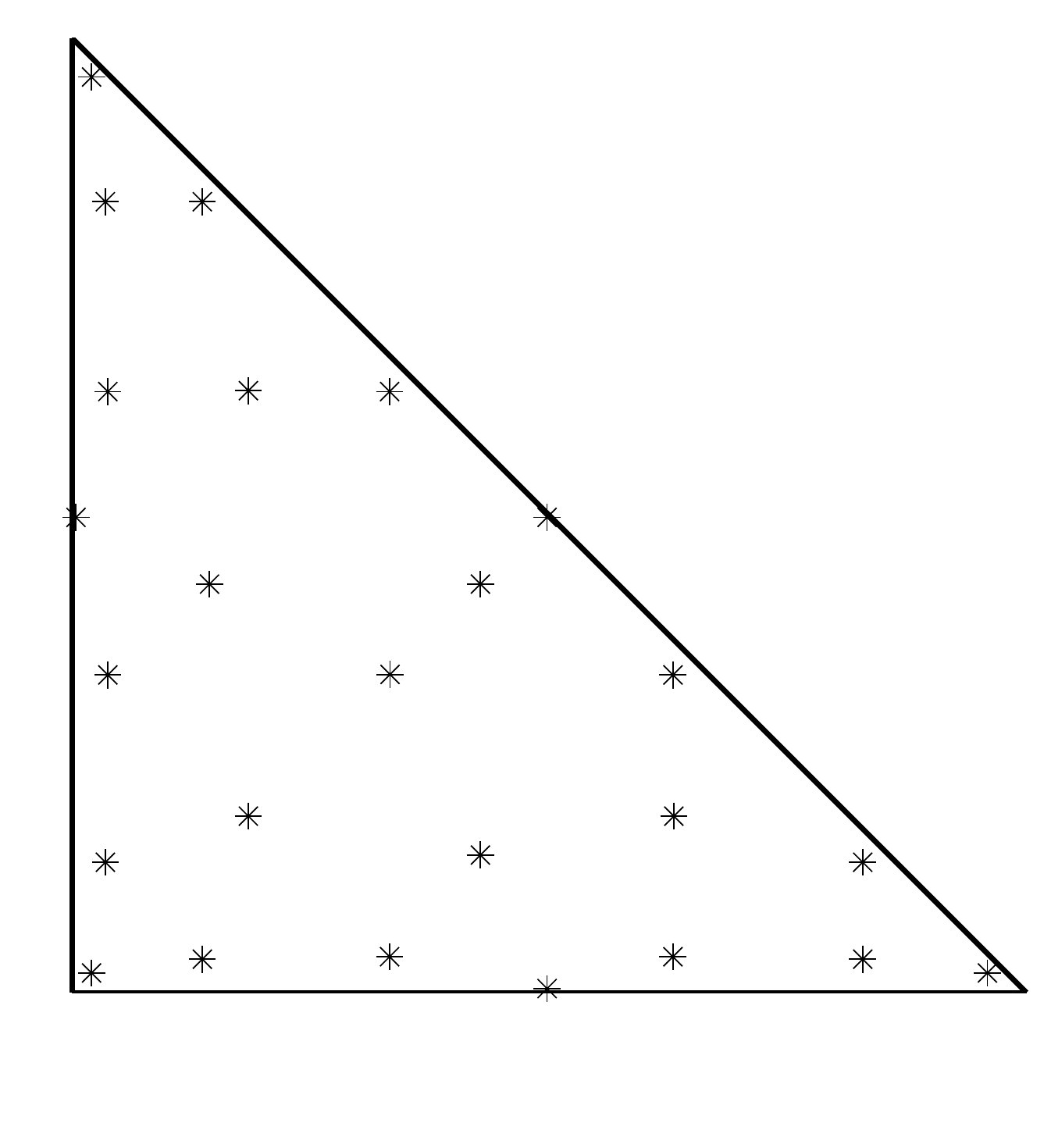}
    \end{minipage}}
  \hspace{0.3cm}
  \subfloat[Gauss quadrature nodes]{%
    \begin{minipage}[c]{0.3\linewidth}
      \centering%
      \includegraphics[trim=4cm 6cm 4cm 7cm,clip=true,width=\textwidth]{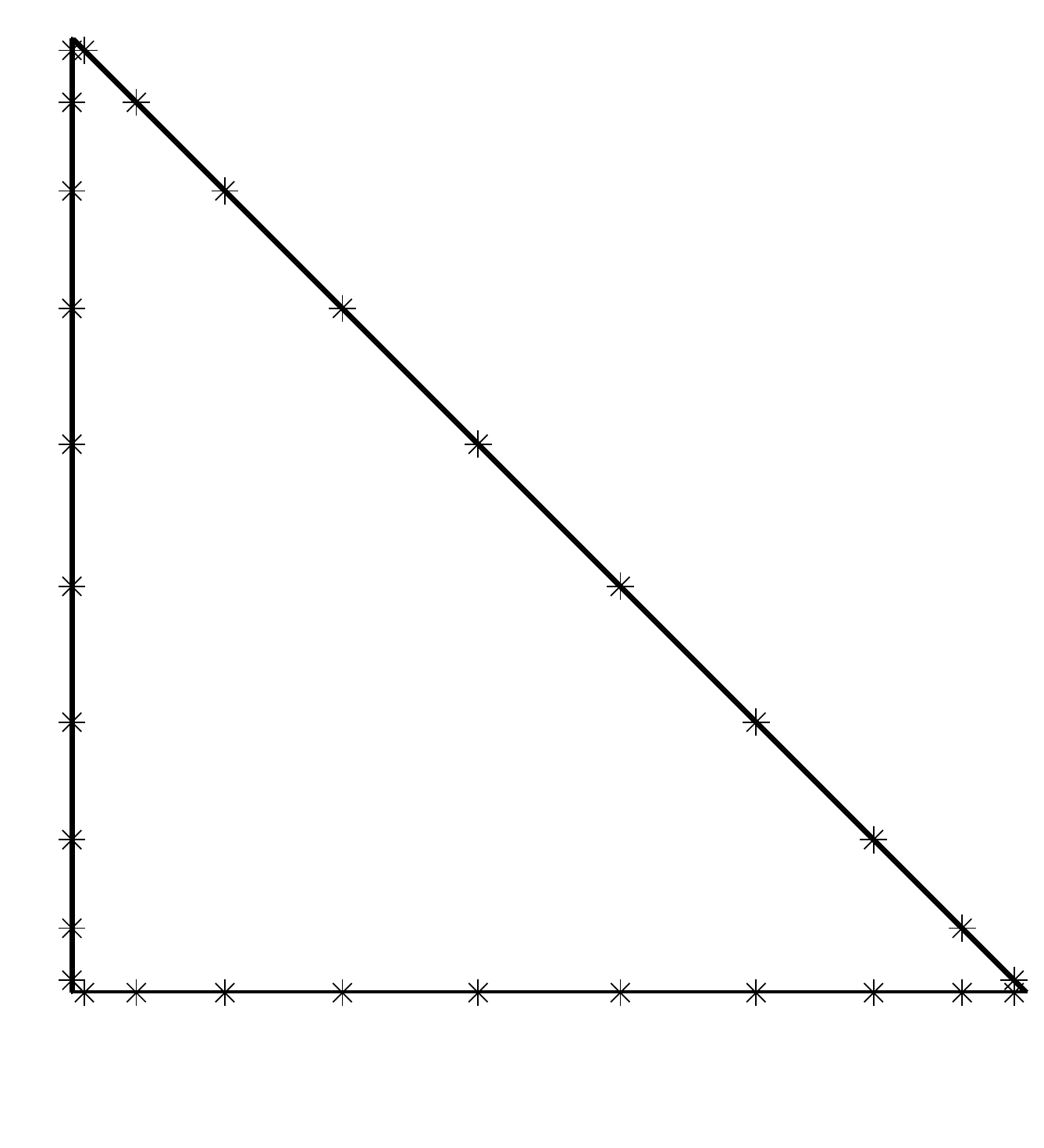}
    \end{minipage}}
    \caption{The distribution of interpolation nodes, cubature integration
nodes, and Gauss quadrature integration nodes for the interpolating polynomial
order 5, and the integration order 10.}
    \label{fig:nodes}
 \end{center}
\end{figure}

Here $(\cdot , \cdot)_{D^{k}}$ represents the inner product taken over the
element $D^{k}$, while $(\cdot , \cdot)_{\partial D^{k}}$ represents the
inner product taken over the boundary of the element $D^{k}$ denoted by $\partial
D^{k}$. We use Lagrange polynomials with Warp \& Blend interpolation nodes \cite{hesthaven2008nodal} as the basis for the polynomial space in each element, and use well-balanced local Lax-Friedrich flux \cite{xing2010positivity} to compute $F^{*}$ and $G^{*}$. The volume integrals are computed using cubature rules for triangles \cite{cools1999monomial}, while the surface integrals are computed using Gauss quadrature rules (see Fig. (\ref{fig:nodes})), leading to a system of ordinary differential equations given by,
\begin{equation}
\label{eq:ode}
\frac{d Q_{H}}{dt} = \mathcal{R}(Q_{H}) = \mathcal{N}(Q_{H}) + \mathcal{S}(Q_{H}^{g,+},Q_{H}^{g,-}),
\end{equation}

Here, $\mathcal{R}$  is the spatial discretization operator and $\mathcal{N}$, $\mathcal{S}$ are nonlinear operators corresponding to volume and surface integrations. $Q_{H}^{g}$ is a vector of  the state variables at the Gauss quadrature nodes. $Q_{H} ^{g,+}$ and $Q_{H}^{g,-}$ represent the positive and negative traces of the solution at the element interfaces. The ODEs  in  Eq. (\ref{eq:ode}) are integrated using a local time-stepping scheme to obtain the solution at a given time $t$. 

\section{Local time stepping}
\label{sec:mrab}
Explicit time integration is conditionally stable; the time step $\Delta t$, has to satisfy Courant-Friedrichs-Lewy  criterion,
\begin{equation}
\label{eq:cfl}
\Delta t = \displaystyle\min_{D^k} \left\{ \frac{H^{k}}{a} C \right\}, \qquad \text{} a = |(u,v)|+\sqrt{gh}.
\end{equation}

Here, $H^k$ is the characteristic length of the element $D^k$, $a$ is the wave speed, and $C$ is a constant that depends on the stability region of the time integration scheme and spatial polynomial order.

For a global scheme, the smallest length scale in the mesh  determines the overall time step, which leads to a very small  allowable time step size. To allow each element to be integrated with its own allowable time step,  we adopt a local time stepping based on the multi-level third order AB scheme \cite{gear1984multirate,godel2010gpu} to nonlinear ODEs. For the shallow water applications, mesh resolutions can vary in few orders of magnitude because a very fine mesh is used near the shore regions to resolve  high frequencies in the waves, while a coarse mesh is sufficient to resolve the flow in the deep oceans. The global allowable time step is few orders of magnitude smaller than the allowable time step for a coarse element.

Local time stepping is done efficiently by grouping the elements into levels based on characteristic lengths and integrating the elements in a level with a fixed time step size. After mesh generation, the individual time step for each element is evaluated and  elements with an allowable time step $\Delta t$ are grouped in to a level $l$ if $2^{l-1} \Delta t _{\text{min}} \le \Delta t < 2^{l} \Delta t _{\text{min}}$. All elements in level $l$ are integrated with a time step  $2^{l-1} \Delta t _{\text{min}}$.

For a given initial condition $Q^{0}$, the time integration procedure is described in the Algorithm (\ref{alg:mrab}).
\begin{algorithm}
  \caption{Multi level Adam-Bashforth explicit time stepping}
 \label{alg:mrab}
  \begin{algorithmic}[1]
   \algnotext{EndIf}
  \algnotext{EndFor}
    \State  $Q_{H} ^{0} = \Pi_{H}Q^{0}$
    \For  {$n=0, 1, 2,....$}
    \For  {$l=Nlevels, Nlevels- 1, ..... 1$}
    \For  {$nstep = 1, 2, ..., Nsteps(l)$}
    \State        $\left. Q_{H} ^{n + \frac{nstep}{Nsteps(l)}} \right |_l = \left. Q_{H} ^{n + \frac{nstep-1}{Nsteps(l)}} \right |_l + \Delta t_{l} \displaystyle\sum_{s=1} ^{3} \alpha _{s} \mathcal{R} \left(Q_{H} ^{n + \frac{nstep-s}{Nsteps(l)}} \right)$

    \EndFor
    \EndFor
    \EndFor
  \end{algorithmic}
\end{algorithm}
$\Pi_{H} Q^{0}$ is the projection of the initial conditions on to the space of piecewise discontinuous polynomials, $n$ is a time level, and $Nlevels$ is the total number of levels in multi-rate scheme. The elements with largest time step are integrated first, followed by  elements with  smaller time steps. $\alpha _{i}$'s are coefficients corresponding to the $3^{rd}$ order Adams-Bashforth linear multi-step method. $\left. Q_{H} ^{n + \frac{nstep}{Nsteps(l)}} \right |_l$ is the numerical solution corresponding to the the elements in level $l$ and time $t =t^{n} + nstep \, \Delta t_{l}$.

From the Eq. (\ref{eq:cfl}), it can be observed that size of the time step for the stability changes over time with the change in characteristic speeds and the elements need to be regrouped in to levels. Grouping the elements into levels carries a nontrivial  setup cost and hence is inefficient to do this every time step. In this work, the allowable time step for each element is fixed throughout the simulations, and our experiments indicate that changes in characteristic speeds are not strong enough to cause instabilities in the solutions and do not pollute the spectral accuracy. However, if needed, elements can be regrouped after every few time steps to satisfy the updated CFL conditions.



\section{Well-balancing, positivity preserving and slope limiting schemes}
\label{sec:PP}

\subsection{Well-balancing}
In the numerical solutions of shallow water equations it is possible to have  nonzero flux derivatives while maintaining a steady state. In these cases, the derivatives are balanced by the source terms from the bathymetry distribution.   The numerical schemes have to exactly balance these terms (exact C-property). We choose a numerical flux proposed in  \cite{xing2006high, xing2010positivity} to satisfy the exact C-property. \begin{eqnarray}
h^{*,\pm} &=& \max(0, \, \, h^{\pm} + B^{\pm} - \max(B^{+}, B^{-})), \nonumber \\
Q^{*,\pm} &=& \begin{pmatrix} h^{*, \pm} \\ h^{*, \pm} u^{\pm} \\ h^{*, \pm} v^{\pm} \end{pmatrix}.
\end{eqnarray}

Here $+$ and $-$ refer to positive and negative trace of field values at the element interfaces. $Q^{*}$ is intermediate field vector. The well balanced fluxes are given by,
\begin{eqnarray}
\label{eq:well_balance}
F^{*, \pm} &=& F_{LF}(Q^{*,+}, Q^{*,-}) + \begin{pmatrix} 0 \\ \frac{g}{2} (h^{\pm})^2 - \frac{g}{2}(h^{*,\pm})^2 \\ 0 \end{pmatrix}, \nonumber \\
G^{*, \pm} &=& G_{LF}(Q^{*,+}, Q^{*,-}) + \begin{pmatrix} 0 \\ 0 \\ \frac{g}{2} (h^{\pm})^2 - \frac{g}{2}(h^{*,\pm})^2 \end{pmatrix},
\end{eqnarray}
where $F_{LF}, G_{LF}$ are Lax-Friedrich fluxes at the element interfaces.

\subsection{Positivity preserving limiter}
Another  concern in the numerical simulation of the shallow water equations is the appearance of dry areas where no water is present. The shallow water equations implicitly assume non-negative depth of water and hence this property has to be ensured by the numerical scheme. Otherwise, the flux Jacobian will have non-real eigenvalues and the PDE will no longer be hyperbolic. For a finite volume method, the water depth is a constant in each element. An element is flagged as dry if the fluid height is below a threshold value and is not considered as part of the simulation until the fluid height reaches the threshold value. This is no longer valid for a high order DG discretization since the height is typically not a constant but a non-monotonic high order polynomial. A wetting drying scheme that conserves mass and momentum is required. Positivity preserving algorithms have been developed under the assumptions of piecewise linear \cite{bunya2009wetting} for the shallow water equations and rectangular elements \cite{zhang2010positivity} for the Euler equations. These were later extended to high order polynomials on rectangular elements for shallow water flows \cite{xing2010positivity}. These schemes can not be extended to triangles due to the assumptions made on numerical integration schemes.  We adopt the approach presented for linear polynomials and triangular meshes in  \cite{bunya2009wetting}, and extend it to an arbitrary order polynomial approximation.

We represent the positive preserving operator with $M \Pi_{H}$, this operator is used along with a TVB limiter $\Lambda \Pi_{H}$ at every intermediate time step in the multi-rate integration. The procedure is outlined in the Algorithm (\ref{alg:mrab_with_limiter}).
\begin{algorithm}[h!]
  \caption{Multi-level Adam-Bashforth  time stepping with PP  and TVB limiting}
\label{alg:mrab_with_limiter}
  \begin{algorithmic}[1]
   \algnotext{EndIf}
  \algnotext{EndFor}
    \State  Set $Q_H ^{0} = \Lambda \Pi_{H} M \Pi_{H} (\Pi_{H} Q^{0})$.
    \For  {$n=0, 1, 2,....$}
    \For  {$l=Nlevel, Nlevel- 1, ....., 1$}
    \For  {$nstep = 1, 2, ..., Nsteps(l)$}
    \State $\left. Q_{H}^{n+\frac{nstep}{Nsteps(l)}} \right|_{l}            = \Lambda \Pi_{H} M \Pi_{H} \left\{            \left. Q_{H}^{n+\frac{nstep-1}{Nsteps(l)}} \right|_{l} + \Delta t_{l} \displaystyle\sum _{s = 1} ^{3} \alpha_{s}  \mathcal{R}\left(Q_{H}^{n+\frac{nstep-s}{Nsteps(l)}} \right) \right\}$

    \EndFor
    \EndFor
    \EndFor
  \end{algorithmic}
\end{algorithm}

 First, we discuss in detail the operator $M \Pi_{H}$ that ensures positivity of the fluid height. To reduce the notational complexity, we represent the numerical solution $Q_{H}$ with $q$. A given polynomial $q = (h, \, hu, \, hv)^{T}$ in an element $D^k$  is modified such that the polynomial corresponding to the fluid height ($h$) is positive ($\ge h_{0}$) in the element, where $h_{0}$ is a threshold value of the fluid height for considering an element/region as dry land. At a time level $n$, the polynomial solution $q^{n} = (h^{n}, \, hu^{n}, \, hv^{n})^{T}$ is modified to obtain $\tilde{q^{n}}$ using Algorithm (\ref{alg:PP}).

\begin{algorithm}[h!]
  \caption{Positivity preserving limiter $M\Pi Q_H := \tilde{q}^n$}
  \label{alg:PP}
  \begin{algorithmic}[1]
  \algnotext{EndIf}
  \algnotext{EndFor}
    \State $h^{n}_{\min} = \displaystyle \min _{j} {h_j}$ \hfill \{ minimum height \}
    \If {$h^{n}_{\min} > \epsilon$} \hfill \{ if nodal values for height are +ve \}
    \State $ \tilde{q}^n _{j} = q^{n},  \qquad \forall \boldsymbol{x} \in D^{k}$ \hfill \{ do not modify the solution \}
    \Else
    \State $q^{n,1} = \Pi_1 ^{n} q^n $ \hfill \{ project the solution to linear polynomial \}
    \If {$\bar{h^n} < h_0$} \hfill \{ if dry element \}
    \State $h^n = h_0$, ${hu}^n = 0$, and ${hv}^n = 0$ \hfill \{ modify the solution \}

    \Else \If {$\bar{h^n} \ge h_0$} \hfill \{ if mean is positive \}
    \State $\tilde{q}^{n}_{j} = \theta (q^{n,1}_{j} - \bar{q}^{n}) + \bar{q}^{n}, \qquad \theta = \min \left\{ 1, \dfrac{\bar{h}^{n, 1} - h_{0}}{\bar{h}^{n}-h_{\min}^{n, 1}}\right\},$

    \State $h_{\min}^{n,1} = \displaystyle \min _{i} h^{n, 1}_{i}$, \hfill $\forall \boldsymbol{x}_i \in$ \{ vertices of the triangle \}

    \EndIf
    \EndIf
    \EndIf
  \end{algorithmic}
\end{algorithm}
It is easy to see that linear polynomial representation of the fluid height $h^{n,1}$ is positive ($\ge h_{0}$) at all interpolation nodes, and so at the integration points.
Note that local conservation of mass and momentum are not violated for wet or partially wet elements. This is achieved by keeping the constant modes in the orthogonal polynomials unaltered during the limiting.

\subsection{TVB slope limiter}
A slope limiter  removing high frequency oscillations in numerical solutions is applied to avoid instabilities due to nonlinear effects. We use Cockburn and Shu's characteristic based TVB limiter \cite{shu1987tvb, cockburn1998runge}  designed for Runge-Kutta methods. However, the TVB limiter does not ensure the positivity of the solution. So we perform a post processing of the solution to ensure positivity of the  fluid height. We  discuss the post processing step without going into the details of the TVB limiter.

\begin{figure}[h!]
 \begin{center}
  \subfloat[TVB limiter]{%
    \begin{minipage}[c]{0.45\linewidth}
      \centering%
         \includegraphics[trim=4cm 15cm 0cm 0cm,clip=true,width=\textwidth]{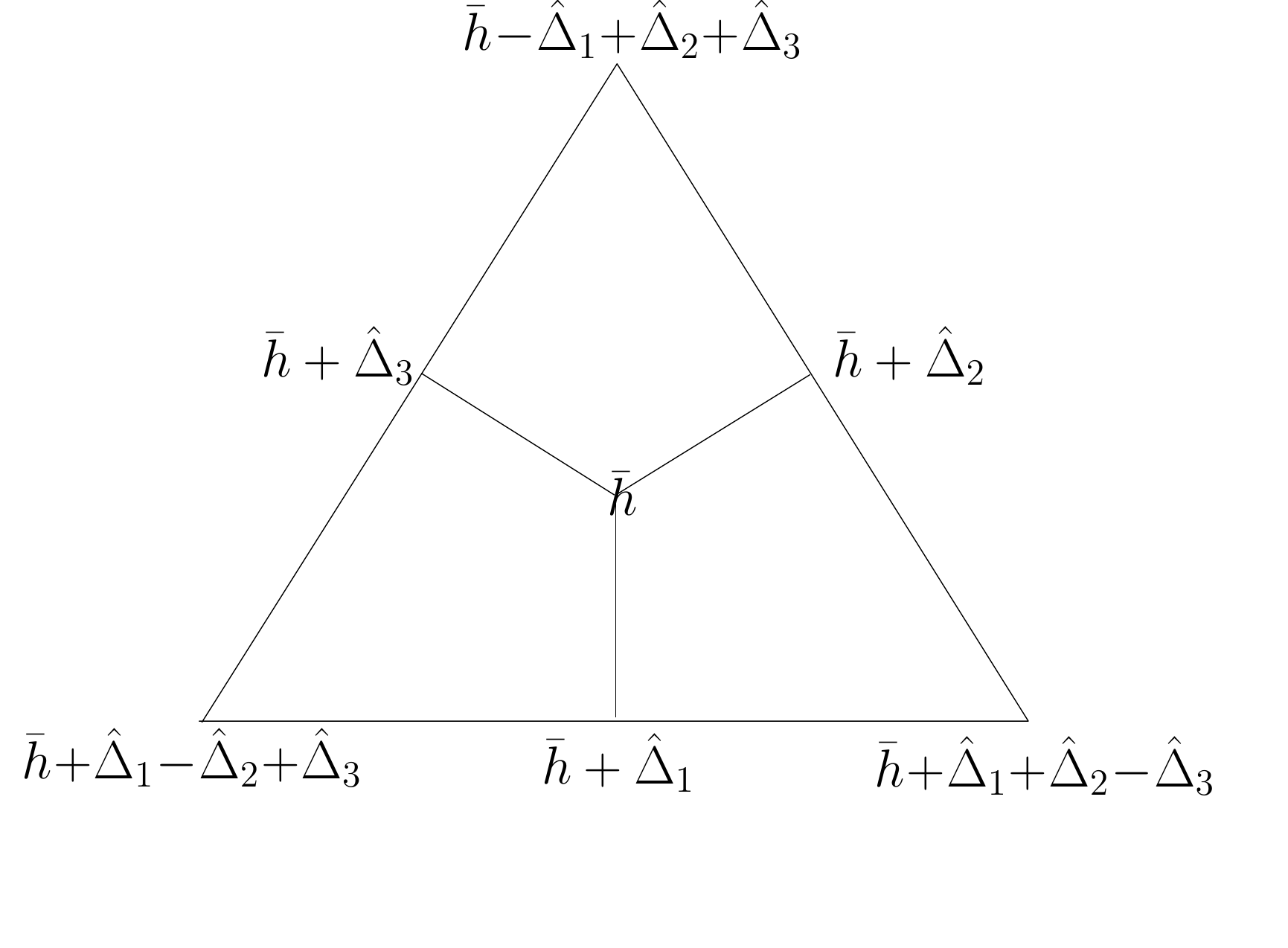}
         \end{minipage}}
  \hspace{1cm}
  \subfloat[Modified TVB limiter]{%
    \begin{minipage}[c]{0.45\linewidth}
         \centering%
         \includegraphics[trim=4cm 15cm 0cm 0cm,clip=true,width=\textwidth]{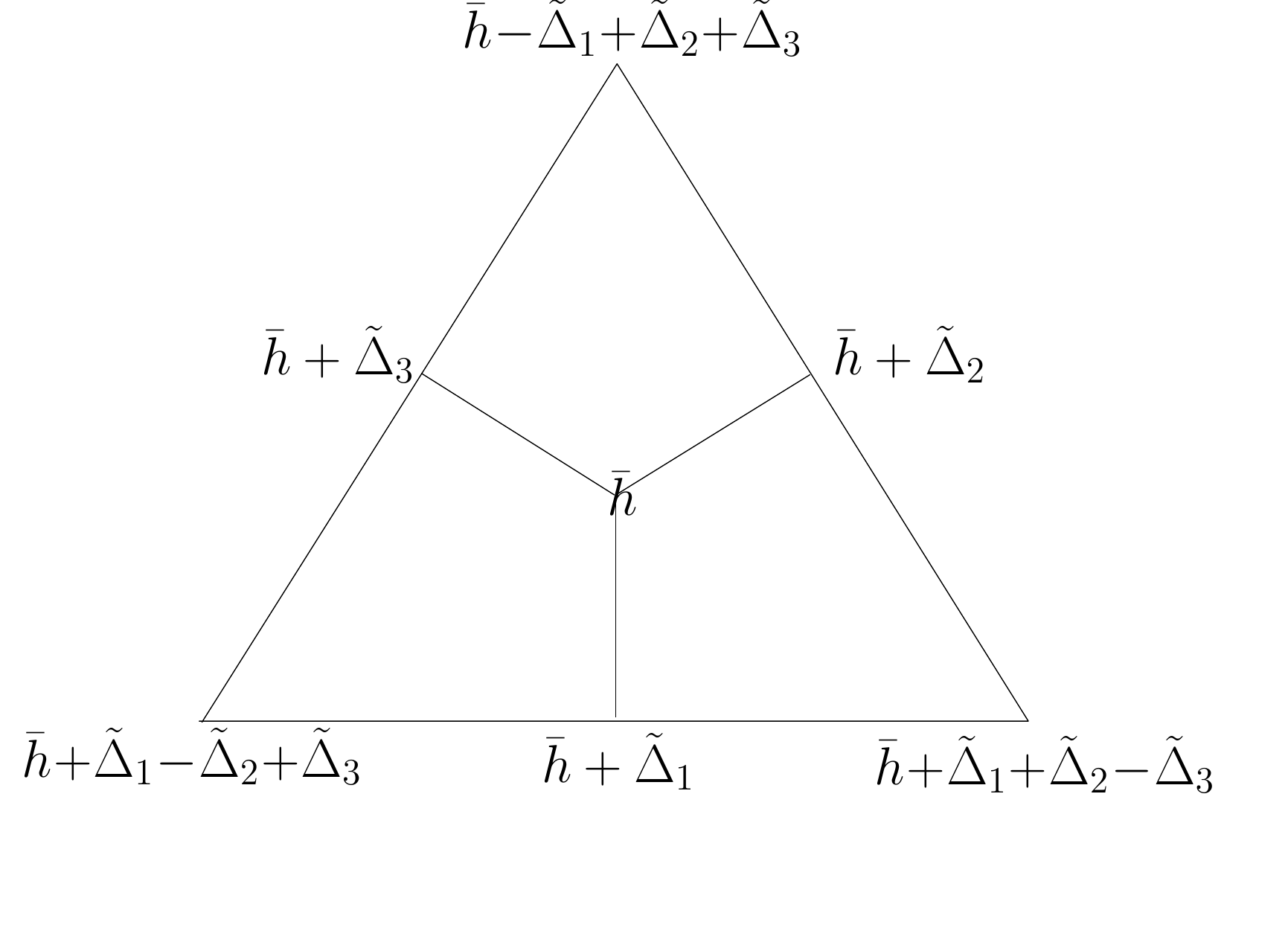}
            \end{minipage}}
\caption{Modification of TVB limiter to ensure positivity of the fluid height}
\label{fig:modified_TVB}
 \end{center}
\end{figure}

TVB slope limiter restricts the local solution to a linear polynomial. As discussed earlier, a linear polynomial attains its maximum and minimum at the vertices. We modify $\hat{\Delta}_{1}, \hat{\Delta}_{2}, \text{and }\hat{\Delta}_{3}$ (refer to \cite{cockburn1998runge} for the notation) to   $\tilde{\Delta}_{1}, \tilde{\Delta}_{2}, \text{and }\tilde{\Delta}_{3}$
without changing the average of these to ensure the positivity of the fluid height at the vertices. The values at the vertices are given in Fig. (\ref{fig:modified_TVB}).
\begin{eqnarray}
\label{eq:modified_TVB}
\tilde{\Delta}_{i} &=& \bar{\Delta} + \theta ( \hat{\Delta}_{i}- \bar{\Delta}), \qquad \text{for } i = \{1,2,3\}, \nonumber \\
\bar{\Delta} &=&  \text{avg}(\hat{\Delta}_{1}, \hat{\Delta}_{2}, \hat{\Delta}_{3}), \nonumber \\
\theta &=& \frac{\bar{h} + \bar{\Delta} - h_{0}}{ \bar{\Delta} - \displaystyle \min_{i,j,k = \{1,2,3\}, i \neq j \neq k} \{ -\hat{\Delta}_i + \hat{\Delta}_j + \hat{\Delta}_k\}}.
\end{eqnarray}

As reported in \cite{ern2008well, bunya2009wetting}, a positive preserving (and/or a wetting drying treatment) and a slope limiter may artificially activate each other causing instability. To avoid this instability, the TVB limiter is not applied for elements that are considered dry. For high order approximations, the TVB limiter may be artificially activated in the immediate neighbors of the dry elements also. Therefore, we do not apply the TVB limiter for dry elements and the immediate neighbors of dry elements as well. The robustness of these limiters and effect on the solution accuracy are demonstrated in the next section.



\section{Verification}
\label{sec:tests}
In this section, we use several test cases with known analytical solutions to  study the discussed limiters.

\subsection{Accuracy test for smooth solution: Couette flow}
\begin{figure}[h!]
\begin{center}
  \subfloat[Couette flow]{%
    \begin{minipage}[c]{0.45\linewidth}
      \centering%
      \includegraphics[trim=2cm 6cm 1cm 6cm,clip=true,width=\textwidth]{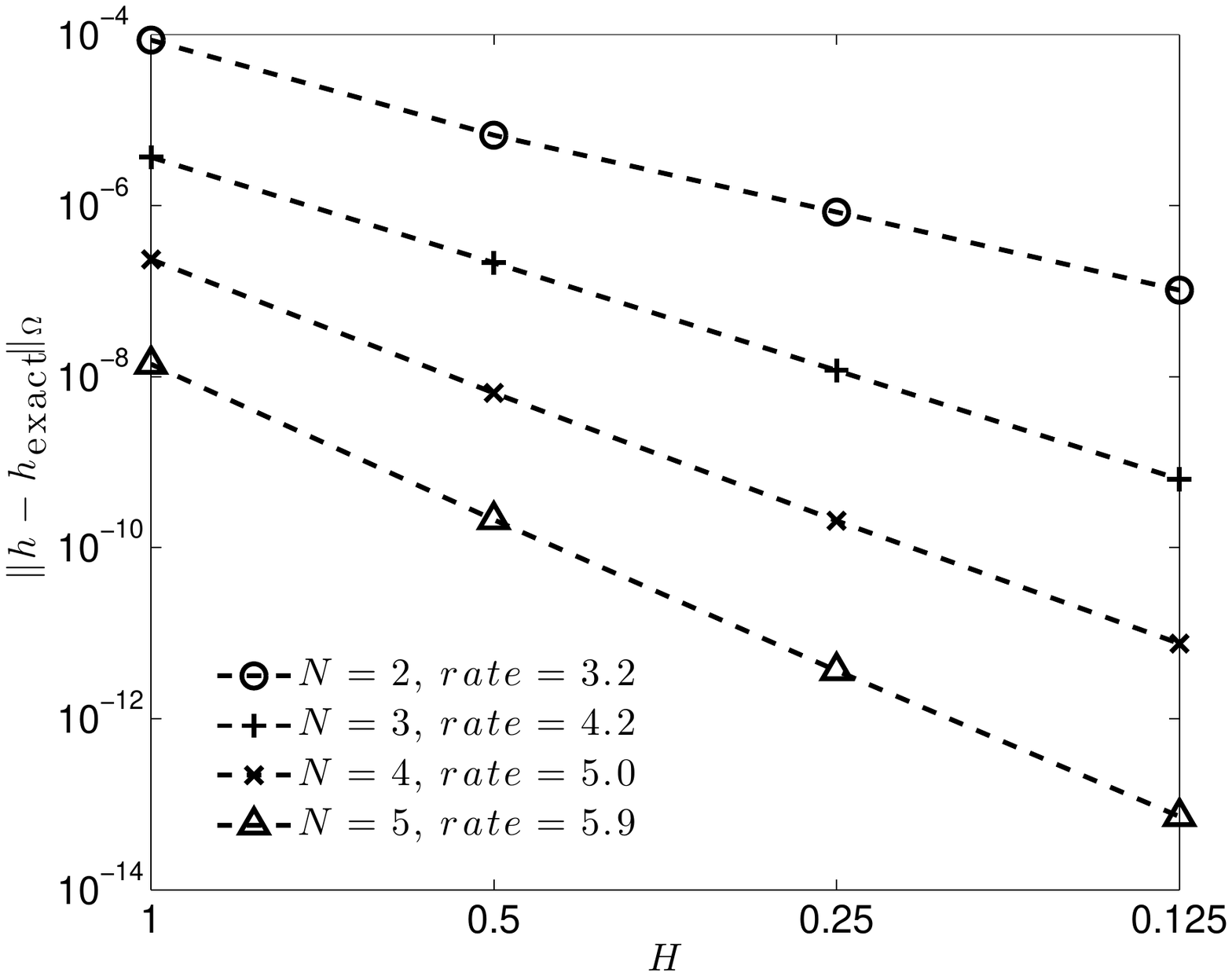}
    \end{minipage}}
  \hspace{0.5cm}
  \subfloat[Translating vortex]{%
    \begin{minipage}[c]{0.45\linewidth}
      \centering%
      \includegraphics[trim=2cm 6cm 1cm 6cm,clip=true,width=\textwidth]{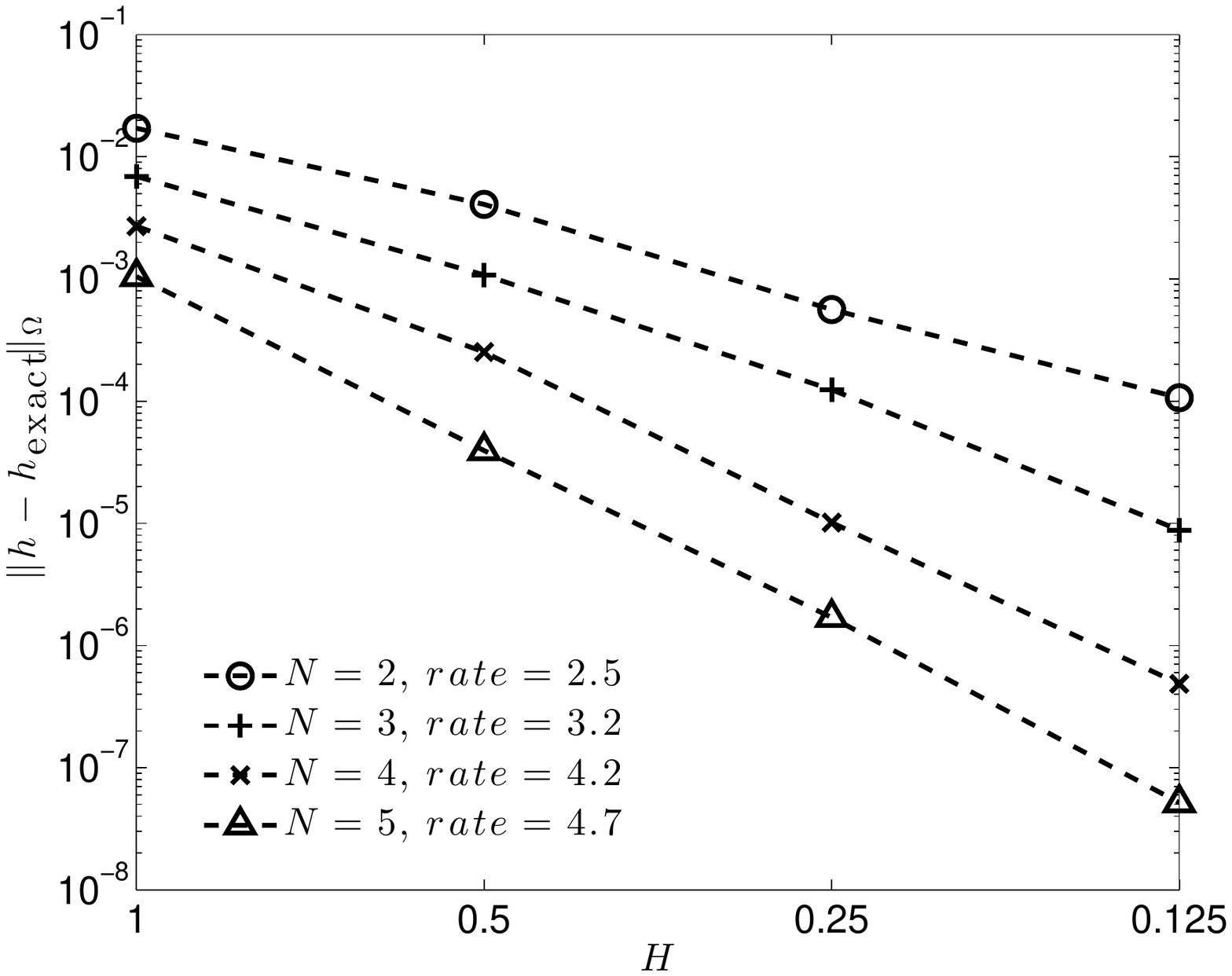}
    \end{minipage}}
\caption{Plot of $L^2$ error vs $H/H_{0}$ for smooth solutions, with polynomial orders $N=2,3,4, \text{ and } 5$}
\label{fig:vortex_couette_convergence}
\end{center}
\end{figure}

We consider Couette flow between two concentric cylinders spinning at different velocities with the exact solution  given by,
\begin{equation}
h = 1, \qquad
u = -\sin(\theta) u_{\theta}, \qquad
v = \cos(\theta) u_{\theta}
\end{equation}
where azimuthal velocity $u_{\theta}$ and bathymetry $B$ are given by $u_{\theta} = \frac{1}{75}\left(-r + \frac{16}{r}\right)$ and  $B = \frac{1}{75^2} \left(\frac{r^2}{2} - 32\log(r) - \frac{128}{r^2}\right)$. Here, $\theta = \tan ^{-1}(\frac{y}{x})$ and $r=\sqrt{x^2+y^2}$. Since the solution corresponds to a steady state, the simulations are started with the exact solution as the initial conditions and run for a long time ($t=10s$) to compute the spatial errors. We observe $L^2$-errors decaying like $\mathcal{O}(H^{N+1/2})$ (see Fig. (\ref{fig:vortex_couette_convergence})). The order of convergence for each polynomial order is computed by determining a polynomial that fits the error in the least square sense.

\begin{figure}[h!]
\begin{center}
  \subfloat[]{%
    \begin{minipage}[c]{0.2\linewidth}
      \centering%
      \includegraphics[trim=4cm 7cm 4cm 7cm,clip=true,width=\textwidth]{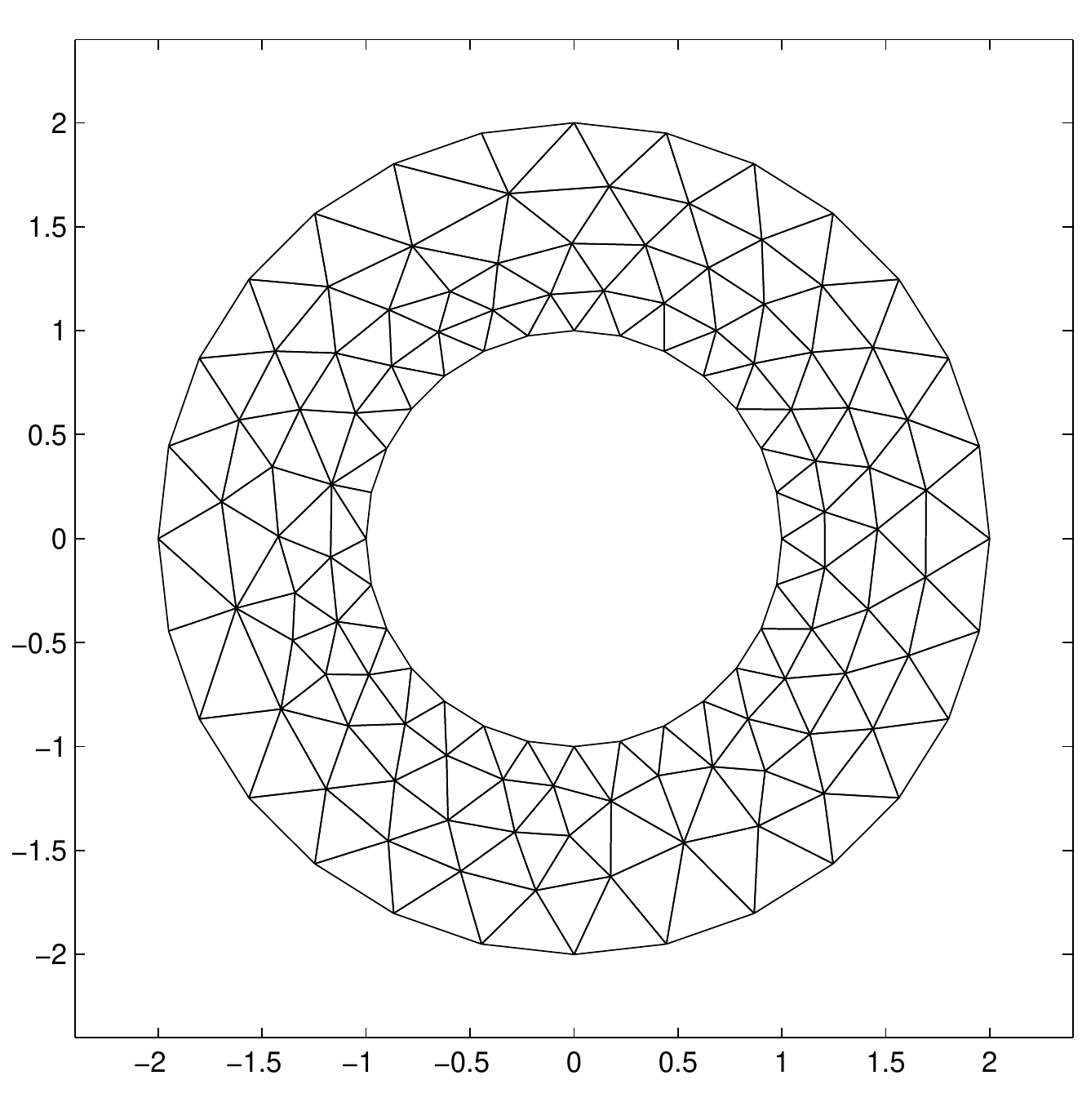}
    \end{minipage}}
  \hspace{0.2cm}
  \subfloat[]{%
    \begin{minipage}[c]{0.2\linewidth}
      \centering%
      \includegraphics[trim=4cm 7cm 4cm 7cm,clip=true,width=\textwidth]{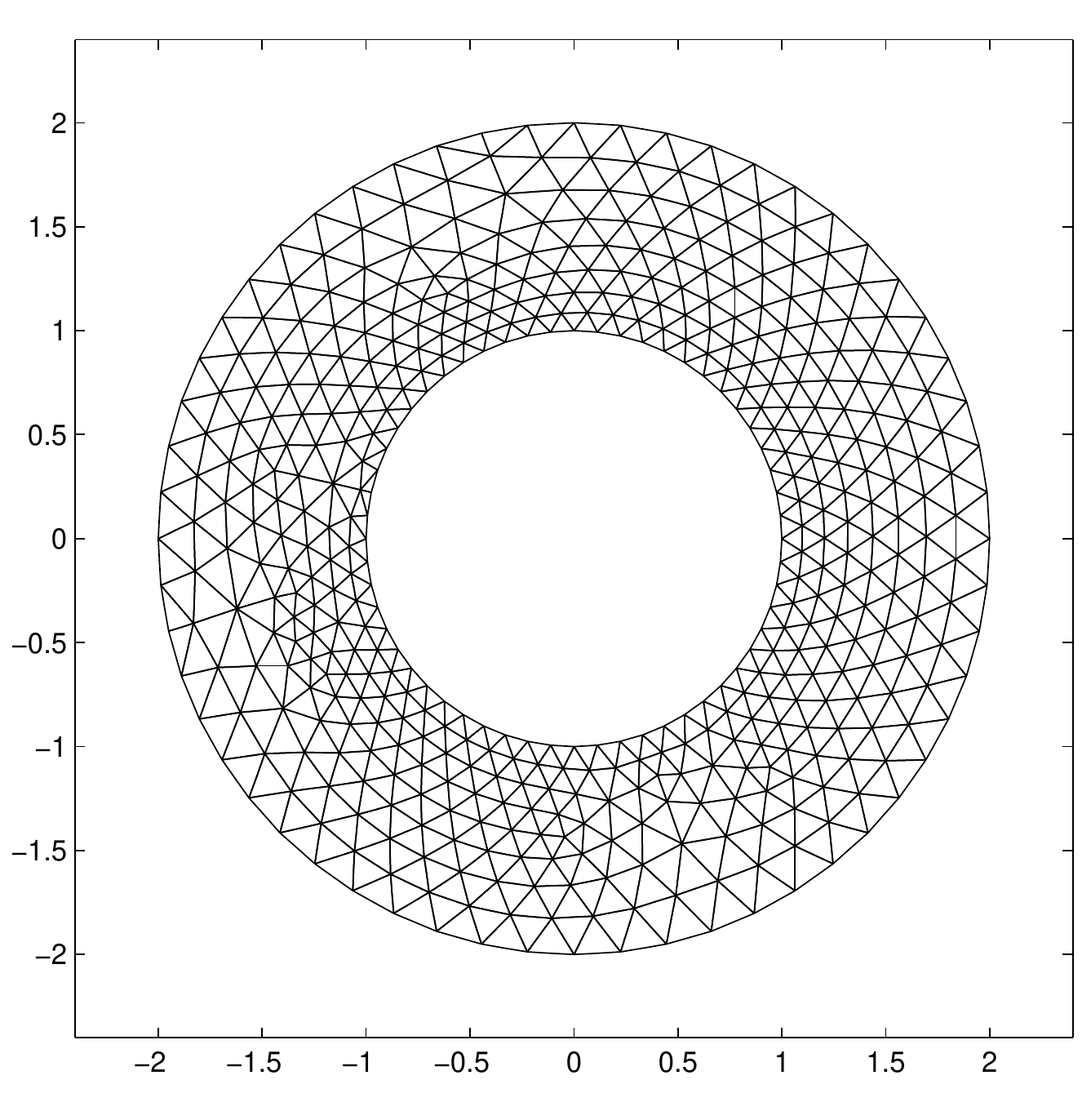}
    \end{minipage}}
  \hspace{0.2cm}
  \subfloat[]{%
    \begin{minipage}[c]{0.2\linewidth}
      \centering%
      \includegraphics[trim=4cm 7cm 4cm 7cm,clip=true,width=\textwidth]{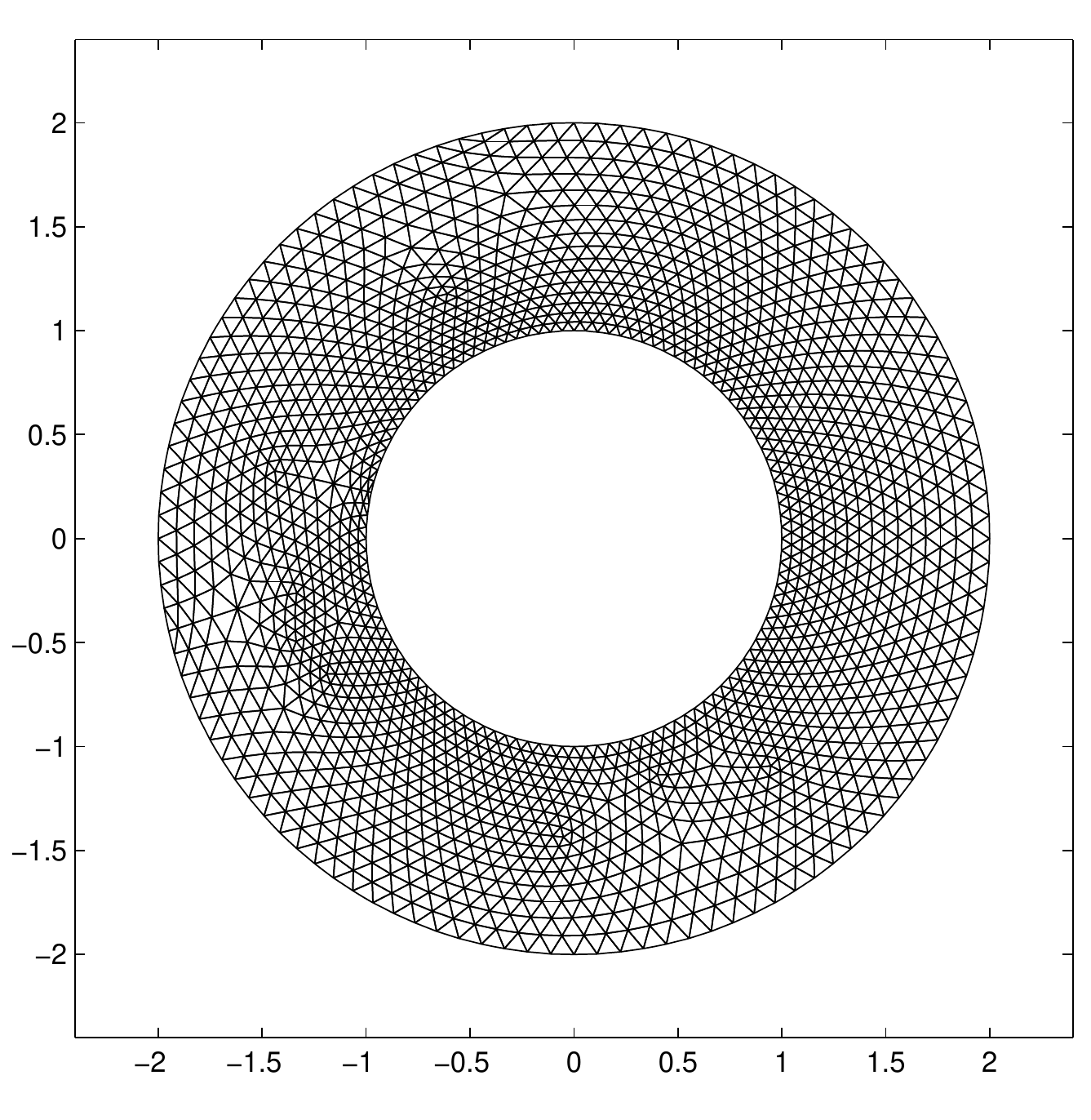}
    \end{minipage}}
  \hspace{0.2cm}
  \subfloat[]{%
    \begin{minipage}[c]{0.2\linewidth}
      \centering%
      \includegraphics[trim=4cm 7cm 4cm 7cm,clip=true,width=\textwidth]{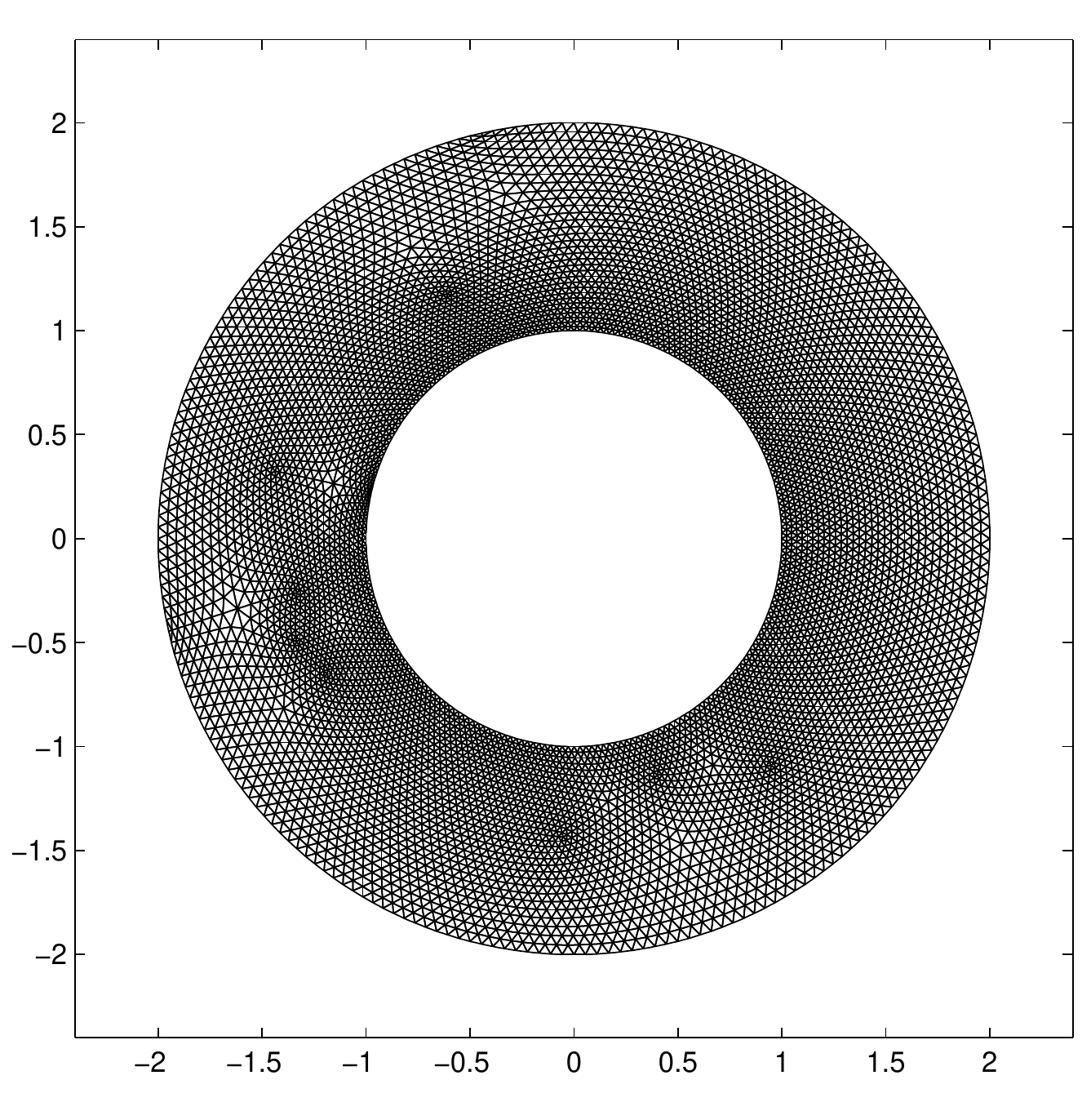}
    \end{minipage}}
  \caption{The sequence of meshes used to perform convergence analysis of rotating Couette flow between two concentric cylinders}
  \label{fig:couette_meshes}
\end{center}
\end{figure}

\subsection{Accuracy test for smooth solution : Translating Vortex}
We use this problem to test the accuracy of the time dependent solutions. An isentropic vortex translates in space with a constant speed and satisfies the two dimensional Euler equations \cite[p.209]{hesthaven2008nodal}. By replacing the density ($\rho$) with the fluid height ($h$), and choosing the gas constant $\gamma = 2$ and gravity $g=2$, an analytical solution for the shallow water equations is obtained.
\begin{figure}[h!]
\begin{center}
  \subfloat[]{%
    \begin{minipage}[c]{0.2\linewidth}
      \centering%
      \includegraphics[trim=2cm 6cm 1cm 6cm,clip=true,width=\textwidth]{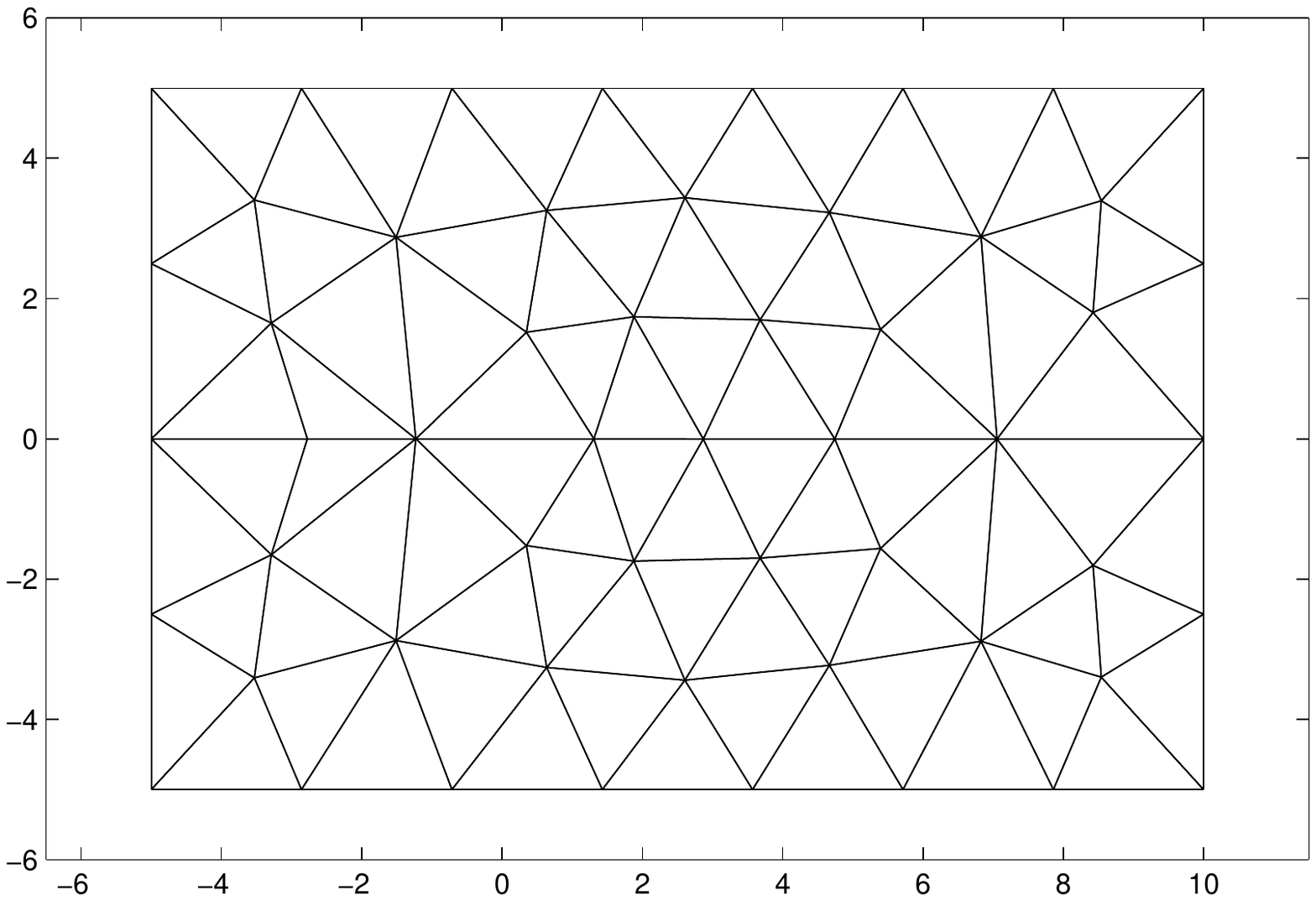}
    \end{minipage}}
  \hspace{0.2cm}
  \subfloat[]{%
    \begin{minipage}[c]{0.2\linewidth}
      \centering%
      \includegraphics[trim=2cm 6cm 1cm 6cm,clip=true,width=\textwidth]{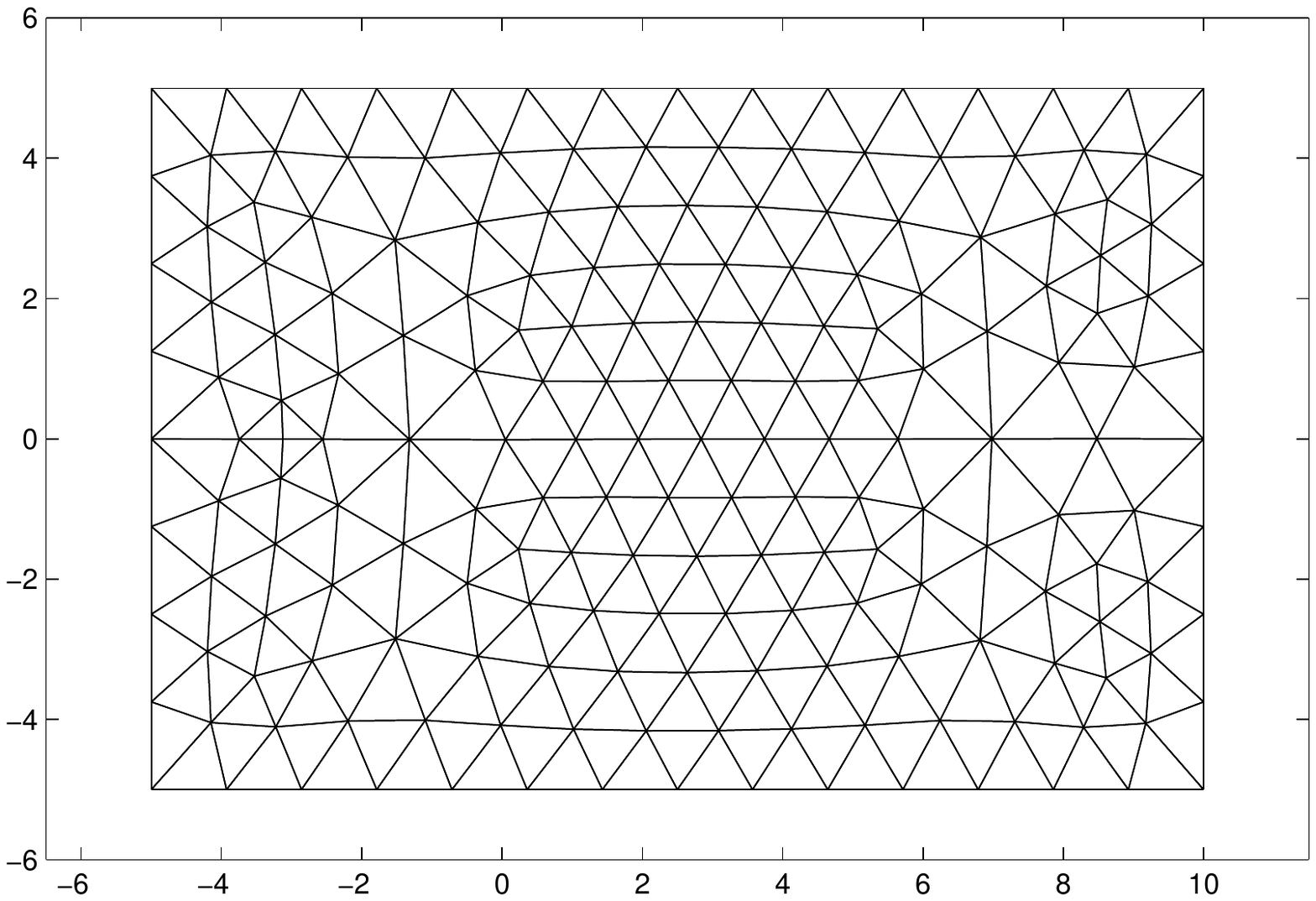}
    \end{minipage}}
  \hspace{0.2cm}
  \subfloat[]{%
    \begin{minipage}[c]{0.2\linewidth}
      \centering%
      \includegraphics[trim=2cm 6cm 1cm 6cm,clip=true,width=\textwidth]{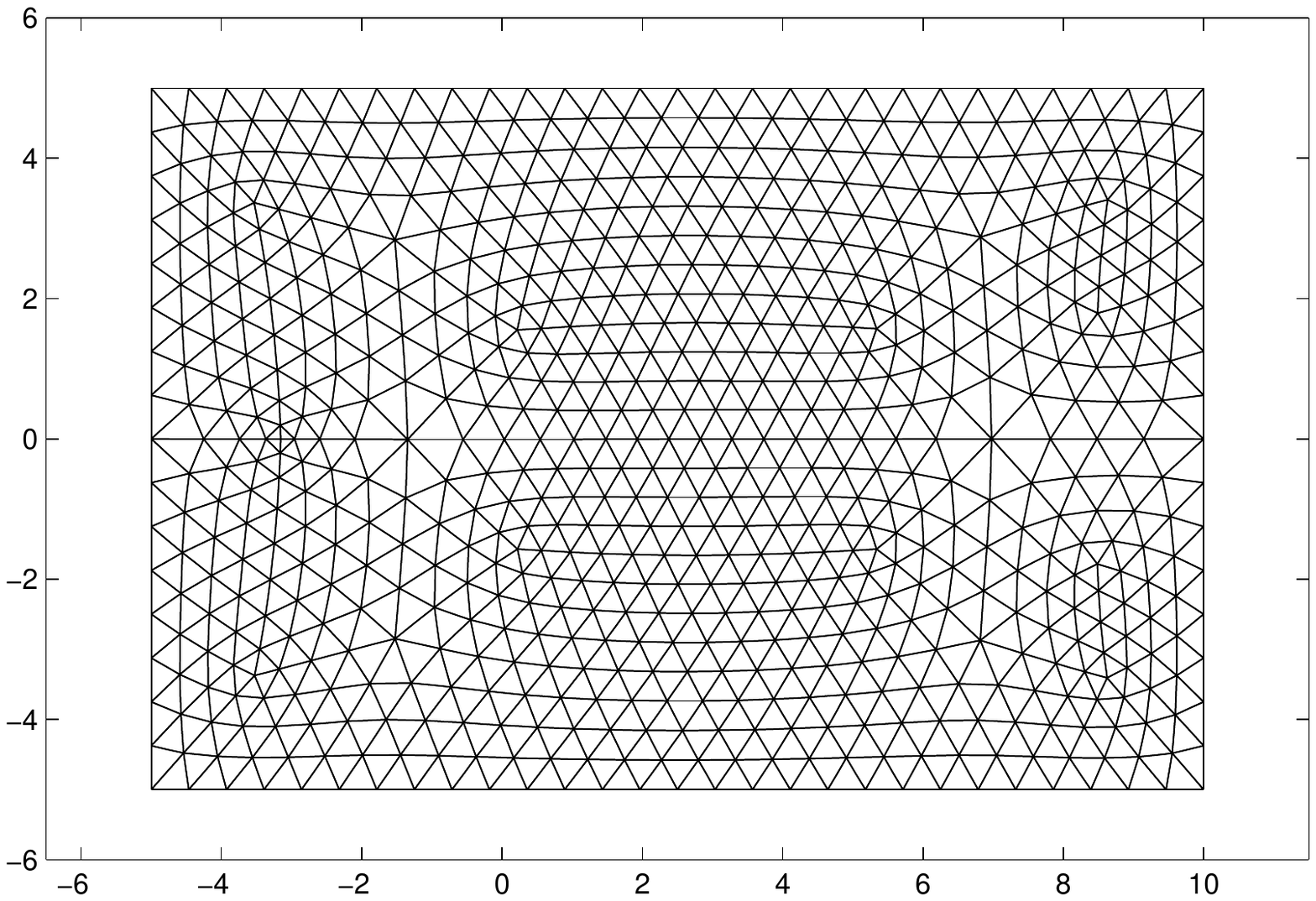}
    \end{minipage}}
  \hspace{0.5cm}
  \subfloat[]{%
    \begin{minipage}[c]{0.2\linewidth}
      \centering%
      \includegraphics[trim=2cm 6cm 1cm 6cm,clip=true,width=\textwidth]{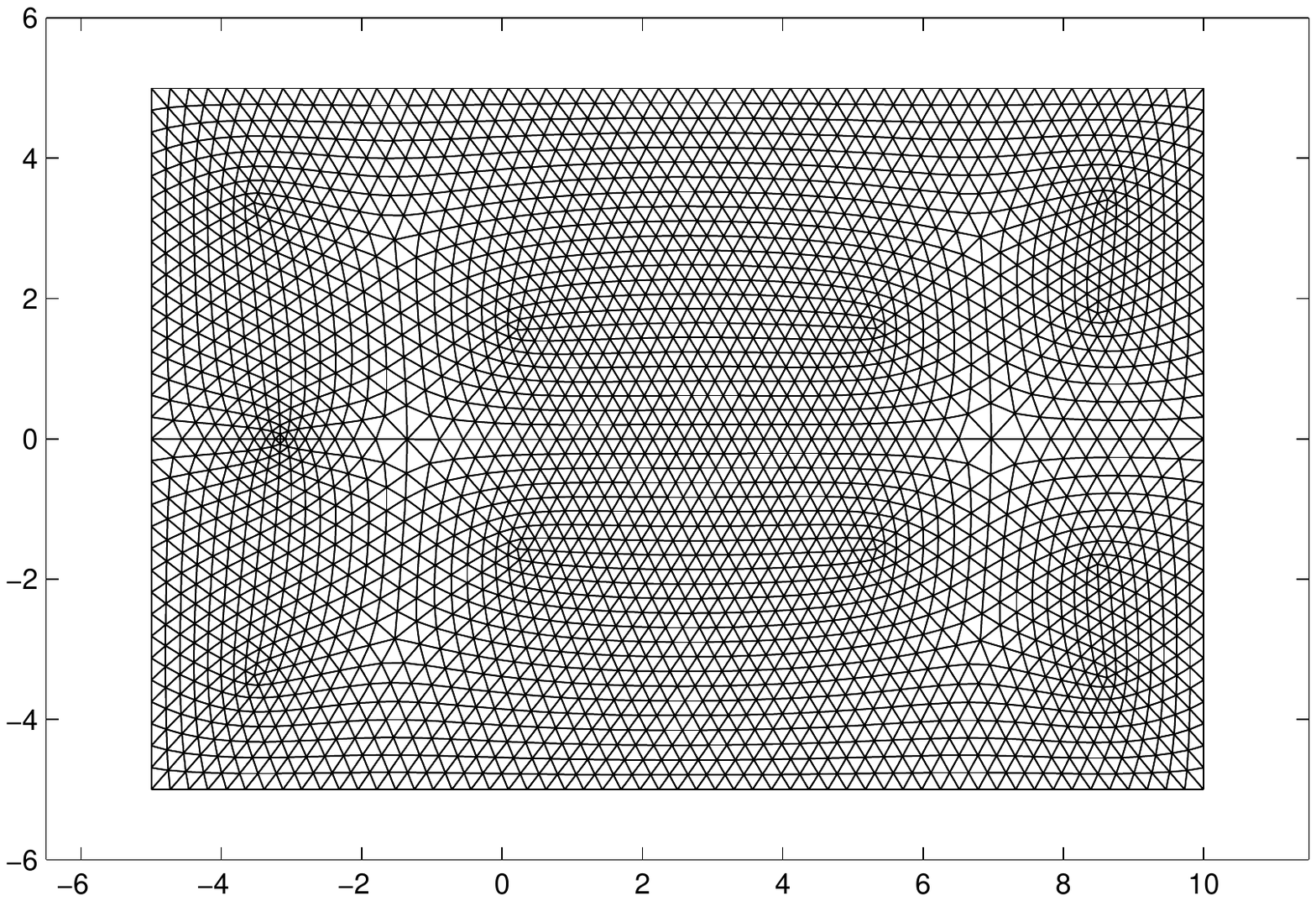}
    \end{minipage}}
  \caption{The sequence of meshes used to perform convergence analysis for a  translating vortex in a rectangular domain}
  \label{fig:vortex_meshes}
\end{center}
\end{figure}
\begin{eqnarray}
h = 1 - \frac{\beta ^ 2}{32 \pi ^2} e^{2(1-r^2)}, \quad
u = 1 - \beta e^{1-r^2} \frac{y - y_{0}}{2 \pi}, \quad
v = \beta e^{1-r^2} \frac{x - t - x_{0}}{2 \pi},
\end{eqnarray}
where $r=\sqrt{(x-t-x_{0})^2 + (y-y_{0})^2}$ and the bathymetry is  a constant. A rectangular domain $[-5,10] \otimes [-6,6] \in \mathbb{R}^{2}$ is chosen as the domain for this problem and both the initial conditions and Dirichlet boundary conditions are set to the exact solution. The sequence of meshes used are given in the Fig. (\ref{fig:vortex_meshes}). Again we observe that $L^2$-error converges like $\mathcal{O}(H^{N+1/2})$  (see Fig. (\ref{fig:vortex_couette_convergence})) as expected for  hyperbolic PDEs \cite{johnson1986analysis}.

\subsection{Parabolic bowl}
We use this test case \cite{ern2008well,thacker1981some} to study the effect of the positive preserving limiter on  solution accuracy. The bathymetry is given by a parabola, $b(x,y) = \alpha r^2$, where $r = \sqrt{x^2 + y^2}$. The exact solution for fluid height is non zero for $r < \sqrt{(X+Y\cos{\omega t})/\alpha (X^2 - Y^2)}$, where $X > 0$, $\lvert Y \rvert < X$ and $\omega ^2 = 8 g \alpha$. The nonzero exact solution is given by,
\begin{eqnarray}
\label{eq:pb_exact}
h(x,y,t) &=& \frac{1}{X+Y\cos(\omega t)} + \alpha (Y^2 - X^2) \frac{r^2}{(X+Y\cos(\omega t))^2}, \nonumber \\
u(x,y,t) &=& -\frac{Y \omega \sin(\omega t)}{X + Y \cos (\omega t)} \frac{x}{2}, \nonumber \\
v(x,y,t) &=& -\frac{Y \omega \sin(\omega t)}{X + Y \cos (\omega t)} \frac{y}{2}.
\end{eqnarray}

\begin{figure}[h!]
\begin{center}
  \subfloat[Global error]{%
    \begin{minipage}[c]{0.45\linewidth}
      \centering%
      \includegraphics[trim=2cm 6cm 1cm 6cm,clip=true,width=\textwidth]{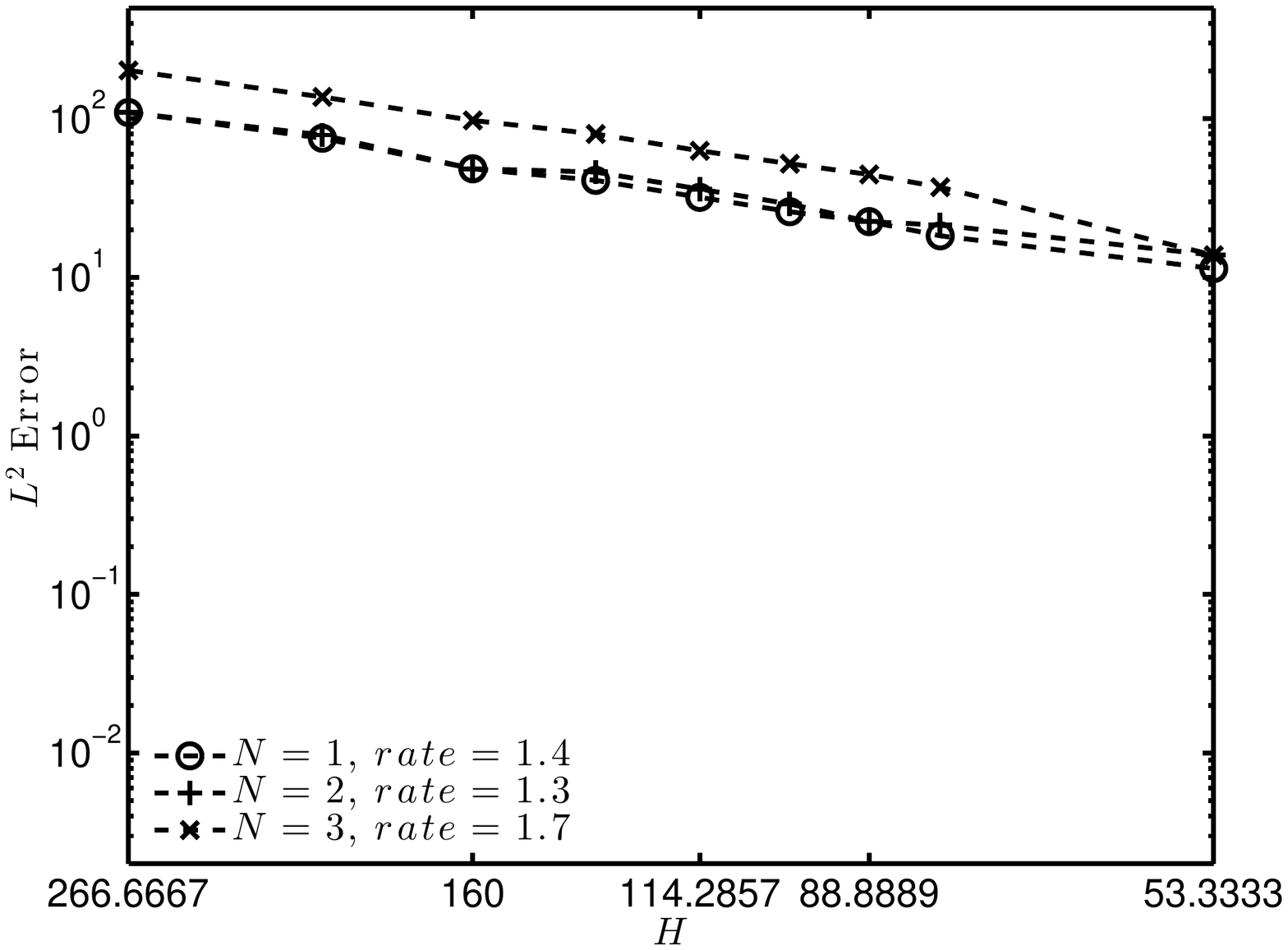}
    \end{minipage}}
  \hspace{0.5cm}
  \subfloat[Local error]{%
    \begin{minipage}[c]{0.45\linewidth}
      \centering%
      \includegraphics[trim=2cm 6cm 1cm 6cm,clip=true,width=\textwidth]{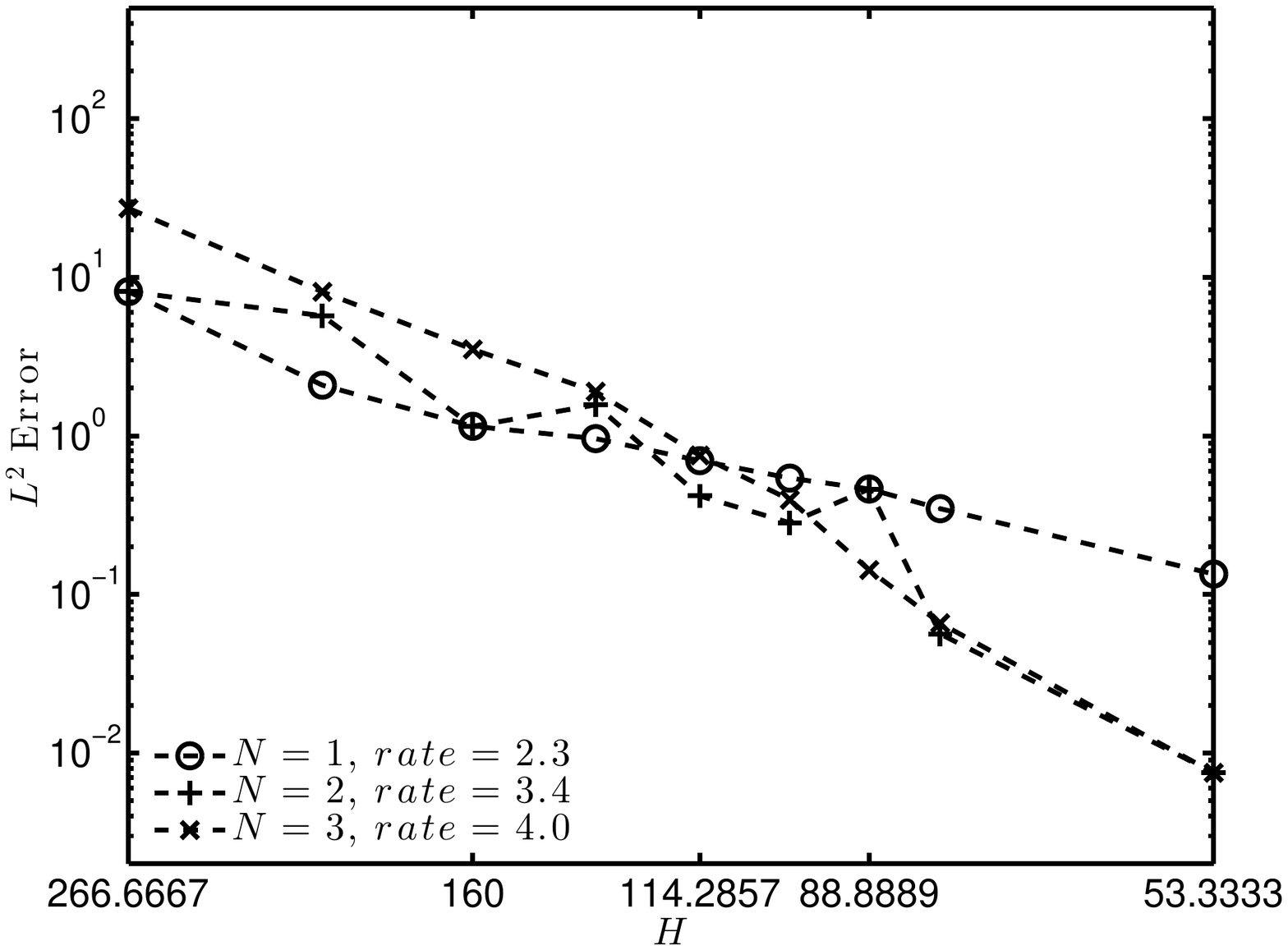}
    \end{minipage}}
\caption{Parabolic bowl test case : (a) global $L^2$ error (b) $L^2$ error for region $r < 500 m$}
\label{fig:parabolic_bowl_local_global_L2}
\end{center}
\end{figure}

\begin{figure}[h!]
\begin{center}
  \subfloat[H=160m, N=2, solution]{%
    \begin{minipage}[c]{0.45\linewidth}
      \centering%
      \includegraphics[trim=0cm 6cm 1cm 6cm,clip=true,width=\textwidth]{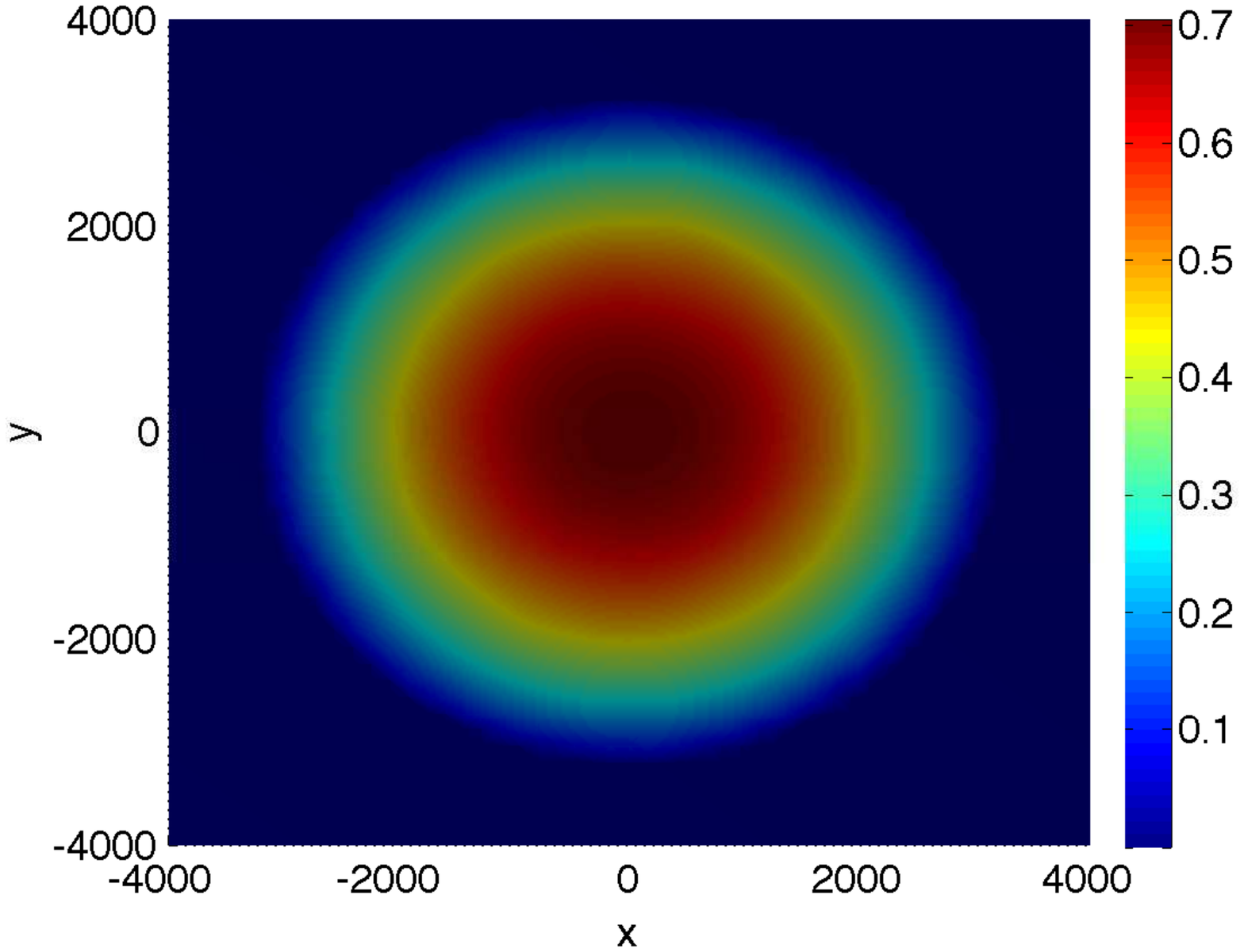}
  \end{minipage}}
  \hspace{0.5cm}
  \subfloat[H=160m, N=2, error]{%
    \begin{minipage}[c]{0.45\linewidth}
      \centering%
      \includegraphics[trim=0cm 6cm 1cm 6cm,clip=true,width=\textwidth]{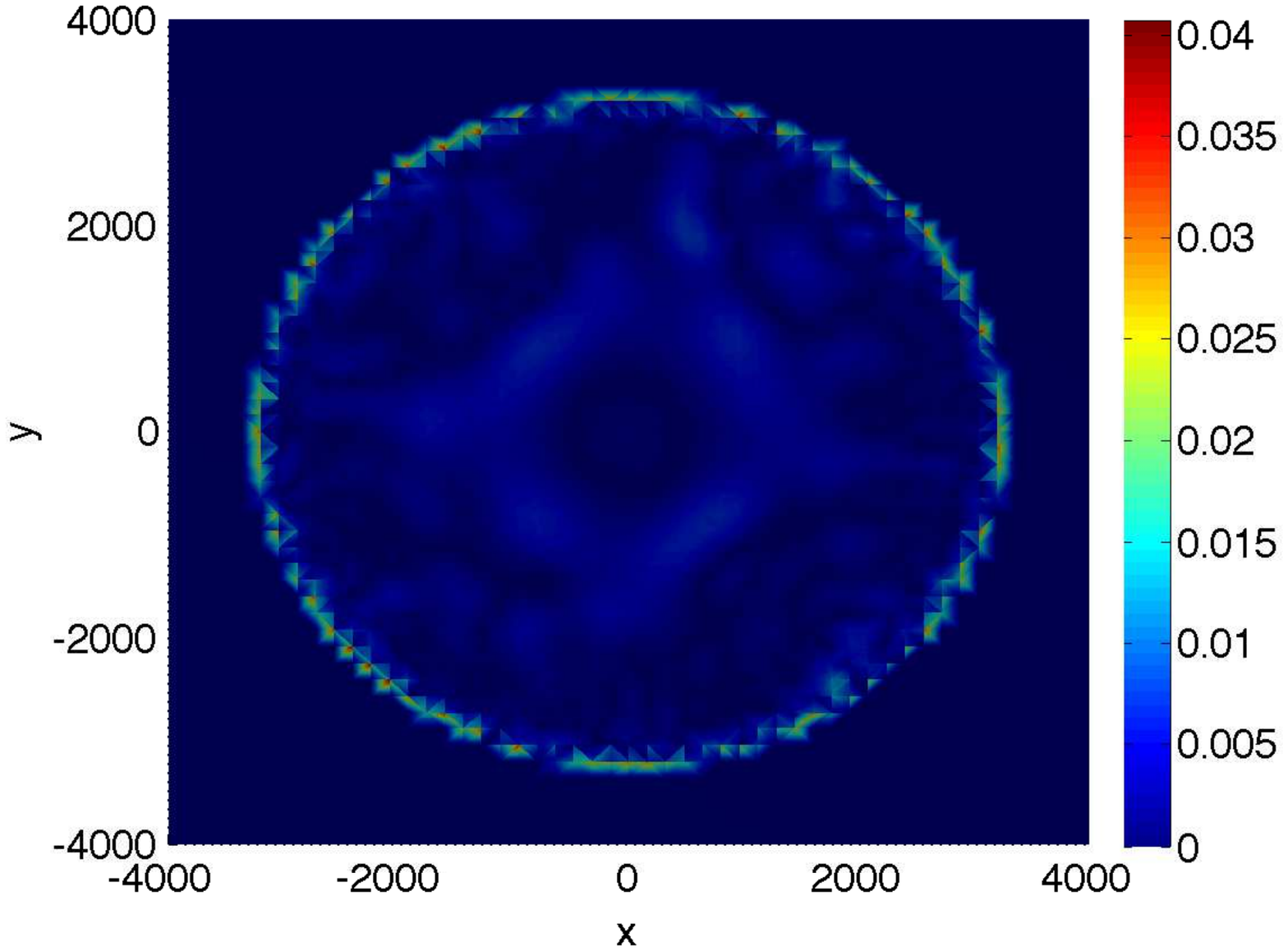}
  \end{minipage}} \\
  \subfloat[H=40m, N=2, solution]{%
    \begin{minipage}[c]{0.45\linewidth}
      \centering%
      \includegraphics[trim=0cm 6cm 1cm 6cm,clip=true,width=\textwidth]{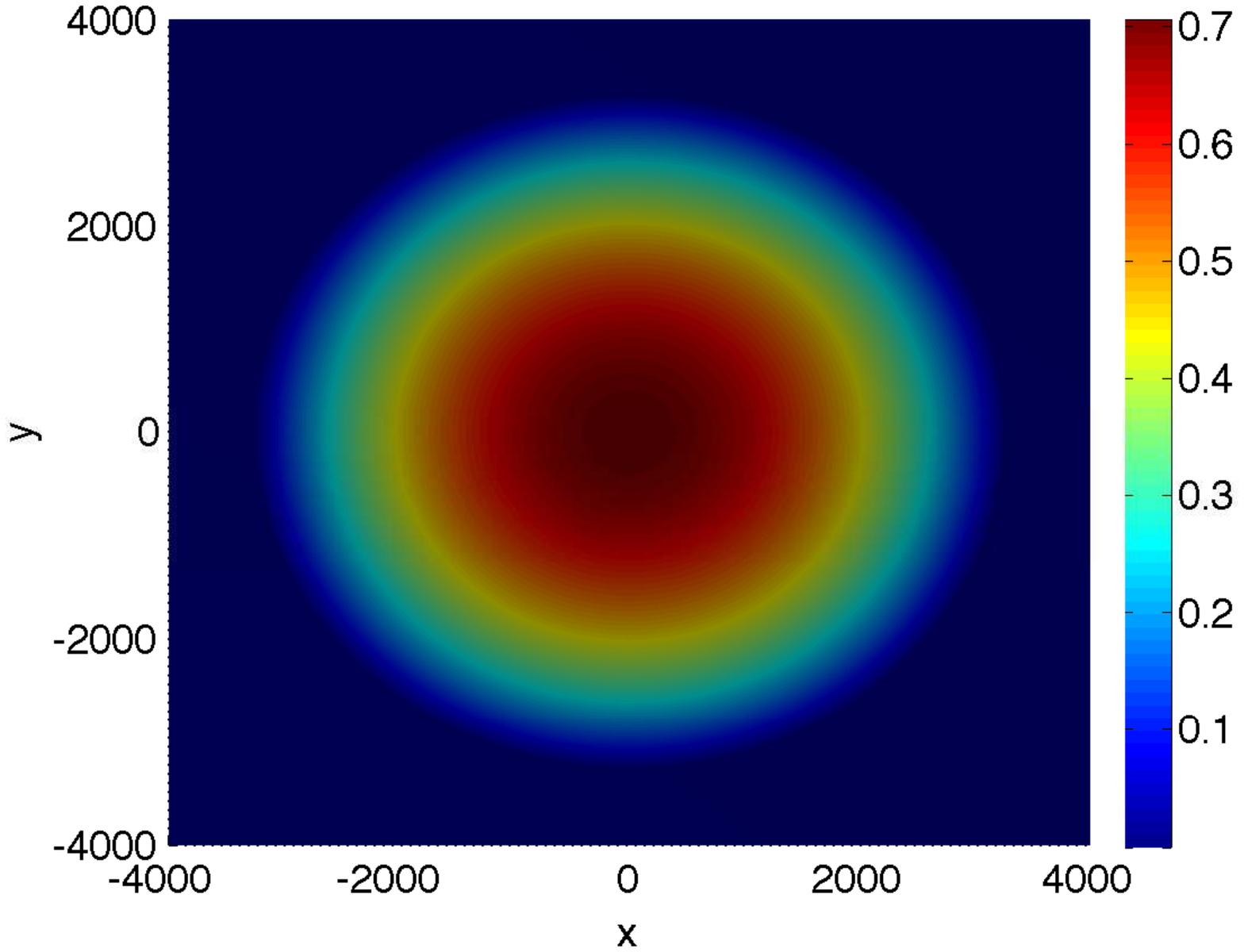}
  \end{minipage}}
  \hspace{0.5cm}
  \subfloat[H=40m, N=2, error]{%
    \begin{minipage}[c]{0.45\linewidth}
      \centering%
      \includegraphics[trim=0cm 6cm 1cm 6cm,clip=true,width=\textwidth]{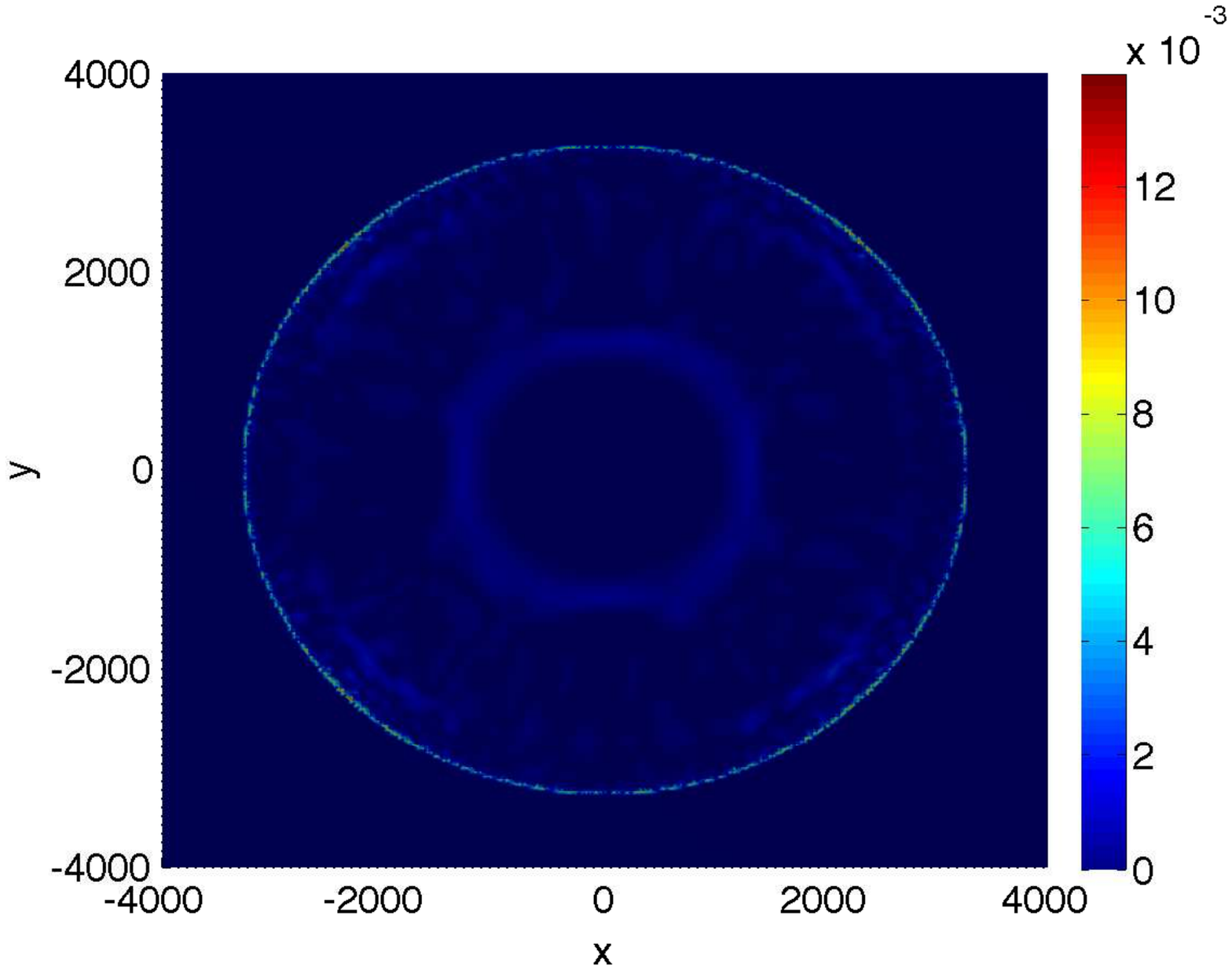}
  \end{minipage}}
  \caption{Parabolic bowl test case : example numerical results at time $t = T/2$}
  \label{fig:parabolic_bowl_errors}
\end{center}
\end{figure}

Here, the constants   are $\alpha = 1.6 \times 10^{-7} m^{-}$, $X = 1 m^{-}$, and $Y = -0.41884 m^{-}$. The solutions are obtained in the square domain with side length $8000m$ centered at the origin. The solution is  continuous, but not continuously differentiable. Fig. (\ref{fig:parabolic_bowl_local_global_L2}) demonstrates the $\mathcal{O}(H^{1.5})$ behavior of the global error in fluid height for approximations with polynomials of order $N=1,2, \text{ and } 3$, while the local errors computed in regions far from the wet/dry front behave like $\mathcal{O}(H^{N+1})$. Fig. (\ref{fig:parabolic_bowl_errors}) indicates that the errors are localized near the wet/dry front and further mesh refinements reduce these local errors.

\subsection{Positivity preserving test : Rarefaction wave}
We use this test case introduced in \cite{ern2008well} to demonstrate the effectiveness of positive preserving limiter. We consider a rectangular domain of $50 \text{m} \times 40 \text{m}$, with a flat bottom. The analytical solution depends on $\xi = \frac{(x-20)}{t}$, and is given by

\begin{equation}
(h, u, v) =
\begin{cases}
(h_{0}, \, 0 , \, 0) & \mbox{if   } \xi < -\sqrt{gh_{0}} \\
(0,     \, 0,  \, 0) & \mbox{if   } \xi > 2\sqrt{gh_{0}} \\
(\frac{1}{9g}(\xi - 2 \sqrt{gh_{0}})^{2},  \, \frac{2}{3}(\xi + \sqrt{gh_{0}}), \, 0) & \mbox{other wise}
\end{cases}
\end{equation}

\begin{figure}[h!]
\begin{center}
  \subfloat[Global  error]{%
    \begin{minipage}[c]{0.45\linewidth}
      \centering%
      \includegraphics[trim=2cm 6cm 1cm 6cm,clip=true,width=\textwidth]{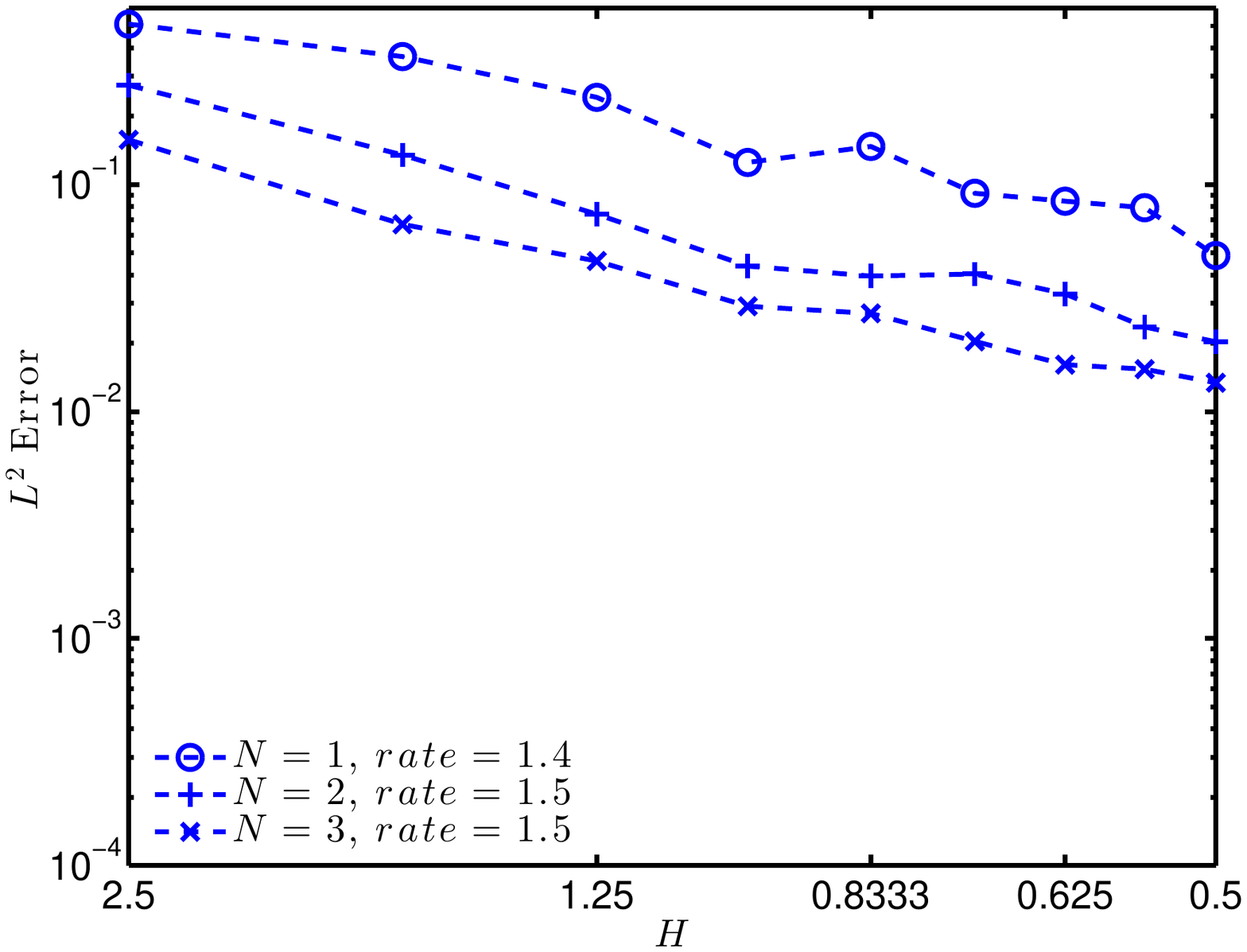}
    \end{minipage}}
  \hspace{0.5cm}
  \subfloat[Local error]{%
    \begin{minipage}[c]{0.45\linewidth}
      \centering%
      \includegraphics[trim=2cm 6cm 1cm 6cm,clip=true,width=\textwidth]{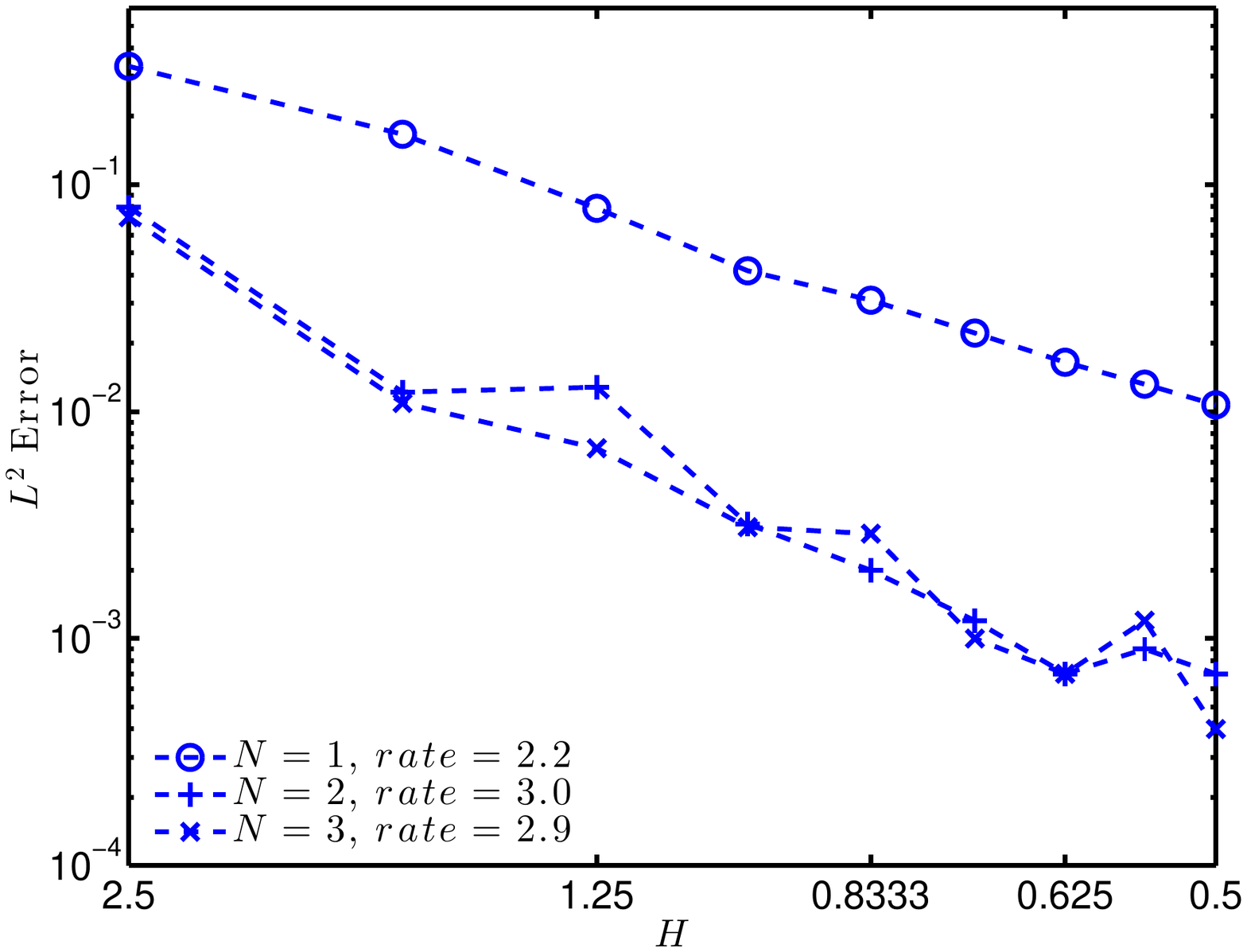}
    \end{minipage}}
\caption{Rarefaction Wave : global and local $L^2$ errors in the solution for fluid height}
\label{fig:rarefaction_global_local_L2}
\end{center}
\end{figure}

We consider $h_0 = 1$ and $g = 1$. The simulations are run until final time $T = 10s$. A CFL number of $0.03$ is used for these simulations in order to minimize transient error and study only  the spatial accuracy. The analytical solution at time $t = 2s$ is used as the initial condition so that the solution is in $C^{0}(\Omega)$. The $L^2$ errors in the  fluid height are presented in Fig. (\ref{fig:rarefaction_global_local_L2}). The TVB\ limiter is not applied for this test case since the solution is smooth enough for a stable computation. We observe that global $L^2$ error in fluid height is $\mathcal{O}(H^{1.4})$, $\mathcal{O}(H^{1.5})$, and $\mathcal{O}(H^{1.5})$ for the polynomial approximations $N=1, \, 2, \, \text{and }3$ respectively.

Since the solution is projected onto linear polynomials, we can expect the global error to behave like $\mathcal{O}(H^{1+1/2})$, which is observed in the estimated order of convergence. To quantify the effect of the positive preserving limiter on the solution accuracy, the $L^2$ errors are measured for the region, $x \in [15, 35]$, where the solution is continuously differentiable. The estimated convergence are $\mathcal{O}(H^{2.2})$, $\mathcal{O}(H^{3.0})$, and $\mathcal{O}(H^{2.9})$ for the polynomials $N=1, \, 2, \, \text{and }3$ respectively (see Fig. (\ref{fig:rarefaction_global_local_L2})). The fluid height is a quadratic polynomial in space and rational polynomial in time, hence increasing the spatial polynomial interpolation beyond $N=2$ does not improve the spatial accuracy which explains the estimated order of convergence for $N=3$.

\begin{figure}[h!]
\begin{center}
  \subfloat[N=1]{%
    \begin{minipage}[c]{0.45\linewidth}
      \centering%
      \includegraphics[trim=0cm 6cm 1cm 6cm,clip=true,width=\textwidth]{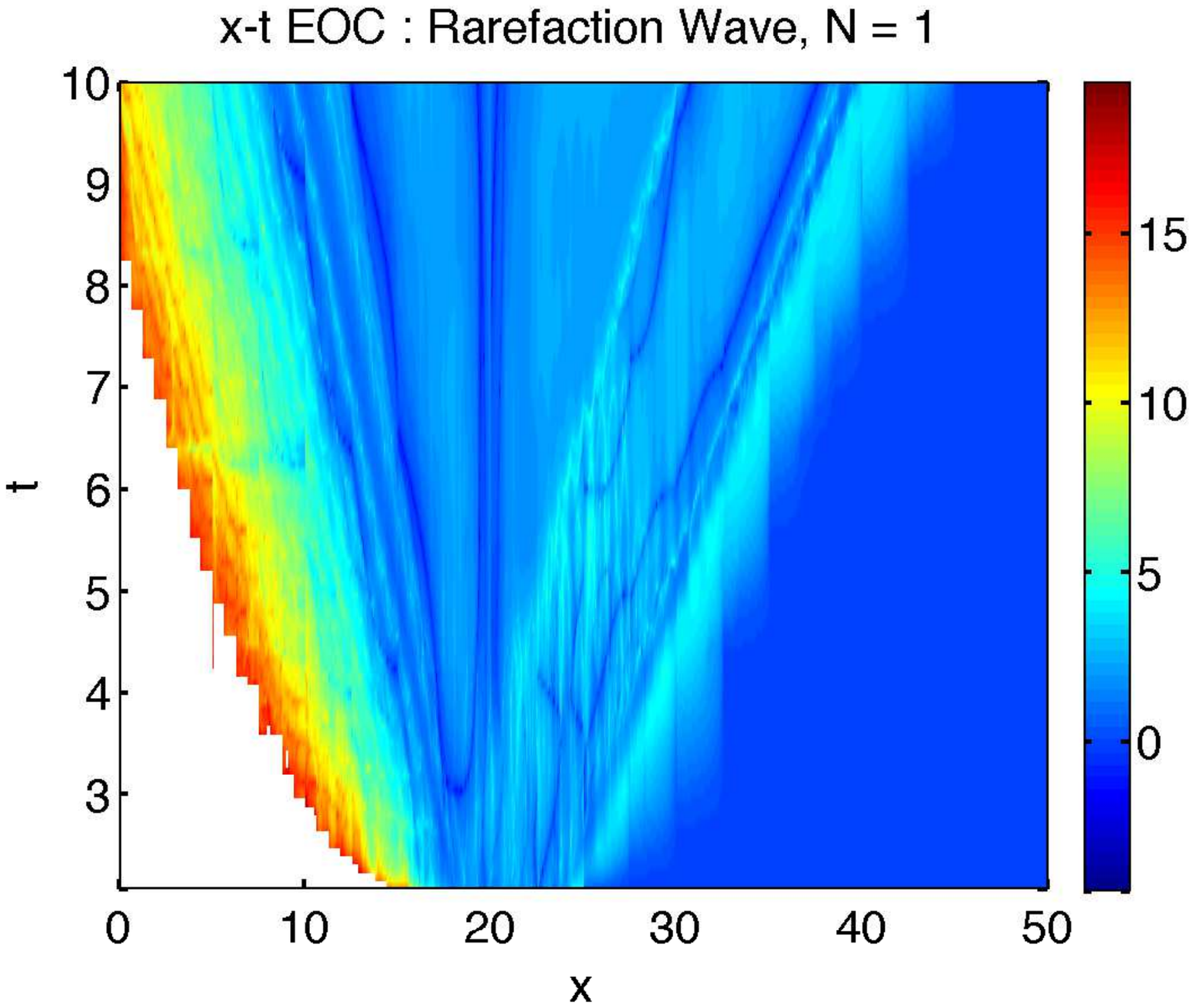}
    \end{minipage}}
  \hspace{0.5cm}
  \subfloat[N=3]{%
    \begin{minipage}[c]{0.45\linewidth}
      \centering%
      \includegraphics[trim=0cm 6cm 1cm 6cm,clip=true,width=\textwidth]{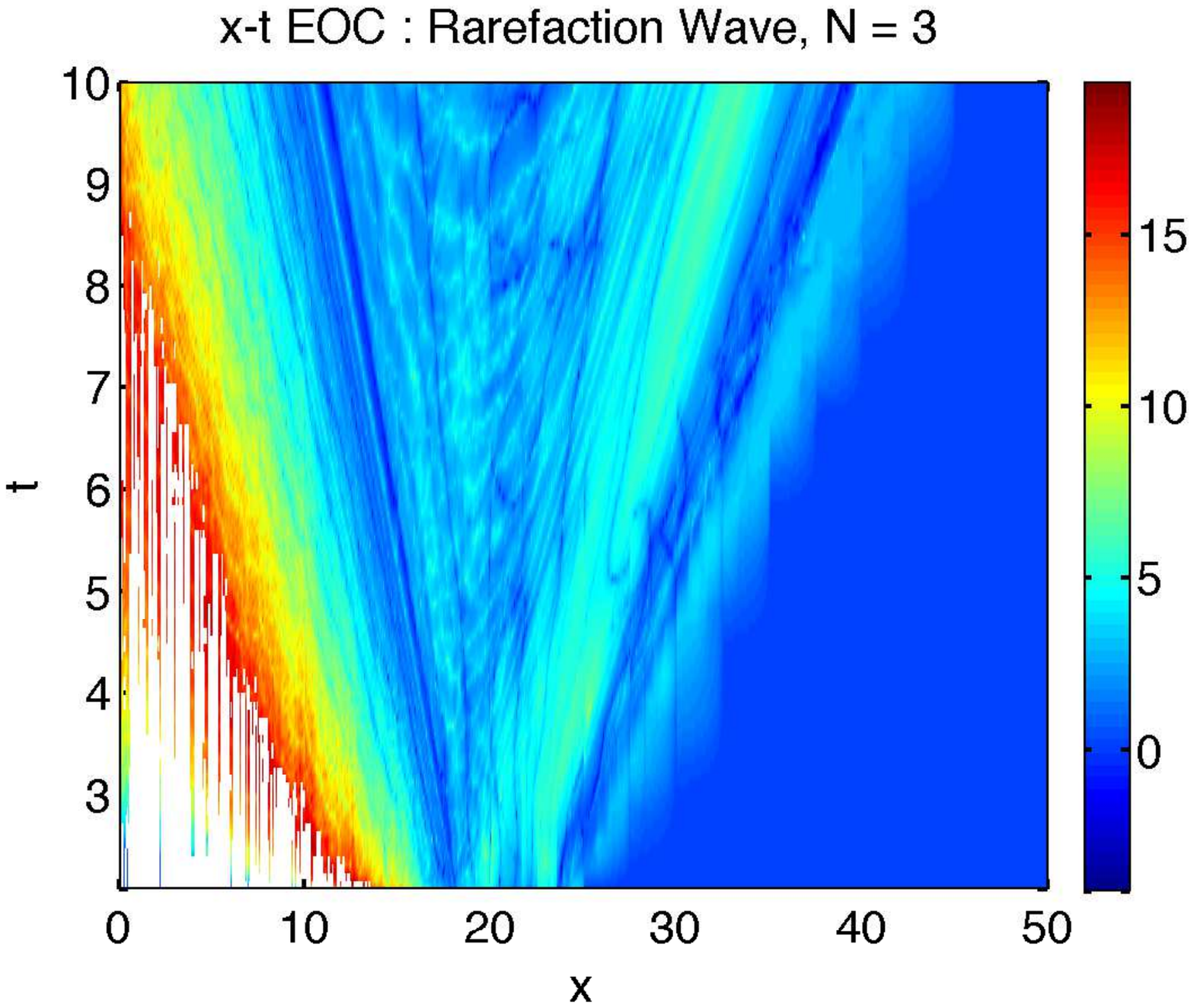}
    \end{minipage}} \\
  \subfloat[Point wise error, N=2 and H = 0.714]{%
    \begin{minipage}[c]{0.45\linewidth}
      \centering%
      \includegraphics[trim=0cm 6cm 1cm 6cm,clip=true,width=\textwidth]{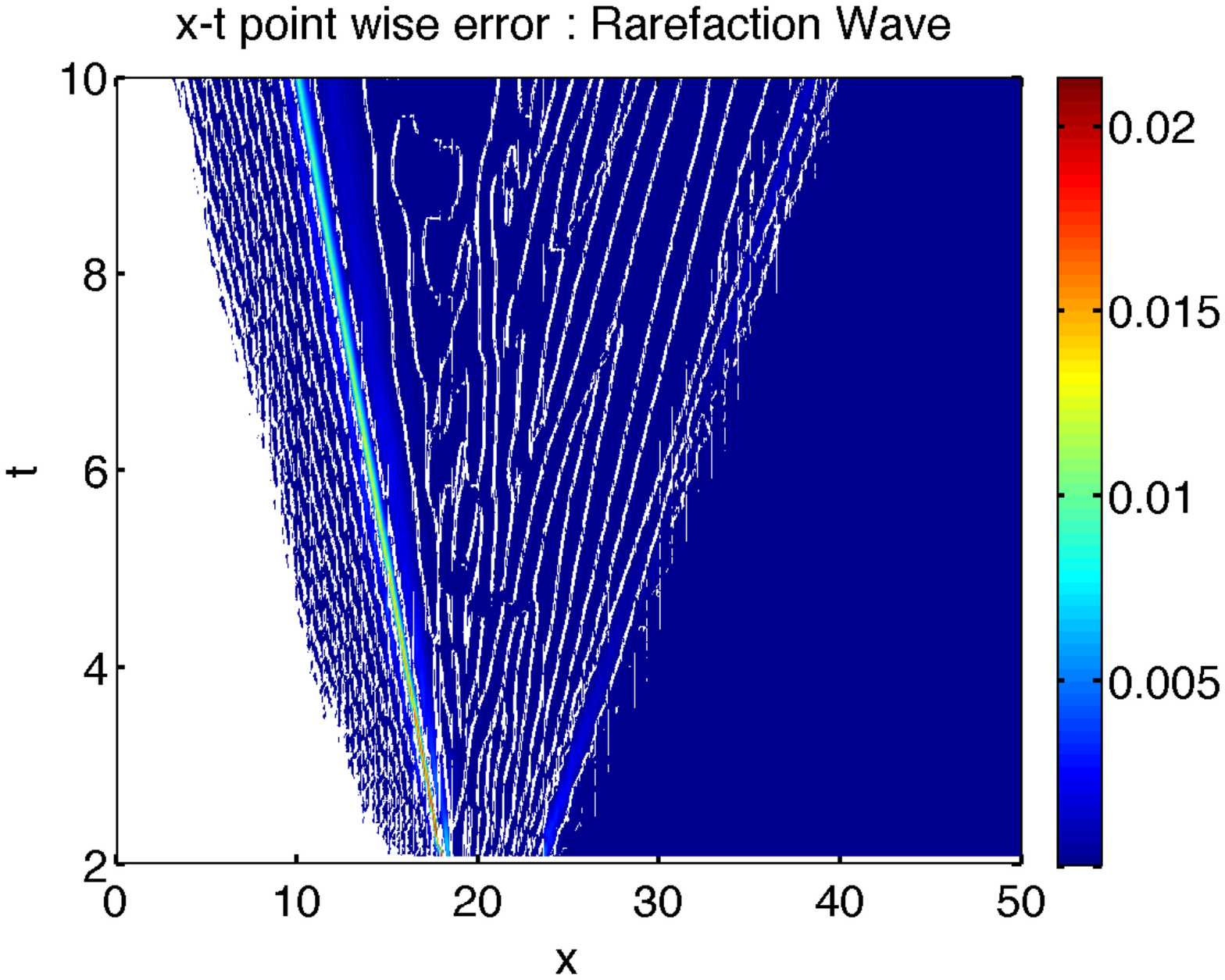}
    \end{minipage}}
\caption{Rarefaction wave : Empirical order of convergence for polynomial orders 1 and 3, and point wise error for polynomial order 2}
\label{fig:rarefaction_xt_EOC}
\end{center}
\end{figure}

There is an increase in error for further refinements. This is due to large oscillations in the numerical solutions, and the slope limiter should be applied to control these oscillations . The point wise EOCs are plotted in the Fig. (\ref{fig:rarefaction_xt_EOC}). It can be observed that the pollution of the solution accuracy is localized to an area of size $\mathcal{O}(T)$, which is the distance traveled by the rarefaction wave in time $T$. The white regions in the plots indicate the error is zero and the light blue color indicate that the error is constant due to the positivity preserving limiter in the dry regions.

\subsection{Limiter test : Two dimensional oscillating lake}
We use this test case proposed in \cite{gallardo2007well, xing2010positivity} to test the effectiveness of the positivity preserving limiter and the modified TVB limiter. We consider a rectangular domain $[-2,2] \otimes [-2,2]$ with a  parabolic bottom topography given by,
\begin{equation}
B(x,y) = h_{0} \frac{x^2 + y^2}{a^2},
\end{equation}
where $h_0$ and $a$ are specified constants. The analytical solution is given by,
\begin{eqnarray}
\label{eq:lake2d}
h(x,y,t) &=& \max \left(0, \frac{\sigma h_0}{a^2} (2x \cos(\omega t) + 2y \sin(\omega t) - \sigma) + h_{0} - b \right), \nonumber \\
u(x,y,t) &=& -\sigma \omega \sin(\omega t), \qquad v(x,y,t) = \sigma \omega \cos(\omega t),
\end{eqnarray}
with frequency $\omega = \sqrt{2gh_0}/a$ and time period $T = 2\pi/\omega$.

\begin{figure}[h!]
\begin{center}
  \subfloat[H=0.0625, N=4, final time = T/2]{%
    \begin{minipage}[c]{0.45\linewidth}
      \centering%
      \includegraphics[trim=0cm 6cm 1cm 6cm,clip=true,width=\textwidth]{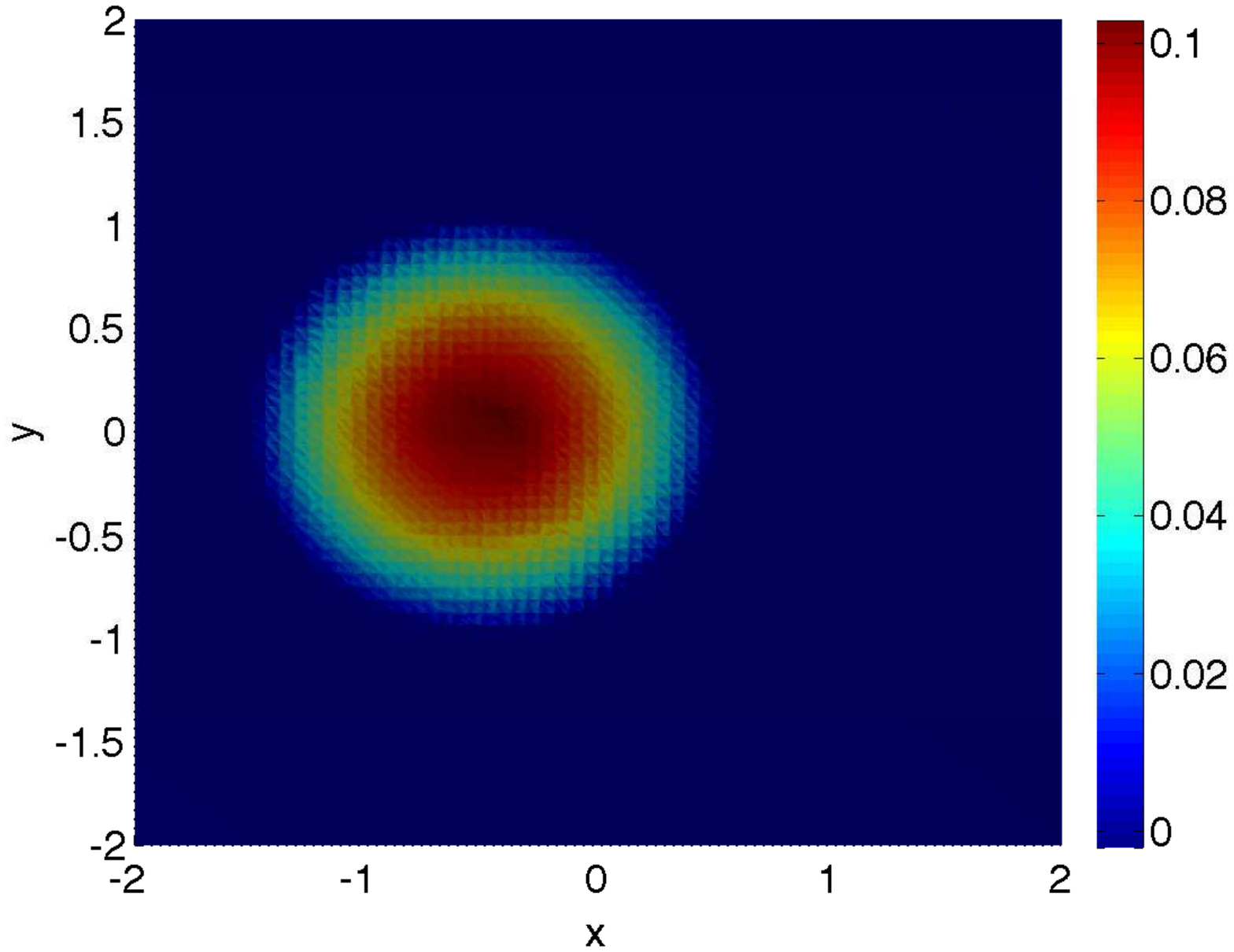}
  \end{minipage}}
  \hspace{0.5cm}
  \subfloat[H=0.03125, N=2, final time = T/2]{%
    \begin{minipage}[c]{0.45\linewidth}
      \centering%
      \includegraphics[trim=0cm 6cm 1cm 6cm,clip=true,width=\textwidth]{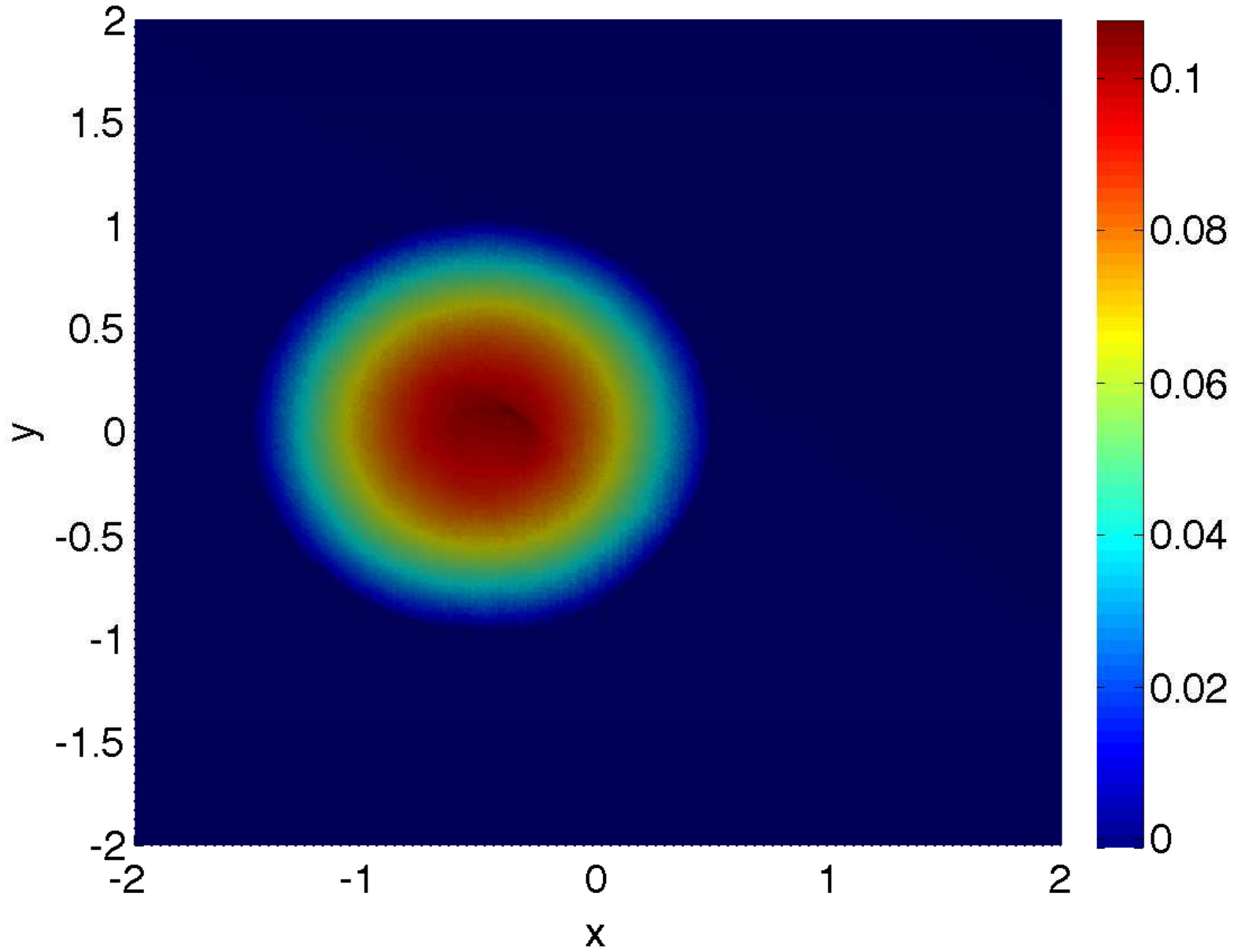}
  \end{minipage}} \\
  \subfloat[H=0.0625, N=4, final time $ = T$]{%
    \begin{minipage}[c]{0.45\linewidth}
      \centering%
      \includegraphics[trim=0cm 6cm 1cm 6cm,clip=true,width=\textwidth]{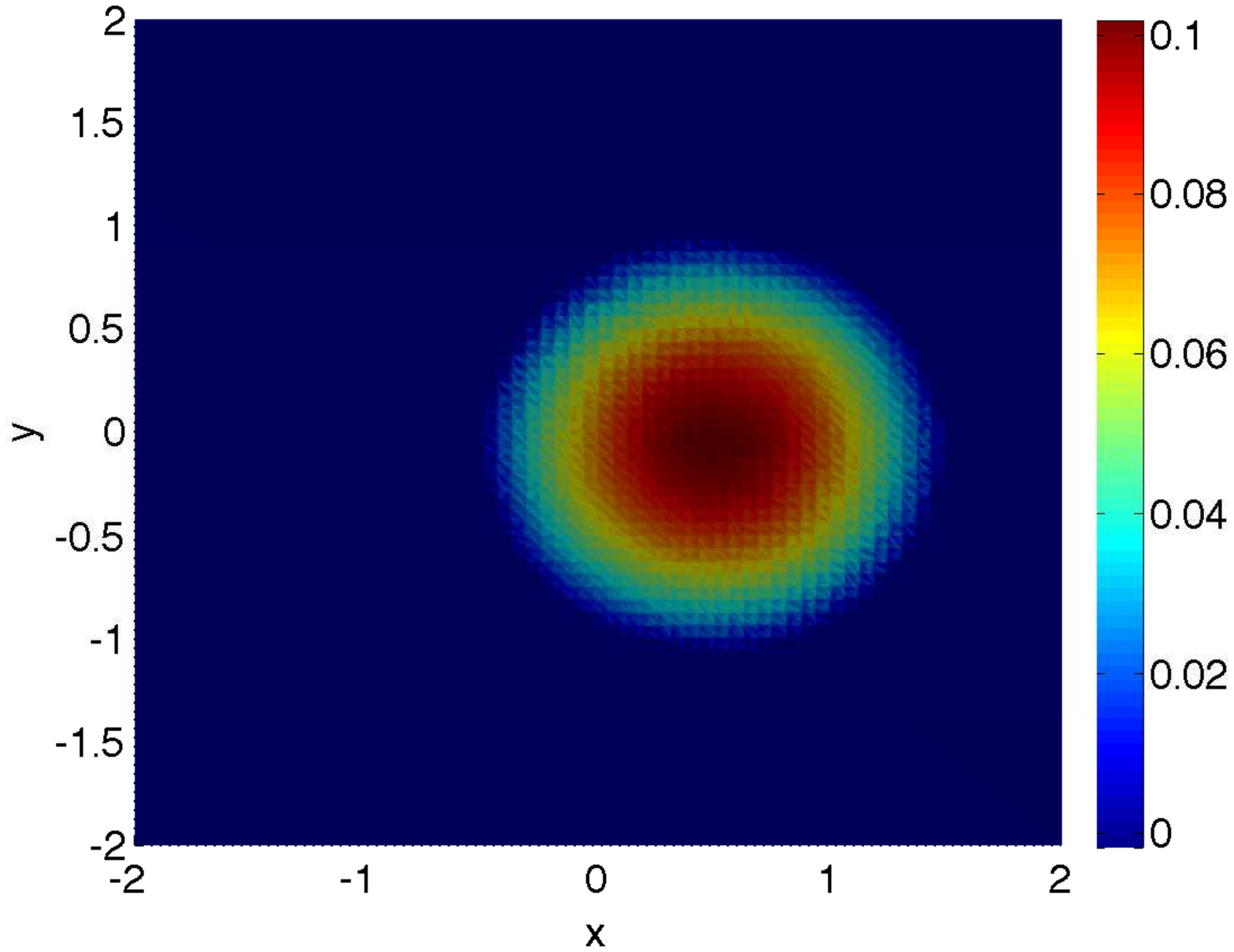}
  \end{minipage}}
  \hspace{0.5cm}
  \subfloat[H=0.03125, N=2, final time $ = T$]{%
    \begin{minipage}[c]{0.45\linewidth}
      \centering%
      \includegraphics[trim=0cm 6cm 1cm 6cm,clip=true,width=\textwidth]{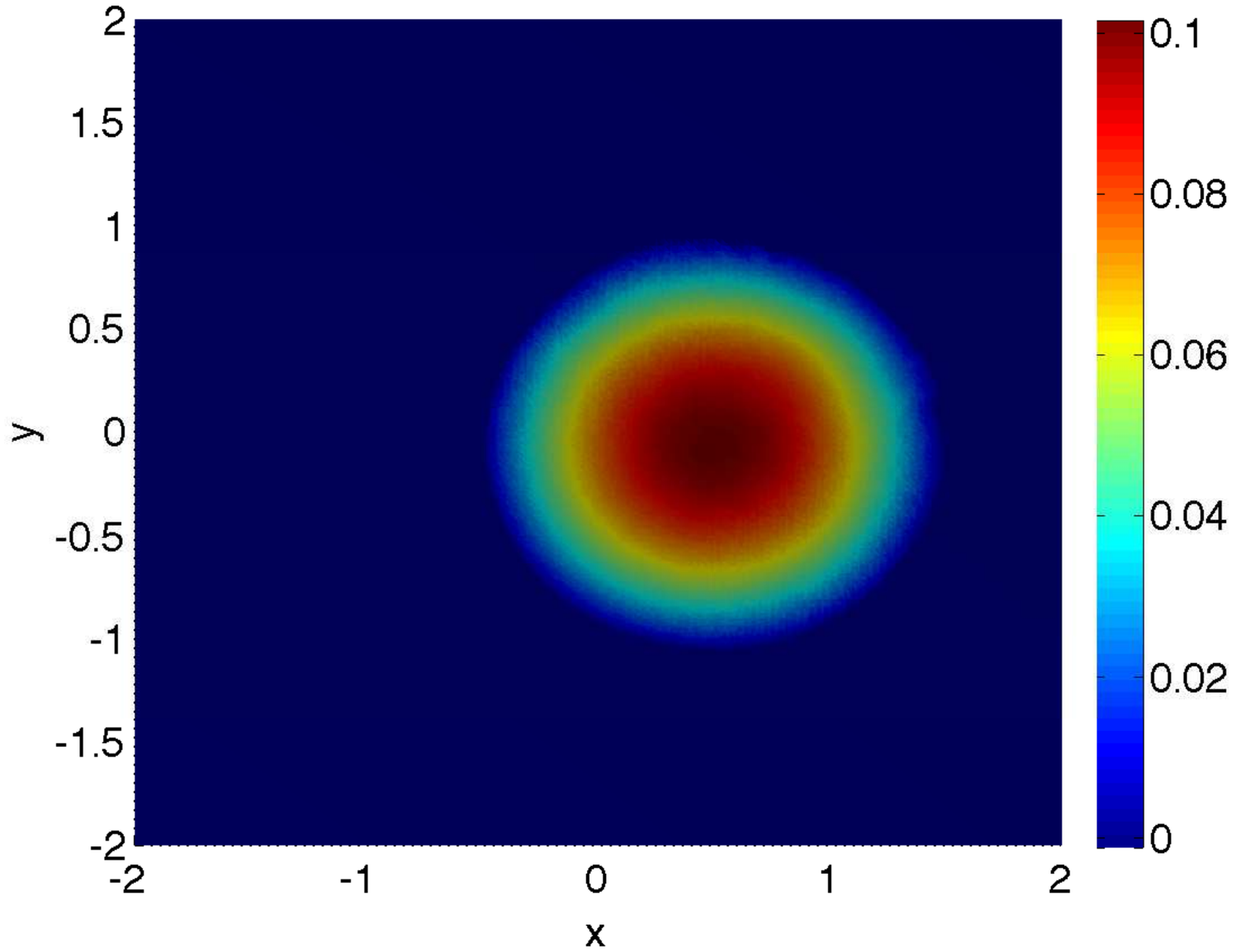}
  \end{minipage}}
  \caption{Oscillating lake 2D : example numerical solutions for fluid height}
  \label{fig:oscillating_lake2D}
\end{center}
\end{figure}

\begin{figure}[h!]
\begin{center}
  \subfloat[H=0.0625, N=4, final time $t = T/2$]{%
    \begin{minipage}[c]{0.45\linewidth}
      \centering%
      \includegraphics[trim=2cm 6cm 1cm 6cm,clip=true,width=\textwidth]{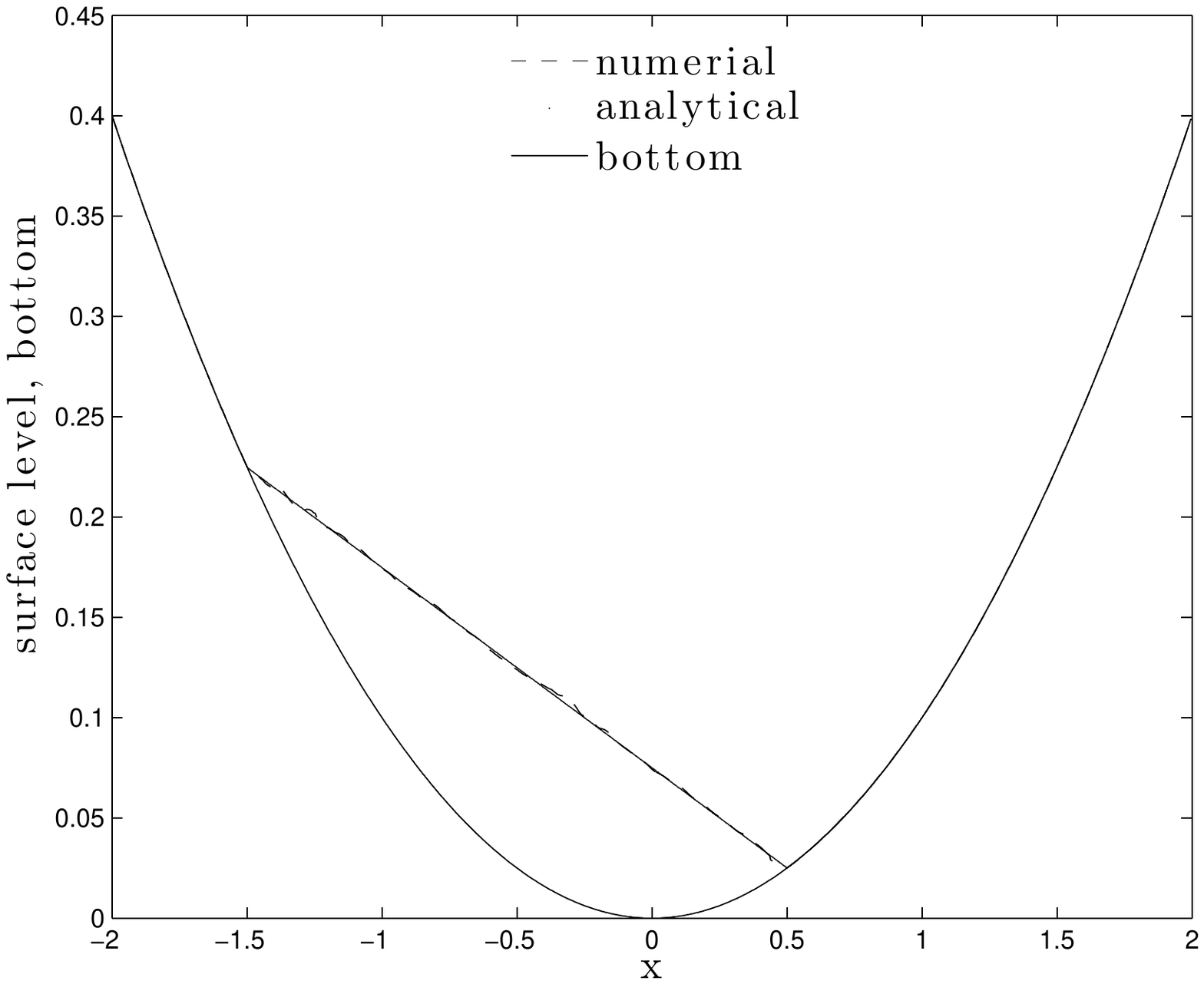}
  \end{minipage}}
  \hspace{0.5cm}
  \subfloat[H=0.03125, N=2, final time $t = T/2$]{%
    \begin{minipage}[c]{0.45\linewidth}
      \centering%
      \includegraphics[trim=2cm 6cm 1cm 6cm,clip=true,width=\textwidth]{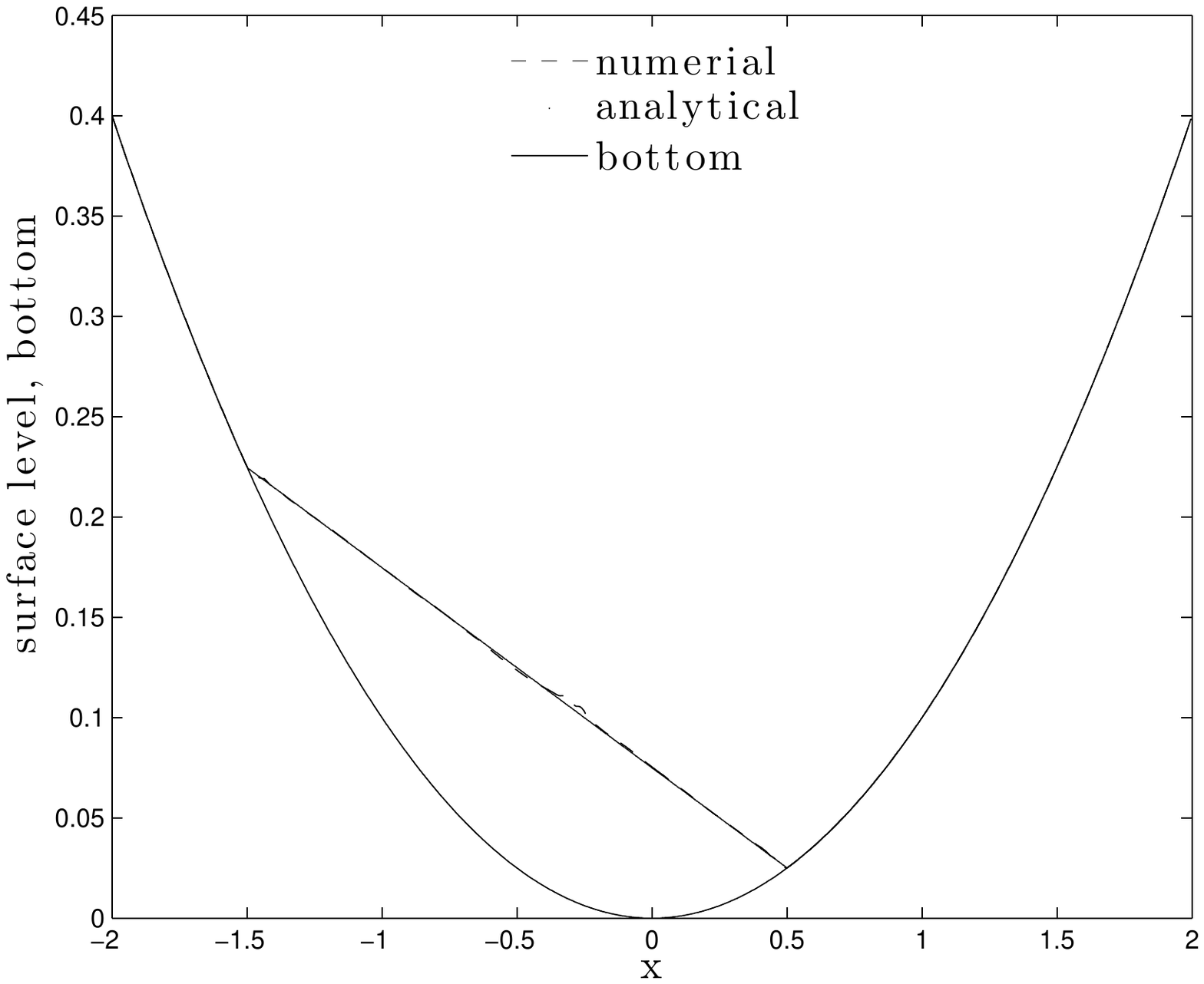}
  \end{minipage}} \\
  \subfloat[H=0.0625, N=4, final time $t = T$]{%
    \begin{minipage}[c]{0.45\linewidth}
      \centering%
      \includegraphics[trim=2cm 6cm 1cm 6cm,clip=true,width=\textwidth]{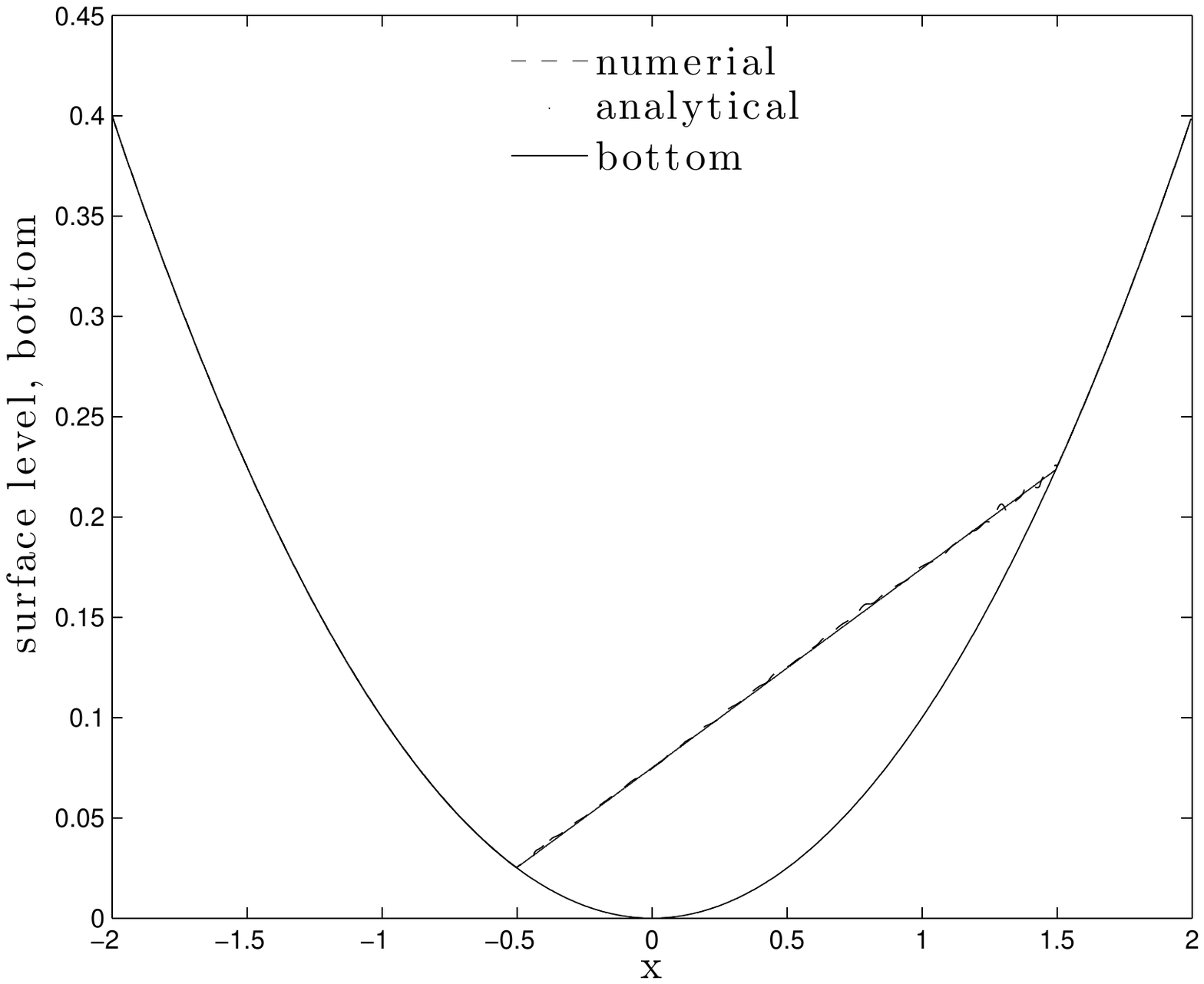}
  \end{minipage}}
  \hspace{0.5cm}
  \subfloat[H=0.03125, N=2, final time $t = T$]{%
    \begin{minipage}[c]{0.45\linewidth}
      \centering%
      \includegraphics[trim=2cm 6cm 1cm 6cm,clip=true,width=\textwidth]{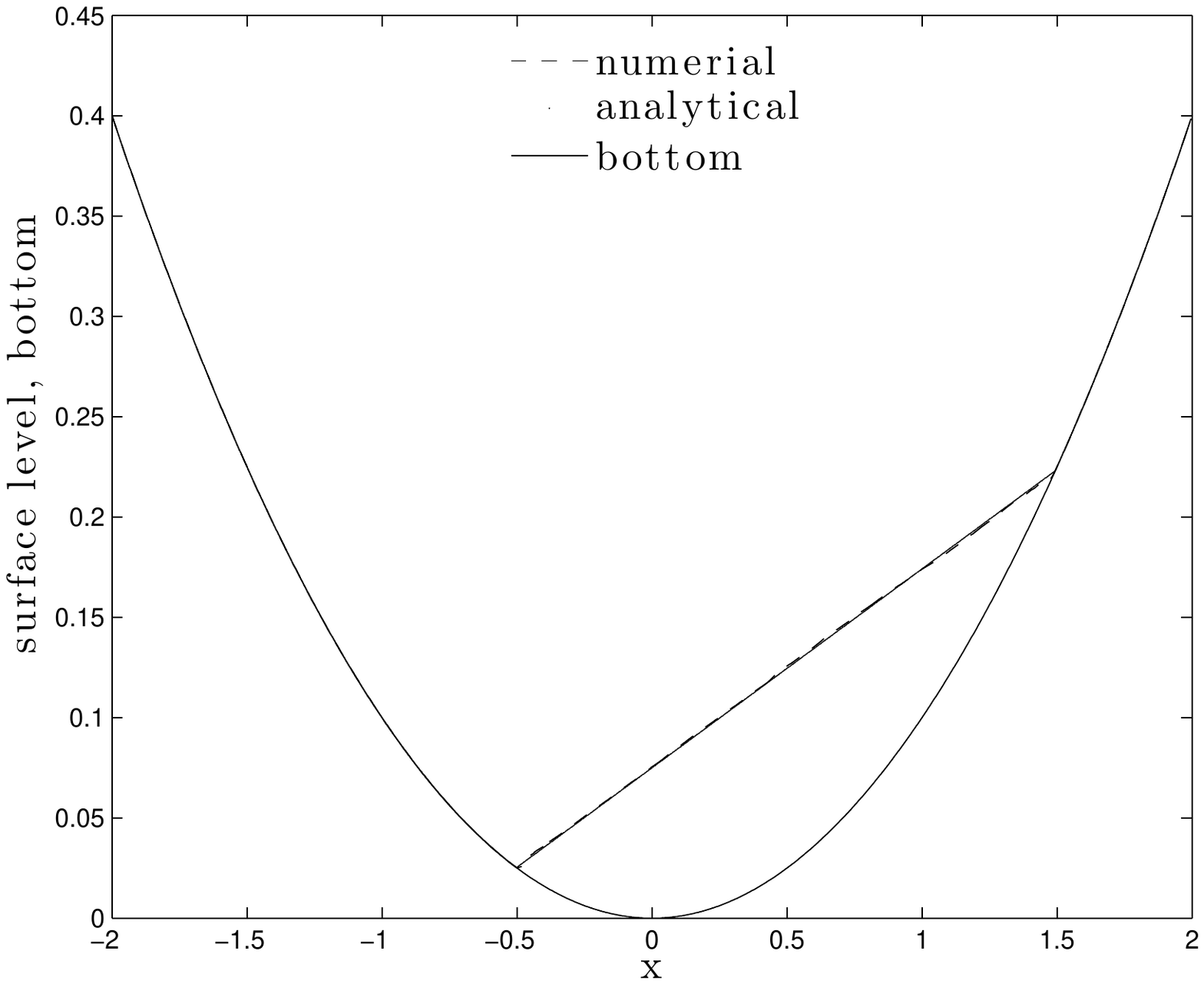}
  \end{minipage}}
  \caption{Oscillating lake 2D : Comparison of numerical solution with analytical solution along the line $y = 0$}
  \label{fig:oscillating_lake2D_y0}
\end{center}
\end{figure}

\begin{figure}[h!]
\begin{center}
  \subfloat[N=1]{%
    \begin{minipage}[c]{0.45\linewidth}
      \centering%
      \includegraphics[width=\textwidth]{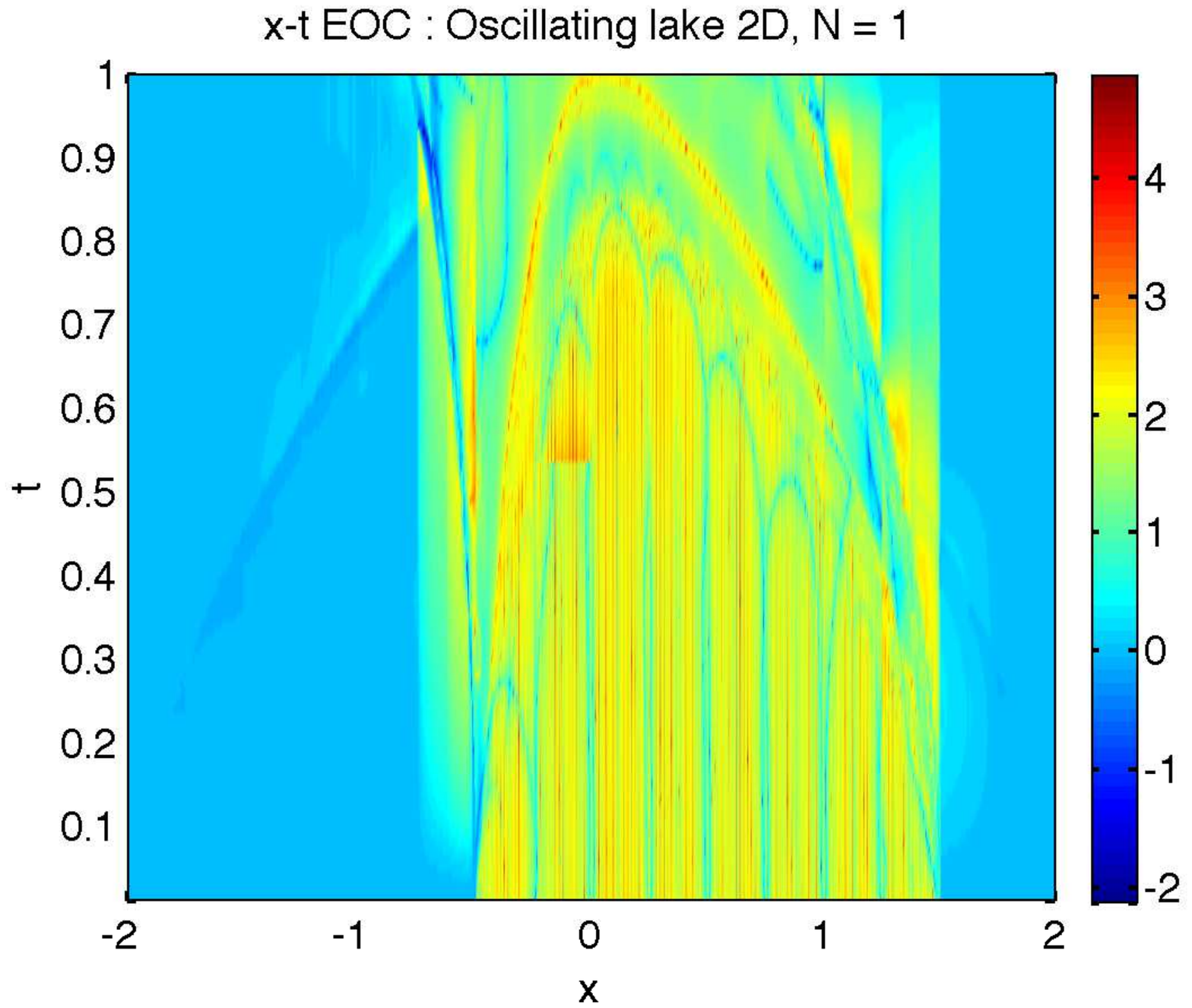}
  \end{minipage}}
  \hspace{0.5cm}
  \subfloat[N=2]{%
    \begin{minipage}[c]{0.45\linewidth}
      \centering%
      \includegraphics[width=\textwidth]{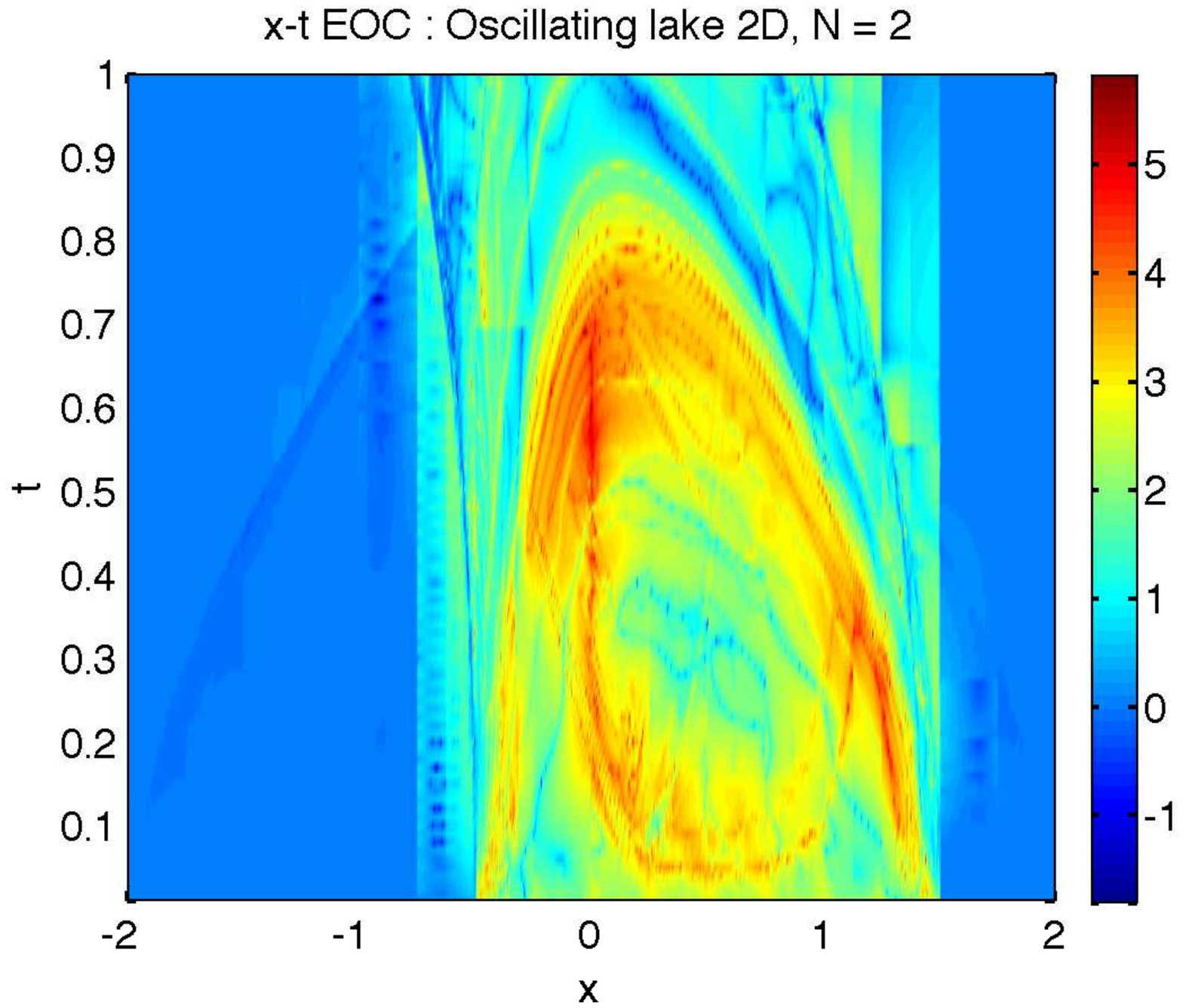}
  \end{minipage}}
\caption{$x-t$ Empirical order of convergence of the numerical solutions. The errors are computed along the line segment $y=0$ for $t \in [0, 1]$.}
\label{fig:oscillating_lake_EOC_xt}
\end{center}
\end{figure}

The constants considered are $a = 1$, $\sigma = 0.5$, and $h_0 = 0.1$. The initial conditions are defined by Eq. (\ref{eq:lake2d}) with $t=0$. Reflecting boundary conditions are used for all the boundaries. Simulations are run until time $T$ with various uniform meshes ($H=0.25, \, 0.125, \, 0.0625, \, \text{and } 0.03125$).
The numerical solutions for $H=0.0625, \, \text{and } 0.03125$, with polynomials of order $3 \, \text{and } 2$ respectively are plotted in the Fig. (\ref{fig:oscillating_lake2D}). The solutions on the line $y=0$ for these simulations are compared with analytical solutions in Fig. (\ref{fig:oscillating_lake2D_y0}). For a fine mesh, the numerical solutions are in good agreement with the analytical solutions. The wet-dry front pollutes the accuracy of the solution. To study the convergence properties of spatial errors, the point wise errors are computed at very small time $t \in [0,1]$. Fig. (\ref{fig:oscillating_lake_EOC_xt}) shows the improvement of point wise EOCs in the solution for $N=2$, compared to $N=1$. These point wise errors are computed along the line $y=0$.

\subsection{Summary of the test cases }
We summarize all the test cases in Table (\ref{tab:summary_test_cases}). For each test case, we describe the regularity of the solution, the limiters used for the simulation, global $L^2$ error estimates, and local $L^2$ estimates for the solutions with insufficient regularity.

\begin{table}[h!]
\begin{center}
\begin{tabular}{llccc}
\hline
Test case              & Solution regularity           & Limiters      & Global Error              & Local Error \\ \hline
Couette flow           & 2D, $C^{\infty}(\Omega)$      & -             & $\mathcal{O}(H^{N+1/2})$  & - \\
Isentropic vortex      & 2D, $C^{\infty}(\Omega, T)$   & -             & $\mathcal{O}(H^{N+1/2})$  & - \\
Parabolic bowl        & 2D, $C^{0}(\Omega, T)$        & PP            & $\mathcal{O}(H^{1+1/2})$  & $\mathcal{O}(H^{N+1/2})$  \\
Rarefaction wave       & 1D, $C^{0}(\Omega, T)$        & PP            & $\mathcal{O}(H^{1+1/2})$  & $\mathcal{O}(H^{N+1})$  \\
Oscillating lake       & 2D, $C^{0}(\Omega, T)$        & PP, TVB       & -                         & - \\ \hline
\end{tabular}
\caption{Summary of the test cases}
\label{tab:summary_test_cases}
\end{center}
\end{table}



\vspace{5cm}

We demonstrated that the rate of decay in the error matches the predicted rate for sufficiently smooth solutions, while accuracy is lost for the problems with insufficient regularity. For these problems, we illustrated that the low EOCs are confined to regions of wave front.\  To demonstrate the impact of the limiters on the solution accuracy, we estimated the point-wise empirical order of convergence for the rarefaction wave test case and observed that the pollution in the accuracy is localized to the regions the wave front. The localized errors can be controlled by adaptively refining the meshes near the irregularities in the solution, however the mesh refinement techniques are beyond the focus of this paper.


\section{GPU Acceleration}
\label{sec:gpu}
This Section describes a mapping of the nodal DG discretization onto  the wide SIMD model of GPUs using OpenCL. See \cite{wen2011gpu, klockner2010high} for a detailed explanation of OpenCL implementation of DG for electromagnetic applications.

There are three major computations; volume integration, surface integration, and time step update. OpenCL kernels \verb=VolumeKernel=, \verb=SurfaceKernel=, and \verb=UpdateKernel= perform these computations respectively.
\begin{equation}
\label{eq:kernels}
\underbrace{\frac{d Q_{H}}{dt}}_{\verb=UpdateKernel=} = \underbrace{\mathcal{N}(Q_{H})}_{\verb=VolumeKernel=} \, \underbrace{ + \, \, \mathcal{S}(Q_{H}^{g,+},Q_{H}^{g,-})}_{\verb=SurfaceKernel=}.
\end{equation}

\begin{figure}[h!]
\begin{center}
\includegraphics[trim=0cm 6cm 1cm 6cm,clip=true,width=0.6\textwidth]{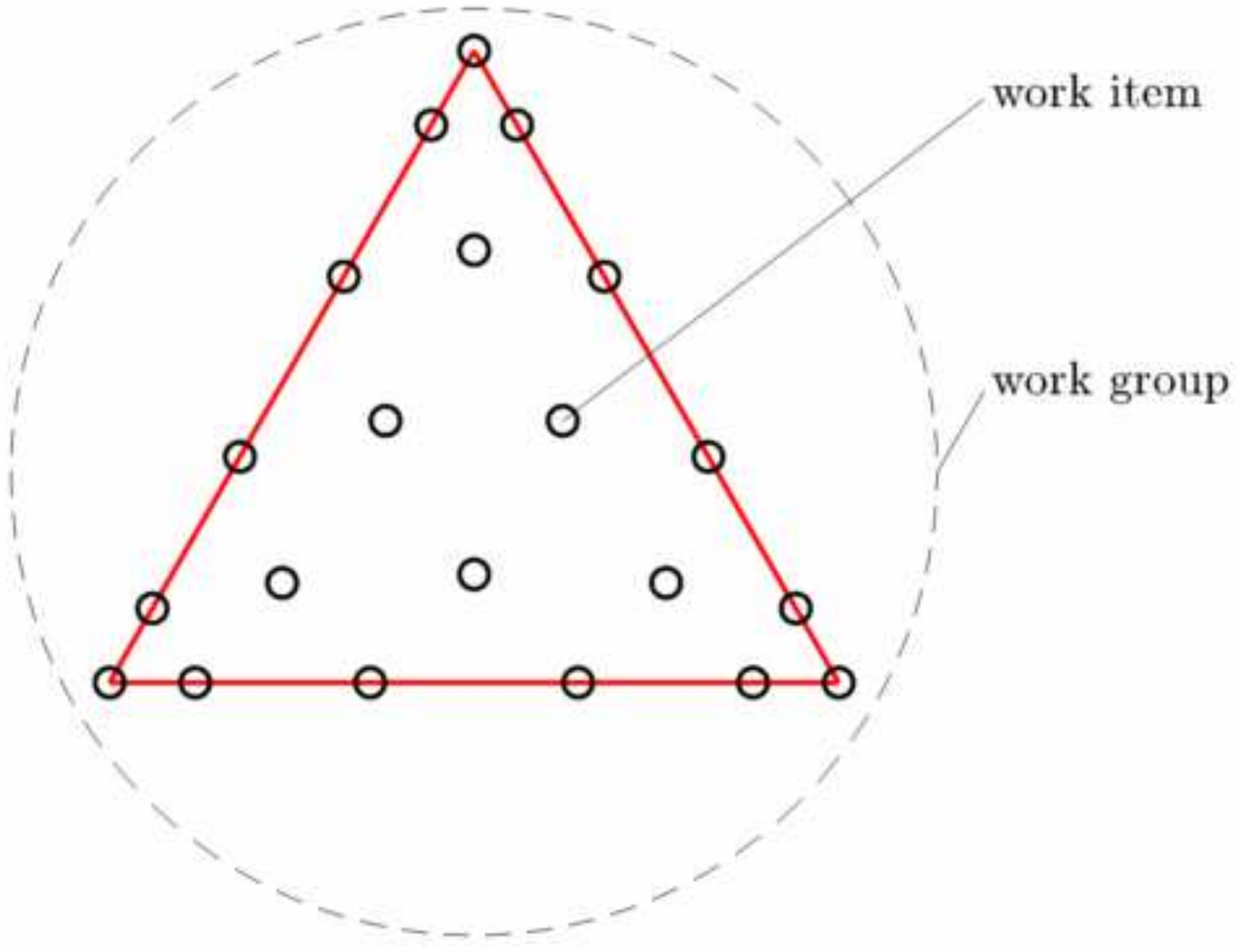}
\caption{Each work item updates the values at an interpolation node, and a work group updates the nodal values in one or more elements.}
\label{fig:work_group}
\end{center}
\end{figure}

An OpenCL work group  computes the integrals of  one or more elements, while one work item computes the contribution from each integration node (see Fig. (\ref{fig:work_group})) in these kernels. For each element, contributions from volume and surface integrals are represented as matrix-vector products of a dense rectangular matrix  and a vector of field values or fluxes. Each work item in a work group computes one entry of the matrix-vector product for every field to avoid the memory conflicts. Here, we explain the implementation  and discuss the performance tuning of these kernels.


\begin{table}[h!]
\begin{center}
\begin{tabular}{lcl}
\hline
Symbol          &   & Definition   \\ \hline
\verb=Nfaces=   &:  & Number of interfaces per triangle ($3$) \\
\verb=N=        &:  & Polynomial order \\
\verb=Np=       &:  & Number of interpolation points \verb=(N+1)*(N+2)/2= \\
\verb=Nfp=      &:  & Number of interpolation points on an interface \verb=N+1=\\
\verb=Ncub=     &:  & Number of cubature nodes on a triangle \\
\verb=Ngauss=   &:  & Number of Gauss integration nodes on a line segment \\\hline
\end{tabular}
\label{tab:notation}
\caption{Notation}
\end{center}
\end{table}

\subsection{Volume Kernel}
The contributions from volume integrals ($\mathcal{N}(Q_{H})$) are computed in this kernel.  $\mathcal{N}(Q_{H})$ is a vector of length \verb=Np= per field for each element $D^{k}$. $\mathcal{N}(Q_{H})_i$ represents $i^{th}$ entry of the vector. The computations are,
\begin{eqnarray}
\mathcal{N}(Q_{H})&=& Pr \times cF1 + Ps \times cF2 + P \times cS,
\end{eqnarray}

\noindent where $cF1$, $cF2$ are numerical flux components in $r-$, $s-$ directions and $cS$ is a vector of source terms at the cubature integration nodes. $P$, $Pr$ and $Ps$ are projection matrices that are pre multiplied with cubature integration weights and depend only on the reference triangle. 

\begin{itemize}
\item \textbf{Interpolation to cubature nodes:} Each work item computes all the field values ($cQ$) at a cubature node. This involves three matrix vector products, all with the same interpolation matrix.  The resulting field values are stored in register memory. \verb=Ncub= work items are assigned for this operation.

\item \textbf{Volume flux evaluation:} Each work item computes all the volume fluxes ($cF1$ and $cF2$) and source terms ($cS$) at a cubature node using the cubature field values stored in register memory, and stores them in a shared memory array for the fluxes. \verb=Ncub= work items are assigned for this operation.

\item \textbf{Projection to interpolation nodes:} Each work item computes the contribution from volume integrals to time derivatives at an interpolation node. This involves eight matrix vector products (three for $Pr \times cF1$, three for $Ps \times cF2$ and two for $P \times cS$) with three projection matrices ($Pr$, $Ps$ and $P$) and three flux vectors that were stored in shared memory in flux evaluation. \verb=Np= work items are assigned for
this operation.
\end{itemize}
 To accommodate the number of work items required for all the computations, a maximum of \verb=Np=  and \verb=Ncub=  work items are requested per element. This requires total number of \verb=Kv*max(Ncub, Np)= work items per work group, where \verb=Kv= is the number of elements processed by a work group.

\subsection{Surface Kernel}
The contributions from surface integrals ( $\mathcal{S}(Q_{H}^{g,+}, Q_{H}^{g,-})$) are computed in this kernel. $\mathcal{S}(Q_{H}^{g,+}, Q_{H}^{g,-})$ is a vector of length \verb=Np= per field for each element $D^{k}$. $\mathcal{S}(Q^{g,+}_H, Q^{g,-}_H)_i$ represents $i^{th}$ entry of the vector. The  computations are,
\begin{eqnarray}
\mathcal{S}(Q^{g,+}_H, Q^{g,-}_H) &=& - L^g \times F_n^{*,g},
\end{eqnarray}

\noindent where $L^g$ is a projection/lifting operator that projects the contribution from the surface integrals to the interpolation nodes. $F_n^{*,g}$ is a vector of stable numerical fluxes at Gauss quadrature nodes. $Q^{g}$, the vector of field values at the Gauss quadrature nodes is computed a priori to this kernel.
\begin{itemize}
\item \textbf{Numerical flux evaluation:} Each work item computes all the numerical fluxes at a Gauss quadrature node ($F_n^{*,g}$). These numerical fluxes are stored in shared memory for further computations. \verb=Ngauss*Nfaces=
 work items are assigned for this operation.



\item \textbf{Lifting the flux to interpolation nodes:} Each work item computes the surface integral contribution at an interpolation node. This involves three matrix vector products with lifting operator ($L^g$). As discussed in volume kernel, the lifting operator is copied to shared memory. \verb=Np= work items are assigned for this operation.




\end{itemize}

To accommodate the number of work items required for all the computations, a maximum of \verb=Np=  and \verb=Ngauss*Nfaces=  work items are used for the computations per element. This requires total number of \verb=Ks*max(Ngauss*Nfaces, Np)= work items per work group, where \verb=Ks= is the number of elements processed by a work group.

\subsection{Update Kernel}
The field values at interpolation nodes ($Q_H$) and Gauss integration nodes ($Q^g_H$) for each element are updated in this kernel. The update kernel is divided in to two parts.
\begin{itemize}
\item \textbf{Update interpolation nodes:} Each work item updates all the field values at an interpolation node. \verb=Np= work items are assigned for this operation.




\item \textbf{Interpolate to Gauss quadrature nodes:} Each work item computes all the field values at a Gauss quadrature node.  This involves three matrix vector products per element, all with a Gauss interpolation matrix. \verb=Ngauss*Nfaces= work items are assigned for this operation.


\end{itemize}

To accommodate the number of work items required for all the computations \verb=Np= work items are requested per element. This requires total number of \verb=Ku*Np= work items per work group. Interpolation to Gauss quadrature nodes is implemented in a separate kernel that requires \verb=Ngauss*Nfaces*Ku= work items per work group. Here, \verb=Ku= is the number of elements processed by a work group.

\subsection{Kernel Tuning}
The performance of the above discussed kernels is very sensitive to the hardware, tuning parameters, optimal usage of shared/local memory. Here, we discuss some of the optimization techniques we adopted and resulted in significant performance improvement.

\begin{figure}
\begin{center}
  \subfloat[Volume kernel on NVIDIA Tesla C2050]{%
    \begin{minipage}[c]{0.45\linewidth}
      \centering%
      \includegraphics[trim=1cm 6cm 1cm 6cm,clip=true,width=1\textwidth]{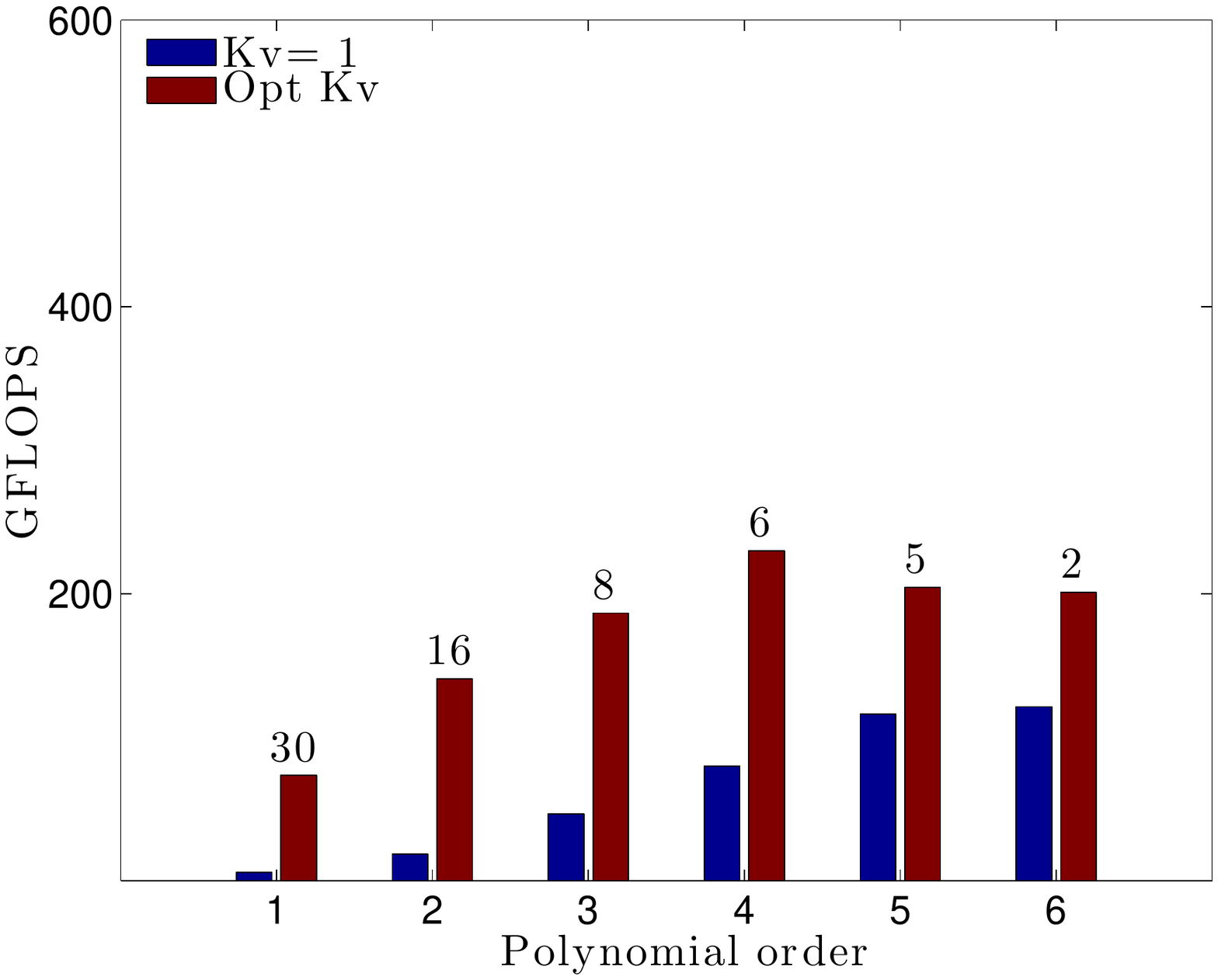}
  \end{minipage}}  \hspace{0.5cm}
  \subfloat[Surface kernel on NVIDIA Tesla C2050]{%
    \begin{minipage}[c]{0.45\linewidth}
      \centering%
      \includegraphics[trim=1cm 6cm 1cm 6cm,clip=true,width=1\textwidth]{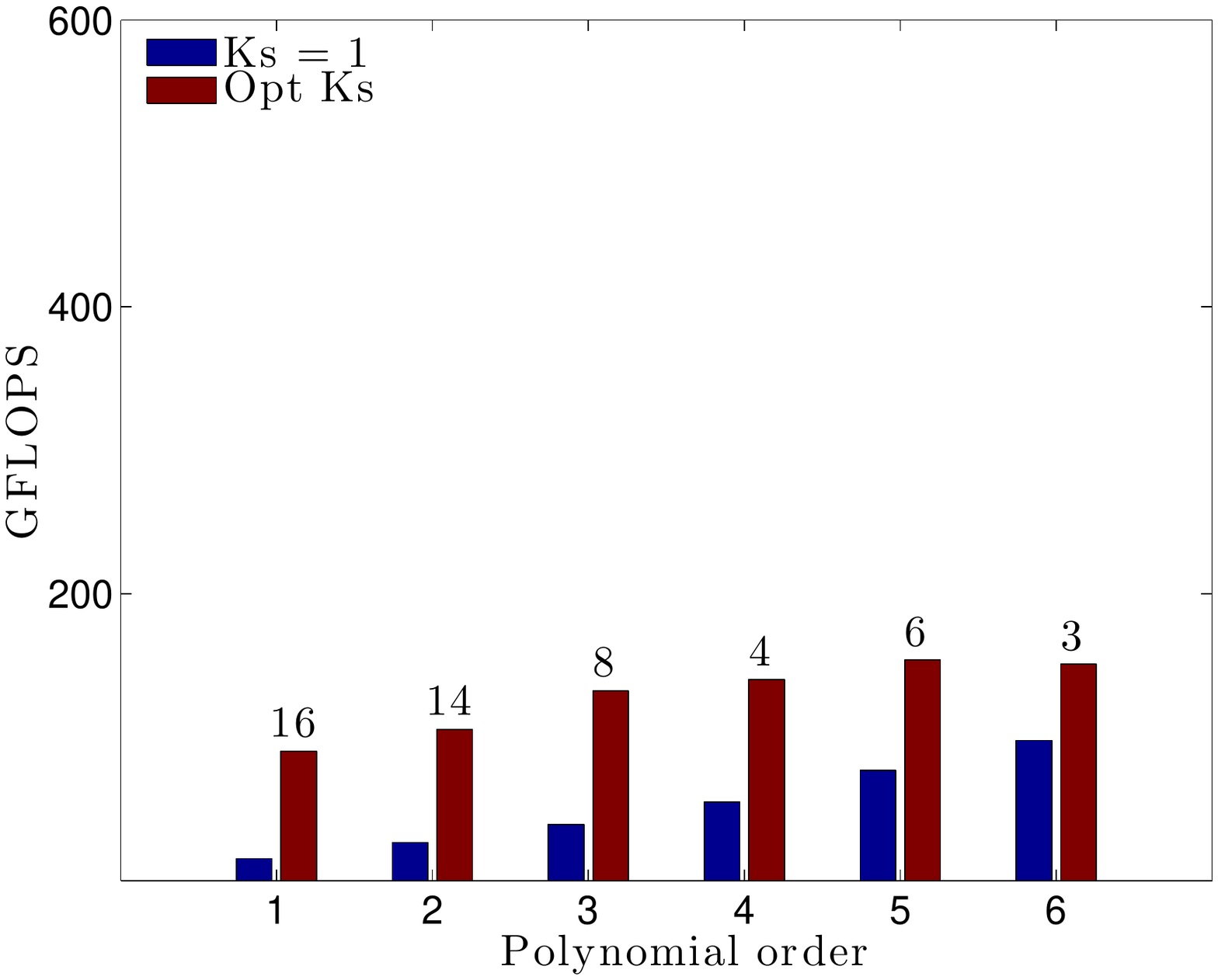}
  \end{minipage}} \\
  \subfloat[Volume kernel on NVIDIA Kepler K20]{%
    \begin{minipage}[c]{0.45\linewidth}
      \centering%
      \includegraphics[trim=1cm 6cm 1cm 6cm,clip=true,width=1\textwidth]{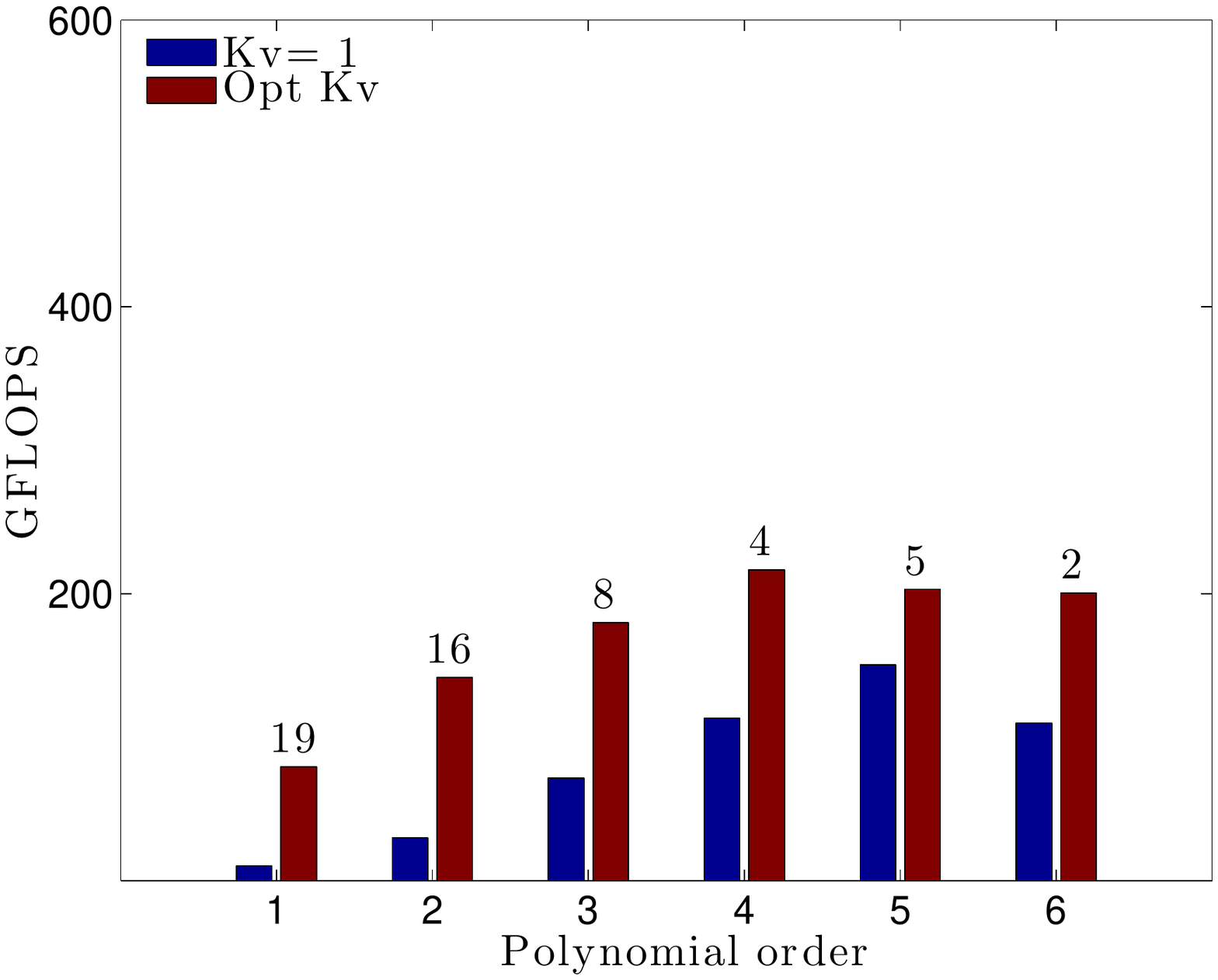}
  \end{minipage}}  \hspace{0.5cm}
  \subfloat[Surface kernel on NVIDIA Kepler K20]{%
    \begin{minipage}[c]{0.45\linewidth}
      \centering%
      \includegraphics[trim=1cm 6cm 1cm 6cm,clip=true,width=1\textwidth]{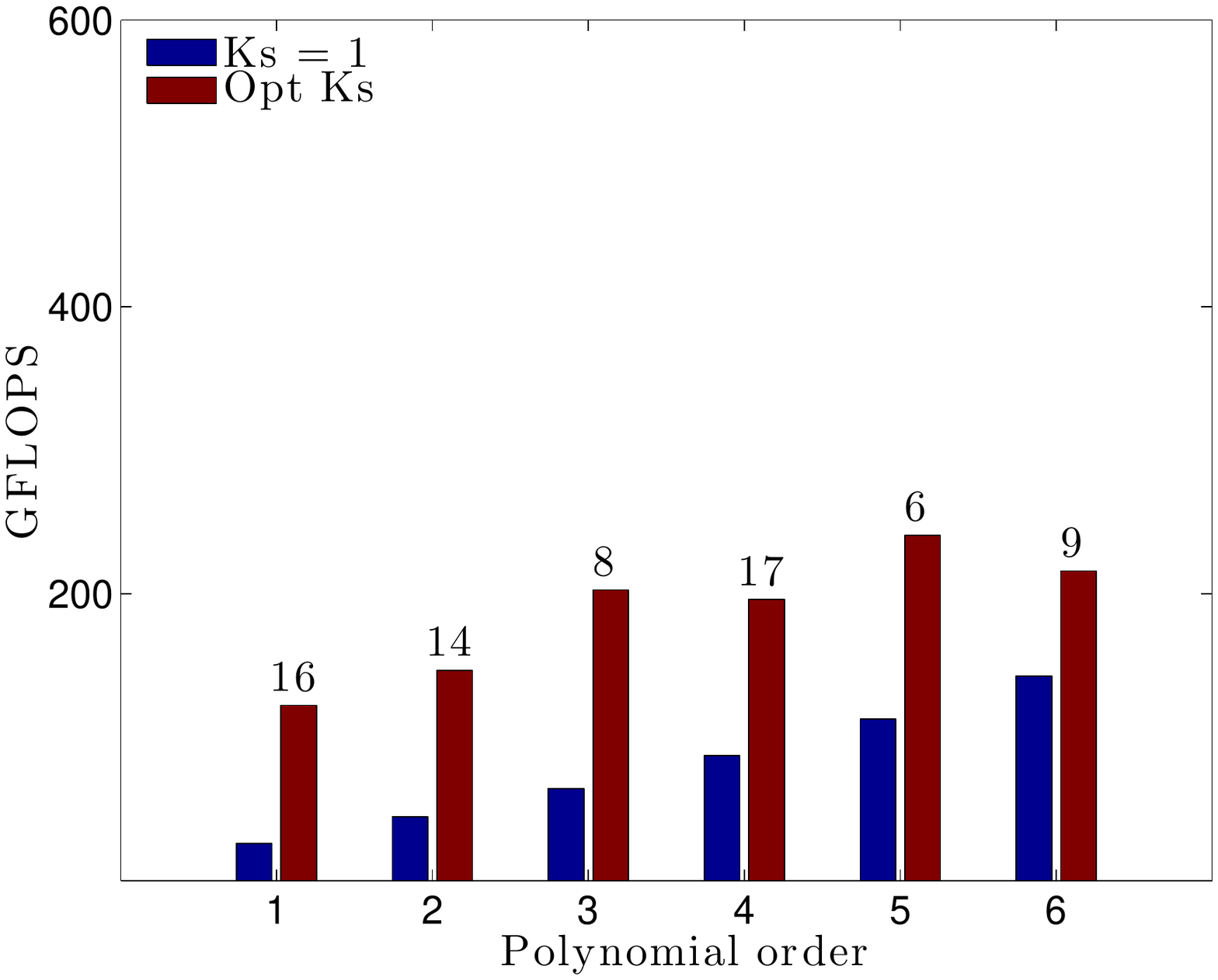}
  \end{minipage}} \\
  \subfloat[Volume kernel on AMD Tahiti 7970]{%
    \begin{minipage}[c]{0.45\linewidth}
      \centering%
      \includegraphics[trim=1cm 6cm 1cm 6cm,clip=true,width=1\textwidth]{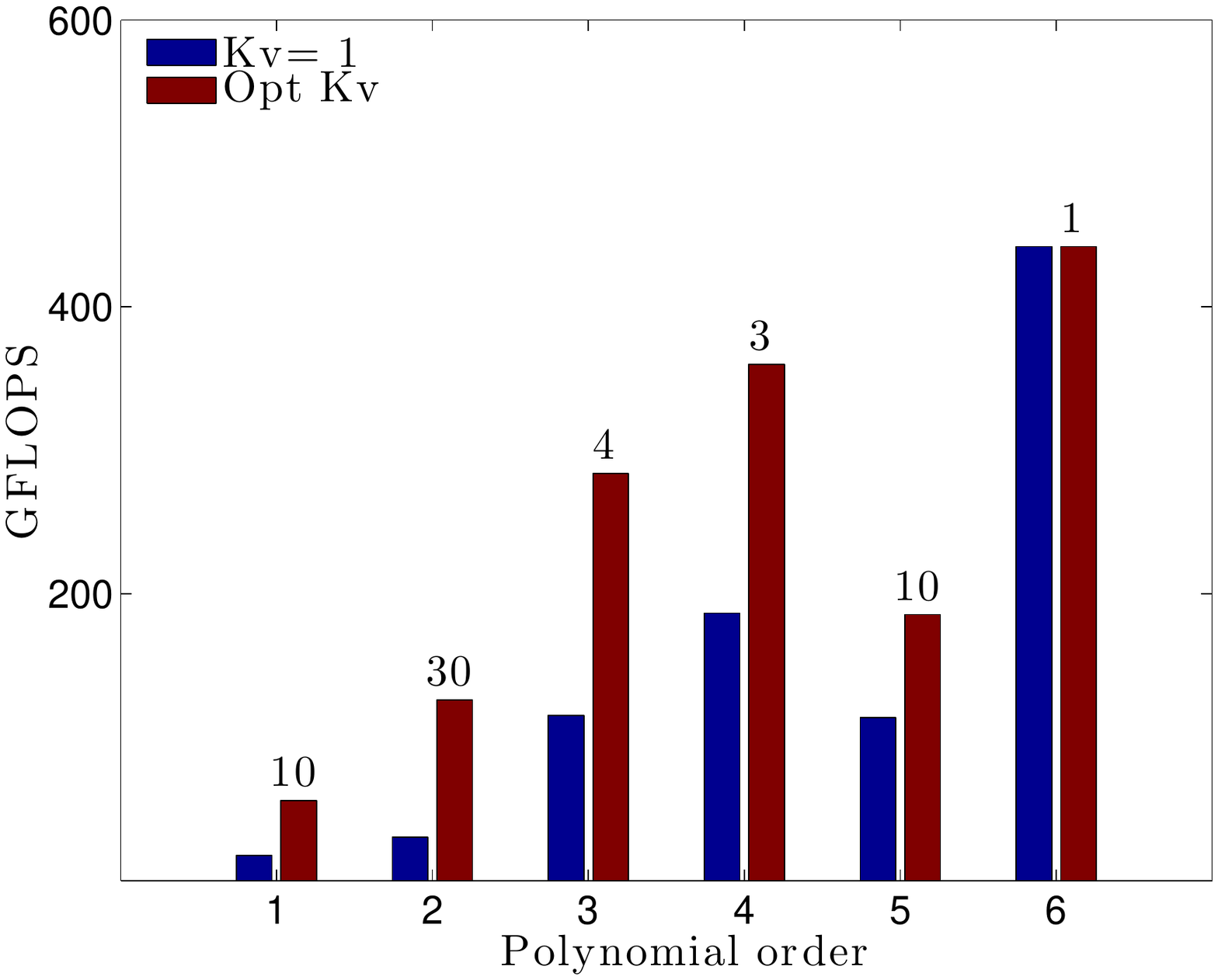}
  \end{minipage}} \hspace{0.5cm}
  \subfloat[Surface kernel on AMD Tahiti 7970]{%
    \begin{minipage}[c]{0.45\linewidth}
      \centering%
      \includegraphics[trim=1cm 6cm 1cm 6cm,clip=true,width=1\textwidth]{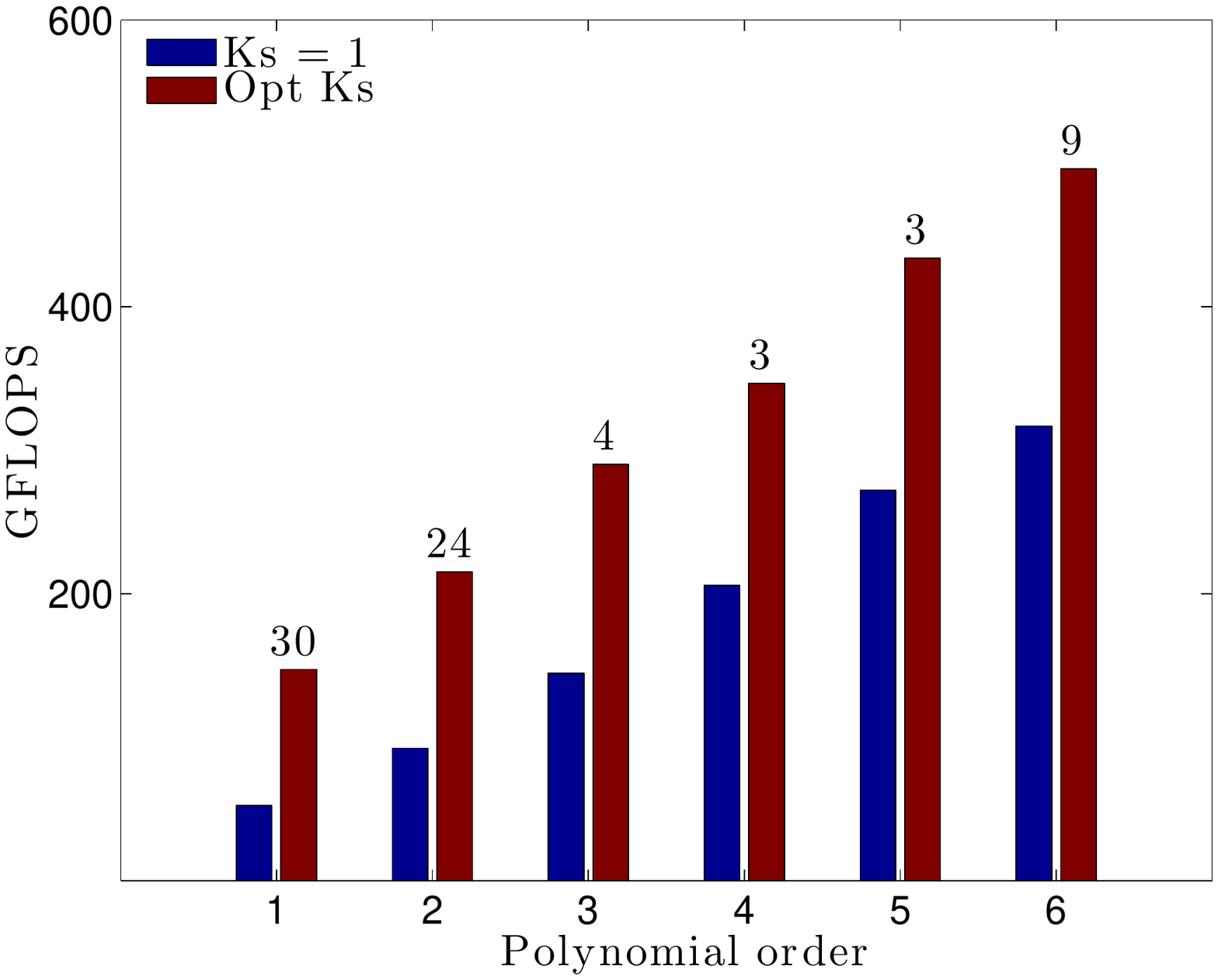}
  \end{minipage}}
  \caption{GFLOPS of native OpenCL volume kernel and surface kernel vs polynomial order. }
  \label{fig:flops_vs_numElements}
\end{center}
\end{figure}

\begin{figure}
 \begin{center}
  \subfloat[Volume kernel on NVIDIA Tesla C2050]{
  \begin{minipage}[c]{0.45\linewidth}
      \centering%
      \includegraphics[trim=1cm 6cm 1cm 6cm,clip=true,width=1\textwidth]{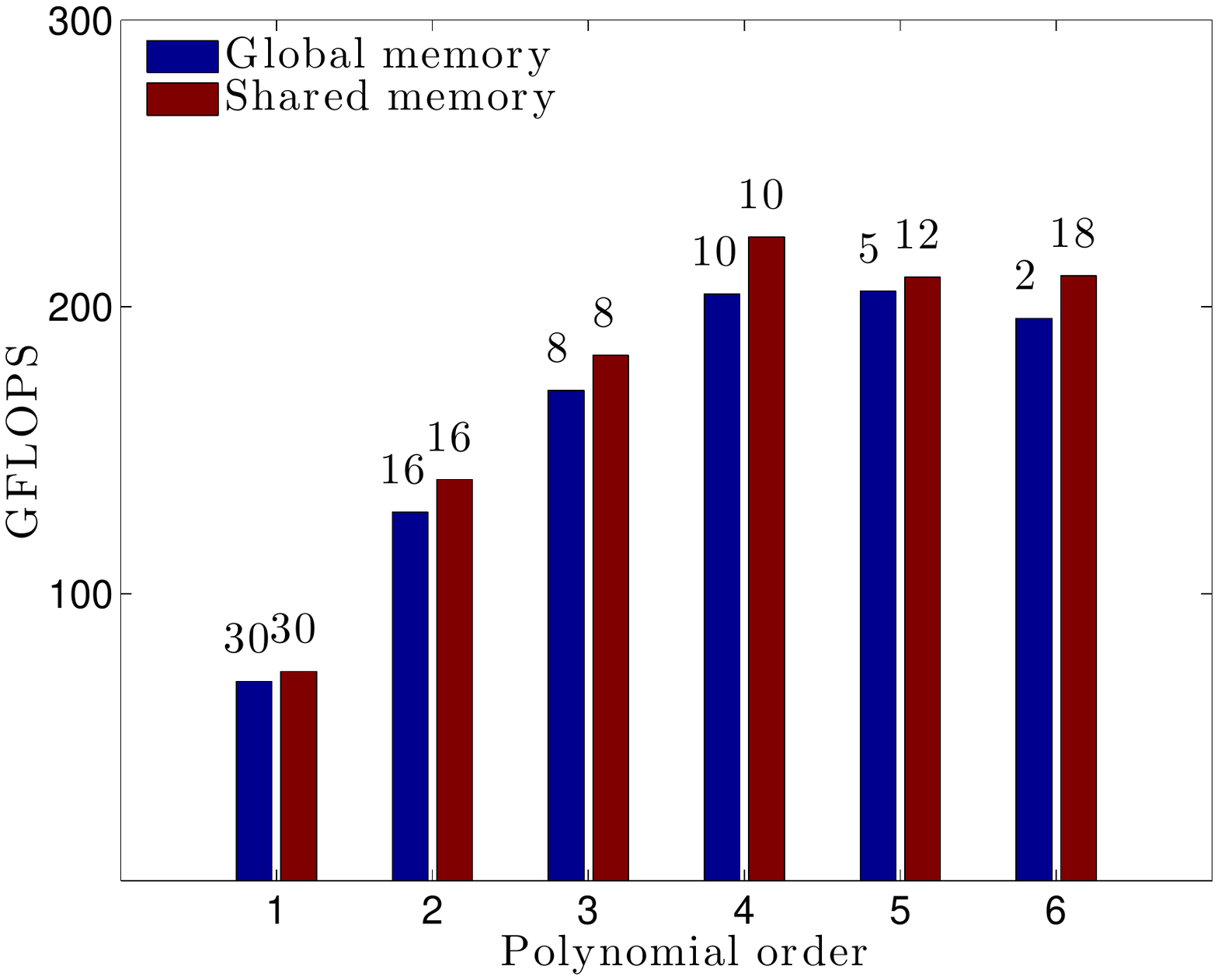}
  \end{minipage}}
  \caption{GFLOPS of native OpenCL volume kernel vs Polynomial order, when the operators are not stored in shared memory (blue) and are stored in shared memory (red).}
  \label{fig:global_vs_shared}
 \end{center}
\end{figure}

\begin{itemize}
\item \textbf{Coalescing:}  The nodal values corresponding to each element, projection operator and interpolation operator are  accessed contiguously from memory to maximize memory bus utilization.
\item \textbf{Padding:} The vector of nodal values for each element is padded with a factor of $4$ to make sure array accesses are aligned in order to minimize bank conflicts. \item \textbf{Unrolling:} Loops are unrolled to reduce the number of instructions to be executed, leading to hiding the latencies in reading the data from  memory.












\item \textbf{Multiple elements per work group:} Multiple elements are processed by each work group to better align the computation with hardware architecture. For low order approximation the speed up is about 5 times compared to using one element processed by a work group. Because of the limited availability of  shared memory and the difference in computation patterns, the optimal number of elements per work group varies for each kernel. The optimal parameters depend on the hardware of the GPU to a large extent. For example, we observe a significant difference in the performance on Tesla C2050 (NVIDIA) compared to Radeon 7970 (AMD). We can see from Fig. (\ref{fig:flops_vs_numElements}) that the performance improvement by using multiple elements processed by each work group is much more significant for Tesla C2050 compared to that of Tahiti 7970.
\item \textbf{Shared/local memory:} The nodal values on each element and/or projection and interpolations  operators  are stored in shared memory in order to reuse the data efficiently within a kernel. Since shared memory is limited in GPUs, using a large shared memory in a kernel will reduce the number of concurrent active work groups, but since the operators do not vary across the elements, they have to be stored only once for all the elements that are processed by a work group. Fig. (\ref{fig:global_vs_shared}) shows the improvement in performance when using shared memory over global memory access for the operators in volume kernel. \end{itemize}

\subsection{OCCA : CUDA vs OpenCL vs OpenMP}
\begin{figure}
\begin{center}
  \subfloat[Volume kernel]{%
    \begin{minipage}[c]{0.45\linewidth}
      \centering%
      \includegraphics[trim=1cm 6cm 1cm 6cm,clip=true,width=1\textwidth]{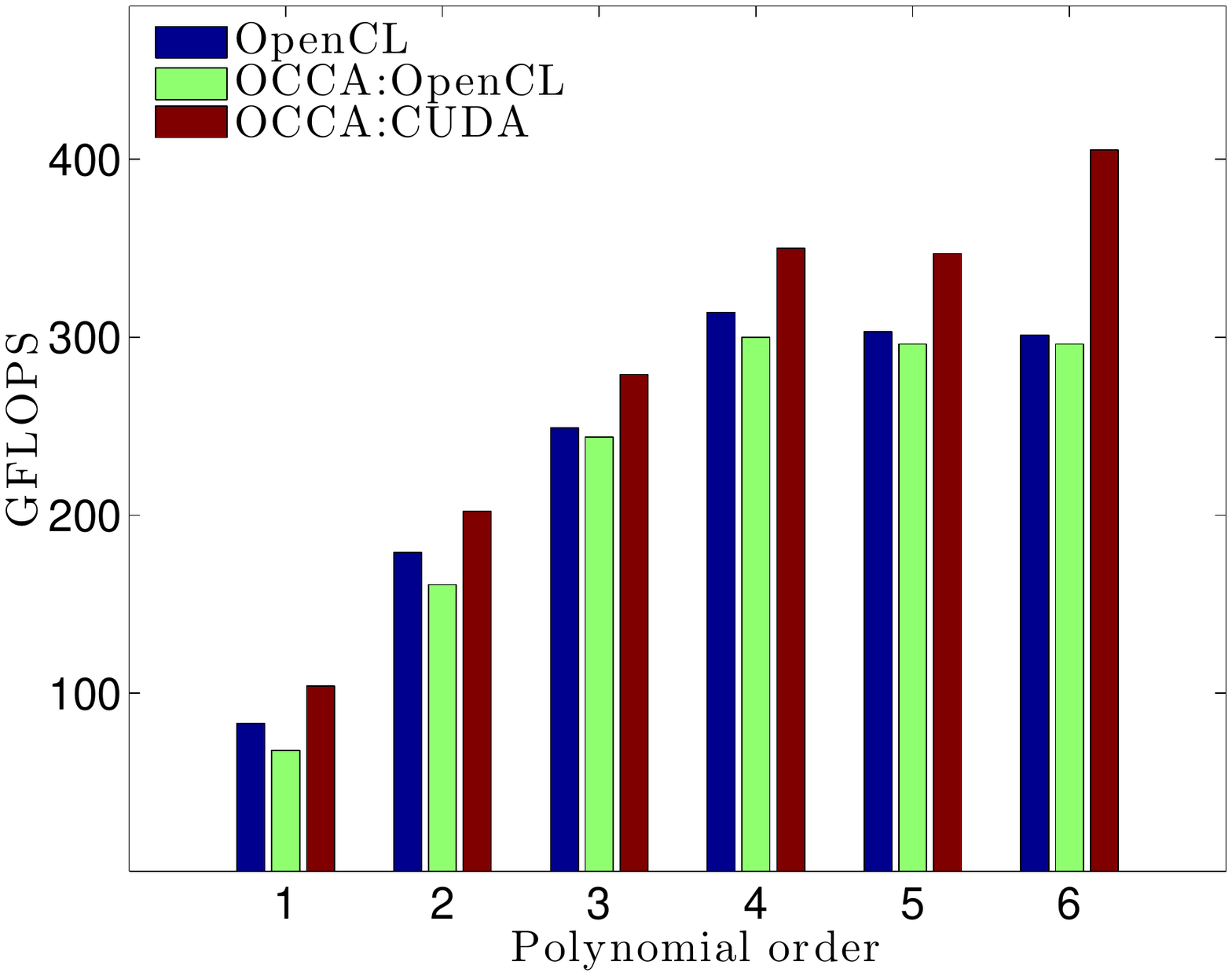}
  \end{minipage}} \hspace{0.5cm}
  \subfloat[Surface kernel]{%
    \begin{minipage}[c]{0.45\linewidth}
      \centering%
      \includegraphics[trim=1cm 6cm 1cm 6cm,clip=true,width=1\textwidth]{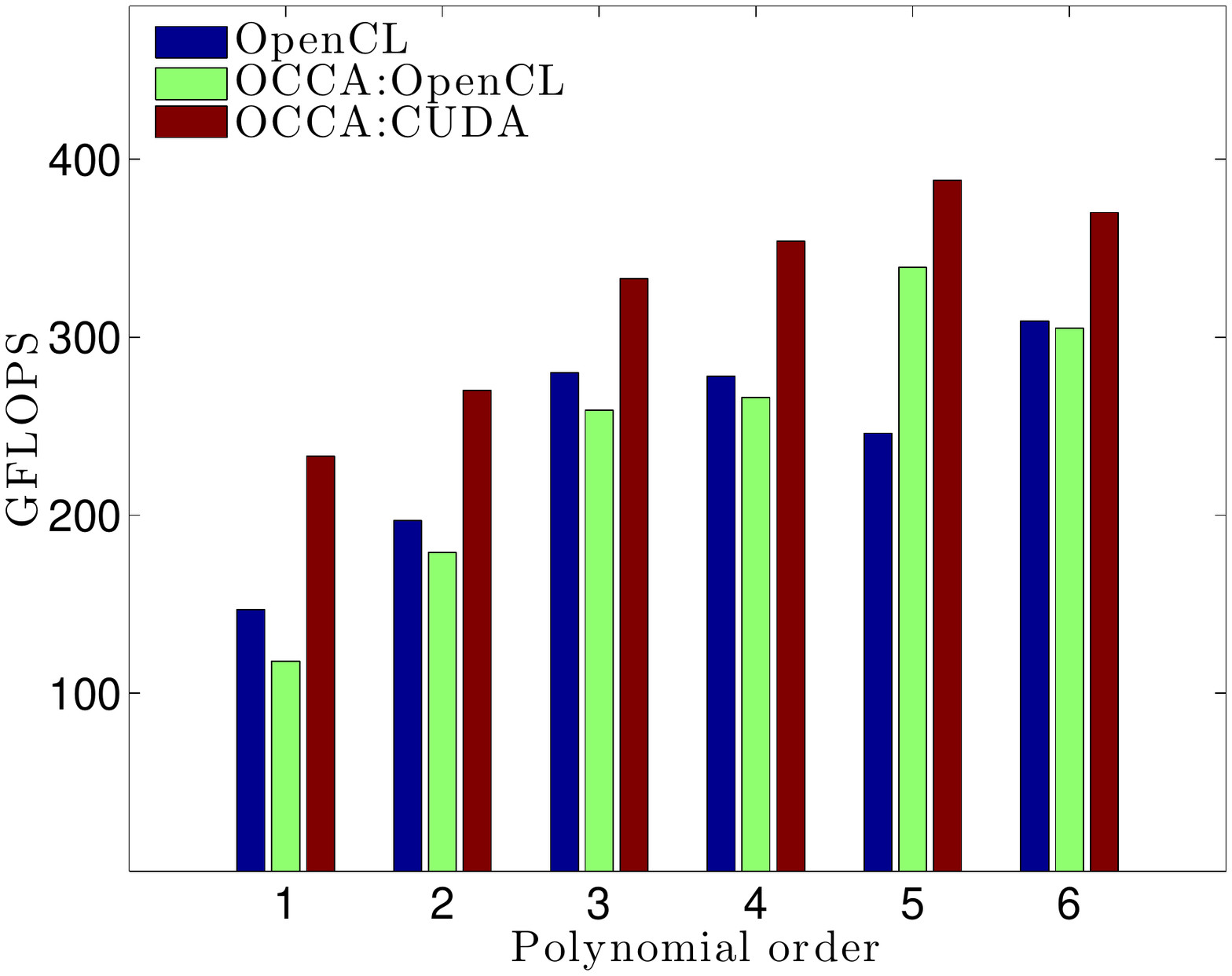}
  \end{minipage}}
  \caption{GFLOPS of  kernels vs polynomial order. A comparison of OpenCL, OCCA:OpenCL, and OCCA:CUDA. Experiments ran on NVIDIA Titan GPU.}
  \label{fig:cl_vs_occaCL_vs_occaCUDA}
\end{center}
\end{figure}

As mentioned earlier, the kernels are written in OCCA, a unified approach to multi-threading languages, to cross compile the device-based work in the OpenCL, CUDA, and OpenMP threading models. In Fig. (\ref{fig:cl_vs_occaCL_vs_occaCUDA}), we compare the performance of the two main kernels, the volume and surface kernels, for each platform together with the original hand-coded OpenCL kernels.
The comparisons are done on an NVIDIA Titan GPU with the optimal kernel tunings for each model as shown in Fig. (\ref{fig:flops_vs_numElements}). Because CUDA ptx compilers can be optimized to it's hardware, the results between CUDA and OpenCL on the Titan are not surprising. However, the main focus is that we observe a similar performance between the OCCA kernels compared to our original native OpenCL kernels.

\begin{figure}
\begin{center}
  \subfloat[Volume kernel]{%
    \begin{minipage}[c]{0.45\linewidth}
      \centering%
      \includegraphics[trim=1cm 6cm 1cm 6cm,clip=true,width=1\textwidth]{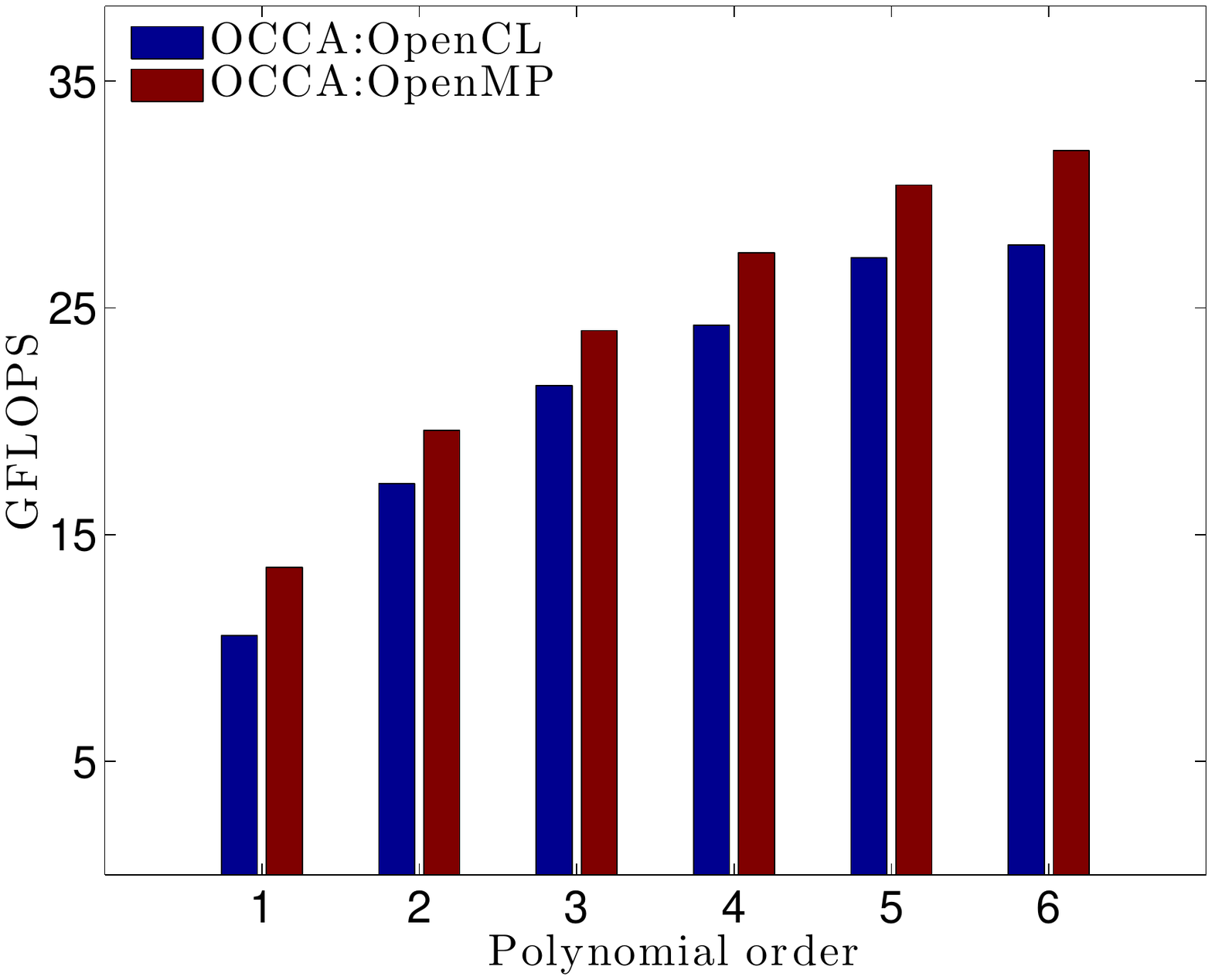}
  \end{minipage}} \hspace{0.5cm}
  \subfloat[Surface kernel]{%
    \begin{minipage}[c]{0.45\linewidth}
      \centering%
      \includegraphics[trim=1cm 6cm 1cm 6cm,clip=true,width=1\textwidth]{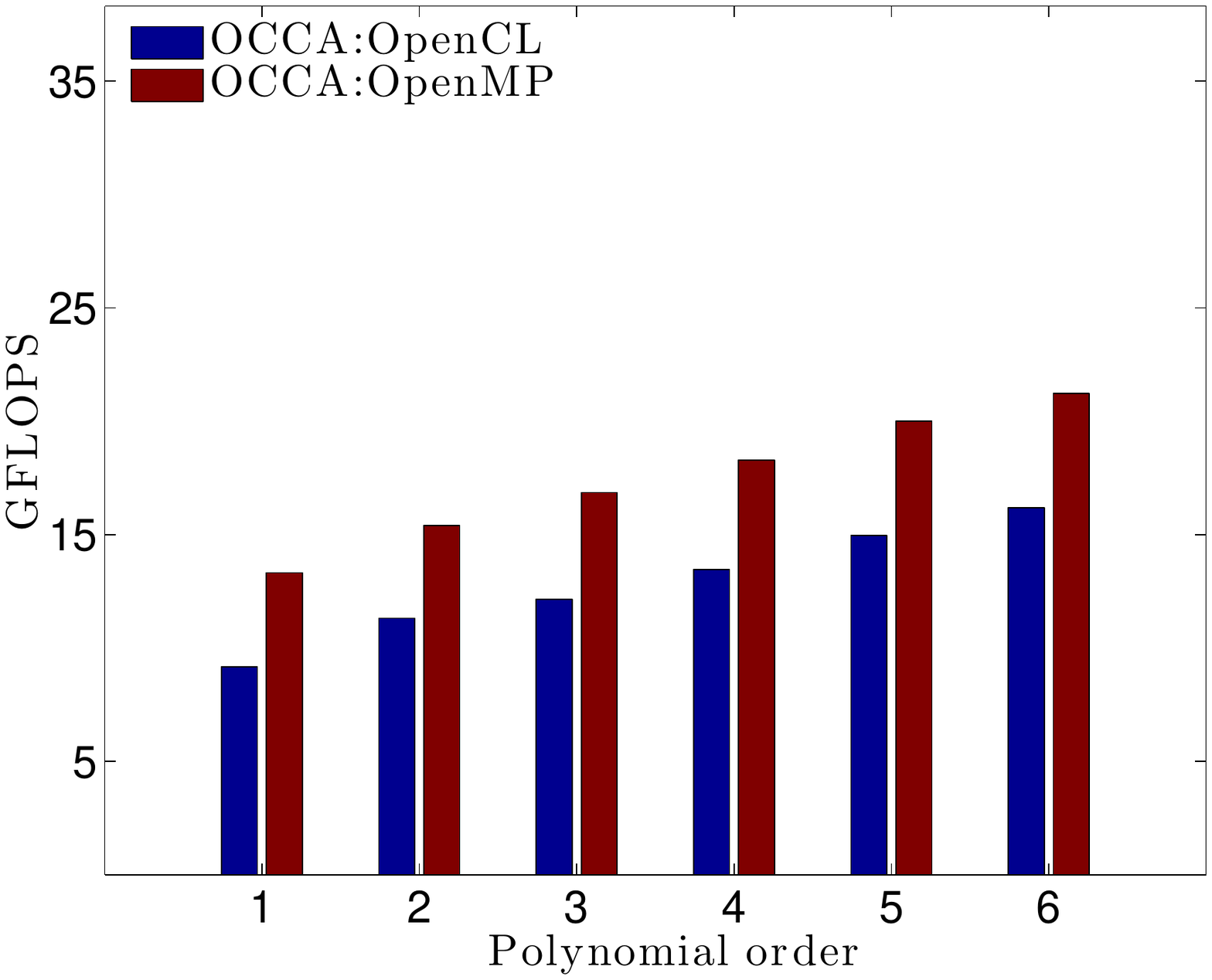}
  \end{minipage}}
  \caption{GFLOPS of  kernels vs polynomial order. A comparison
of OCCA:OpenCL and OCCA:OpenMP. Experiments ran on Intel(R) Core(TM) i7-3930K CPU .}
  \label{fig:occaCL_vs_occaOMP}
\end{center}
\end{figure}

In Fig. (\ref{fig:occaCL_vs_occaOMP}), we compare the performance of OCCA kernels when they are cross compiled with OpenCL and OpenMP on CPU. We observe that OCCA with OpenMP outperforms OCCA with OpenCL.


\section{Conclusions and future work}
A GPU accelerated discontinuous Galerkin algorithm for shallow water equations is presented. This algorithm is further accelerated with a multi-rate time stepping scheme. We presented a modified TVB limiter and a positivity preserving limiter for high order approximations. These algorithms are tested for accuracy, robustness and efficiency using several standard test cases in the literature. We compared the performance of OCCA\ kernels when they are cross compiled with threading models OpenCL, CUDA, and OpenMP at runtime.  
Our future work will include developing full three dimensional incompressible Navier-Stokes models using high order discontinuous Galerkin methods for better understanding of the tsunamis near coastal regions and accelerating such simulations through GPU hardware.


\section*{Acknowledgements}
The authors gratefully acknowledge travel grants from Pan-American Advanced Studies Institute. The authors also acknowledge the grant from DOE and ANL (ANL Subcontract No. 1F-32301 on DOE grant No. DE-AC02-06CH11357), and support from Shell (Shell Agreement No. PT22584), NVIDIA, and AMD.

\clearpage
\bibliographystyle{siam}
\bibliography{refs}

\begin{thebibliography}{10}

\bibitem{aizinger2002discontinuous}
{\sc V.~Aizinger and C.~Dawson}, {\em A discontinuous {G}alerkin method for
  two-dimensional flow and transport in shallow water}, Advances in Water
  Resources, 25 (2002), pp.~67--84.

\bibitem{berger2011geoclaw}
{\sc M.~Berger, D.~George, R.~LeVeque, and K.~Mandli}, {\em The {G}eo{C}law
  software for depth-averaged flows with adaptive refinement}, Advances in
  Water Resources, 34 (2011), pp.~1195--1206.

\bibitem{bunya2009wetting}
{\sc S.~Bunya, E.~Kubatko, J.~Westerink, and C.~Dawson}, {\em A wetting and
  drying treatment for the {R}unge--{K}utta discontinuous {G}alerkin solution
  to the shallow water equations}, Computer Methods in Applied Mechanics and
  Engineering, 198 (2009), pp.~1548--1562.

\bibitem{casulli1990semi}
{\sc V.~Casulli}, {\em Semi-implicit finite difference methods for the
  two-dimensional shallow water equations}, Journal of Computational Physics,
  86 (1990), pp.~56--74.

\bibitem{cockburn1998runge}
{\sc B.~Cockburn and C.-W. Shu}, {\em The {R}unge--{K}utta discontinuous
  {G}alerkin method for conservation laws {V}: multidimensional systems},
  Journal of Computational Physics, 141 (1998), pp.~199--224.

\bibitem{cools1999monomial}
{\sc R.~Cools}, {\em Monomial cubature rules since "{S}troud": a
  compilation-part 2},  (1999).

\bibitem{ern2008well}
{\sc A.~Ern, S.~Piperno, and K.~Djadel}, {\em A well-balanced {R}unge--{K}utta
  discontinuous {G}alerkin method for the shallow-water equations with flooding
  and drying}, International journal for numerical methods in fluids, 58
  (2008), pp.~1--25.

\bibitem{eskilsson2004triangular}
{\sc C.~Eskilsson and S.~Sherwin}, {\em {A} triangular spectral/hp
  discontinuous {G}alerkin method for modelling 2{D} shallow water equations},
  International Journal for Numerical Methods in Fluids, 45 (2004),
  pp.~605--623.

\bibitem{gallardo2007well}
{\sc J.~M. Gallardo, C.~Par{\'e}s, and M.~Castro}, {\em On a well-balanced
  high-order finite volume scheme for shallow water equations with topography
  and dry areas}, Journal of Computational Physics, 227 (2007), pp.~574--601.

\bibitem{gear1984multirate}
{\sc C.~W. Gear and D.~Wells}, {\em Multirate linear multistep methods}, BIT
  Numerical Mathematics, 24 (1984), pp.~484--502.

\bibitem{giraldo2008high}
{\sc F.~Giraldo and T.~Warburton}, {\em A high-order triangular discontinuous
  {G}alerkin oceanic shallow water model}, International journal for numerical
  methods in fluids, 56 (2008), pp.~899--925.

\bibitem{godel2010gpu}
{\sc N.~Godel, S.~Schomann, T.~Warburton, and M.~Clemens}, {\em {GPU}
  accelerated {A}dams--{B}ashforth multirate discontinuous {G}alerkin {FEM}
  simulation of high-frequency electromagnetic fields}, Magnetics, IEEE
  Transactions on, 46 (2010), pp.~2735--2738.

\bibitem{hesthaven2008nodal}
{\sc J.~Hesthaven and T.~Warburton}, {\em {N}odal discontinuous {G}alerkin
  methods: algorithms, analysis, and applications}, vol.~54, Springer Verlag,
  2008.

\bibitem{johnson1986analysis}
{\sc C.~Johnson and J.~Pitk{\"a}ranta}, {\em An analysis of the discontinuous
  {G}alerkin method for a scalar hyperbolic equation.}, Math. Comput., 46
  (1986), pp.~1--26.

\bibitem{klockner2010high}
{\sc A.~Klockner}, {\em High-performance high-order simulation of wave and
  plasma phenomena}, PhD thesis, Brown University, 2010.

\bibitem{klockner2009nodal}
{\sc A.~Klockner, T.~Warburton, J.~Bridge, and J.~Hesthaven}, {\em {N}odal
  discontinuous {G}alerkin methods on graphics processors}, Journal of
  Computational Physics, 228 (2009), pp.~7863 -- 7882.

\bibitem{kubatko2006hp}
{\sc E.~Kubatko, J.~Westerink, and C.~Dawson}, {\em hp {D}iscontinuous
  {G}alerkin methods for advection dominated problems in shallow water flow},
  Computer Methods in Applied Mechanics and Engineering, 196 (2006),
  pp.~437--451.

\bibitem{leveque2011tsunami}
{\sc R.~LeVeque, D.~George, and M.~Berger}, {\em Tsunami modelling with
  adaptively refined finite volume methods}, Acta Numerica, 20 (2011),
  pp.~211--289.

\bibitem{merrill1997finite}
{\sc D.~Merrill}, {\em Finite difference and pseudospectral methods applied to
  the shallow water equations in spherical coordinates}, PhD thesis, University
  of Colorado, 1997.

\bibitem{navon1988review}
{\sc I.~Navon}, {\em A review of finite-element methods for solving the
  shallow-water equations}, B. Schrefl er and OC Zienkiewicz, ed., Computer
  Modelling in Ocean Engineering, Balkeman, Rotterdam,  (1988), pp.~273--278.

\bibitem{shu1987tvb}
{\sc C.-W. Shu}, {\em {TVB} uniformly high-order schemes for conservation
  laws}, Mathematics of Computation, 49 (1987), pp.~105--121.

\bibitem{thacker1981some}
{\sc W.~C. Thacker}, {\em Some exact solutions to the nonlinear shallow-water
  wave equations}, Journal of Fluid Mechanics, 107 (1981), pp.~499--508.

\bibitem{wen2011gpu}
{\sc W.~H. Wen-mei}, {\em {GPU} {C}omputing {GEM}s {J}ade {E}dition}, Morgan
  Kaufmann, 2011.

\bibitem{xing2006high}
{\sc Y.~Xing and C.~Shu}, {\em High order well-balanced finite volume {WENO}
  schemes and discontinuous {G}alerkin methods for a class of hyperbolic
  systems with source terms}, Journal of Computational Physics, 214 (2006),
  pp.~567--598.

\bibitem{xing2010positivity}
{\sc Y.~Xing, X.~Zhang, and C.~Shu}, {\em Positivity-preserving high order
  well-balanced discontinuous {G}alerkin methods for the shallow water
  equations}, Advances in Water Resources, 33 (2010), pp.~1476--1493.

\bibitem{zhang2010positivity}
{\sc X.~Zhang and C.~Shu}, {\em On positivity-preserving high order
  discontinuous {G}alerkin schemes for compressible {E}uler equations on
  rectangular meshes}, Journal of Computational Physics, 229 (2010),
  pp.~8918--8934.

\end{thebibliography}
\end{document}